%
%

\magnification=1100

\font\titfont=cmr10 at 12 pt

\font\headfont=cmr10 at 12 pt



\overfullrule=0in

\def\boxit#1{\hbox{\vrule
 \vtop{%
  \vbox{\hrule\kern 2pt %
     \hbox{\kern 2pt #1\kern 2pt}}%
   \kern 2pt \hrule }%
  \vrule}}

  \def\harr#1#2{\ \smash{\mathop{\hbox to .3in{\rightarrowfill}}\limits^{\scriptstyle#1}_{\scriptstyle#2}}\ }

 \def\GG{{{\bf G} \!\!\!\! {\rm l}}\ }

\def\GL{{\rm GL}}

\def\uu{{\underline u}}
\def\ou{{\overline u}}
\def\seq{subequation }

\def\bbJt{{\bf J}^2}

\def\J{J}  
\def\bbf{{\bf F}}
\def\bbj{{\bf J}}
    
\def\jt{\J^2}
\def\jtx{\jt_x}

\def\ss{\subset}

\def\half{\hbox{${1\over 2}$}}
\def\smfrac#1#2{\hbox{${#1\over #2}$}}
\def\oa#1{\overrightarrow #1}

\def\dist{{\rm dist}}

\def\Hess{{\rm Hess}}

\def\tr{{\rm tr}}
\def\max{{\rm max}}
\def\min{{\rm min}}

\def\Hom{{\rm Hom\,}}
\def\det{{\rm det}}

\def\Sym{{\rm Sym}^2}

\def\arr{\longrightarrow}

\def\rn{\bbr^n}

\def\Int{{\rm Int}}

\def\Symn{{\Sym(\rn)}}

 \def\cm{{\cal M}}

\def\Theorem#1{\medskip\noindent {\bf THEOREM \bf #1.}}
\def\Prop#1{\medskip\noindent {\bf Proposition #1.}}
\def\Cor#1{\medskip\noindent {\bf Corollary #1.}}
\def\Lemma#1{\medskip\noindent {\bf Lemma #1.}}
\def\Remark#1{\medskip\noindent {\bf Remark #1.}}

\def\Def#1{\medskip\noindent {\bf Definition #1.}}

\def\Ex#1{\medskip\noindent {\bf Example \bf    #1.}}

\def\pf{\medskip\noindent {\bf Proof.}\ }
\def\qed{\hfill  $\vrule width5pt height5pt depth0pt$}
\def\n{\nabla}

   \def\cp{{\cal P}}

   \def\cm{{\cal M}}
   \def\cn{{\cal N}}

\def\cp{{\cal P}}
\def\cf{{\cal F}}

\def\vf{\varphi}

\def\wt{\widetilde}

\def\and{\qquad {\rm and} \qquad}
\def\arr{\longrightarrow}

\def\bbr{{\bf R}}\def\bbh{{\bf H}}\def\bbo{{\bf O}}
\def\bbc{{\bf C}}

\def\bbp{{\bf P}}

\def\bbf{{\bf F}}
\def\bbs{{\bf S}}

\def\bbs{{\bf S}}

\def\a{\alpha}
\def\b{\beta}
\def\d{\delta}
\def\e{\epsilon}
\def\f{\phi}
\def\g{\gamma}

\def\l{\lambda}
\def\o{\omega}

\def\s{\sigma}
\def\x{\xi}

\def\D{\Delta}
\def\L{\Lambda}
\def\G{\Gamma}
\def\O{\Omega}

\def\bo{\partial \Omega}

\def\Symn{\Sym(\rn)}
 
\def\USC{{\rm USC}}
\def\fa{{\rm\ \  for\ all\ }}

\def\cpt{\wt{\cp}}
\def\ft{\wt F}
\def\ob{\overline{\O}}

\def\Fa{{\oa F}}

 \def\LAG{{\rm LAG}}

\centerline{\titfont DIRICHLET DUALITY AND }
\medskip

\centerline{\titfont  THE NONLINEAR DIRICHLET PROBLEM}
\medskip

\centerline{\titfont   ON RIEMANNIAN MANIFOLDS}
\bigskip

\centerline{\titfont F. Reese Harvey and H. Blaine Lawson, Jr.$^*$}
\vglue .9cm
\smallbreak\footnote{}{ $ {} \sp{ *}{\rm Partially}$  supported by
the N.S.F. }

\vskip .5in
\centerline{\bf ABSTRACT} \bigskip
  \font\abstractfont=cmr10 at 10 pt
{{\parindent= .53in
\narrower\abstractfont \noindent

In this paper we study the Dirichlet problem for fully nonlinear second-order
equations on a riemannian manifold. As in our previous paper  [HL$_4$] we define 
equations via closed subsets of the 2-jet bundle where each equation 
has a natural dual equation.  Basic existence and
uniqueness theorems are established in a wide variety  of settings. 
However, the emphasis is on starting with a constant coefficient equation
as a model, which then universally determines an equation on every 
riemannian manifold which is equipped with a topological reduction 
of the structure group to the invariance group of the model. 
 For example, this covers  {\sl all branches} of the homogeneous complex 
Monge-Amp\`ere equation on an {\sl almost complex} hermitian manifold $X$.

In general, for an equation  $F$ on a manifold   $X$  and a smooth domain  $\O\ss\ss X$,
 we give geometric conditions  which imply  that the Dirichlet problem on  $\O$  
is uniquely solvable for all  continuous  boundary functions. We begin by 
introducing a weakened form of  comparison which has the advantage 
that local implies global.  We then introduce two fundamental concepts.
The first is the notion of a monotonicity cone $M$ for $F$.
If  $X$   carries a global  $M$-subharmonic
function, then weak comparison implies full comparison.  
The second notion is that of boundary $F$-convexity, which is defined
in terms of the asymptotics of $F$ and is used to define barriers.
In combining these notions the   Dirichlet problem 
becomes uniquely solvable
as claimed.

This article also introduces   the notion of {\sl local affine
jet-equivalence}  for subequations.  It is used in treating the cases
above, but gives results for a much broader spectrum of equations on manifolds, 
including inhomogeneous equations and the Calabi-Yau equation 
on almost complex hermitian manifolds.

A considerable portion of the paper is concerned with specific examples.
They include a wide variety of equations which make sense on any riemannian
manifold, and many which hold universally on almost complex or quaternionic 
hermitian manifolds, or topologically calibrated manifolds.

}}

\vfill\eject

\centerline{\bf TABLE OF CONTENTS} \bigskip

{{\parindent= .1in\narrower\abstractfont \noindent

\qquad 1. Introduction.\smallskip

\qquad 2.      $F$-Subharmonic  Functions.
\smallskip

\qquad 3.   Dirichlet Duality and the Notion of a Subequation. \smallskip

\qquad 4.    The Riemannian Hessian -- A Canonical Splitting of $J^2(X)$.
  
\smallskip

\qquad 5.  Universal  Subequations on Manifolds with Topological $G$-Structure.

\smallskip

\qquad 6.   Jet-Equivalence of  Subequations.
\smallskip

\qquad 7. Strictly $F$-Subharmonic Functions.
   
\smallskip

\qquad 8.  Comparison Theory -- Local to Global.
   
\smallskip

\qquad 9.  Strict Approximation and  Monotonicity Subequations.
   
\smallskip

\qquad 10.  A Comparison Theorem for $G$-Universal Subequations. 
   
\smallskip

\qquad 11. Strictly $F$-Convex Boundaries and Barriers.

\smallskip

\qquad 12.  The Dirichlet Problem -- Existence.

\smallskip

\qquad 13.  The Dirichlet Problem -- Summary Results.

\smallskip

\qquad 14.  Universal Riemannian Subequations.

\smallskip

\qquad 15.  The Complex and Quaternionic Hessians.

\smallskip

\qquad 16.  Geometrically Defined Subequations -- Examples.

\smallskip

\qquad 17.  Equations Involving The Principal Curvatures of the Graph.

\smallskip

\qquad 18. Applications of  Jet-Equivalence -- Inhomogeneous Equations.

\smallskip

\qquad 19. Equations of Calabi-Yau type in the Almost Complex Case.

}}

\vskip .3in

{{\parindent= .3in\narrower

\noindent
{\bf Appendices: }\medskip

A. \   Equivalent Definitions of $F$-Subharmonic.
\smallskip

B. \  Elementary Properties of $F$-Subharmonic Functions.
\smallskip

C. \  The Theorem on Sums.
\smallskip

D. \  Some Important Counterexamples.
\smallskip

}}

\vfill\eject


\centerline{\headfont \ 1.\  Introduction}
\medskip

In a recent article [HL$_4$] the authors studied the Dirichlet problem for fully nonlinear
 equations of the form
 $
 {\bf F}(\Hess\, u) = 0
 $
 on smoothly bounded domains in $\rn$.  Our approach employed a duality
 which enabled  us to geometrically characterize   domains for which one has  existence 
 and uniqueness of solutions    for all continuous boundary data.
 These results covered, for example, all branches of the real, complex and quaternionic Monge-Amp\`ere
 equations, and all branches of the special Lagrangian potential equation.
 
Here we shall   extend these results in several ways.   First, all results in [HL$_4$] are shown
 to carry over to riemannian manifolds with an appropriate {\sl topological} reduction of the 
 structure group. For example, we treat the complex Monge-Amp\`ere equation
 on almost complex manifolds with hermitian metric.
 Second, the results are  extended to equations involving the full 2-jet of functions.
 There still remains a basic notion of duality, and a geometric form of  {$F$-boundary convexity}.  
Existence  and uniqueness theorems are established, and a large number of examples are examined in
 detail.

 In [HL$_4$]   our  approach was to replace  the function ${\bf F}(\Hess\, u)$
   by  a closed subset $F\ss \Symn$   of the symmetric  $n\times n$ matrices 
 subject only to the {\sl positivity condition }
$$
F+\cp\ \ss\ F
\eqno{(1.1)}
$$
 where $\cp \equiv \{A\in \Symn: A\geq 0\}$  (cf. [K]).  
 Such an $F$ was called a  {\sl Dirichlet set} in [HL$_4$]  but will be   called a {\sl subequation} here.
 A   $C^2$-function $u$ on a domain $\O$ is 
 {\sl $F$-subharmonic}  if $ \Hess_x u\in   F$ for all $x\in\O$, and it is  {\sl $F$-harmonic} 
  if $ \Hess_x u \in\partial F$ for all $x$. (The reader might note the usefulness
  of this approach in treating other branches of the equation $\det(\Hess\, u)=0$.)
  The important point is to extend these notions
  to upper semi-continuous $[-\infty,\infty)$-valued functions.
We will follow the approach in   [HL$_4$]   by using the {\sl dual subequation}
$$
\wt F \ \equiv \ -(\sim F)\ =\ \sim(-F).
\eqno{(1.2)}
$$
The class of subaffine functions, or $\cpt$-subharmonic functions, played a key role in [HL$_4$]
 since the positivity condition $F+\cp   \ss F$ is equivalent
to $F+\ft\ss \cpt$.

Each  subequation  $F$ has an associated  {\sl asymptotic interior} $\oa F$, and
using this we introduced a notion of strictly  {\sl $\Fa$-convex} boundaries.
For the $q$th branch of the complex Monge-Amp\`ere equation for example, this is just
classical $q$-pseudo-convexity.    
It is then shown in [HL$_4$]  that for any domain $\O$ whose boundary 
is strictly $\Fa$- and $\overrightarrow{{\wt{F}}}$-convex,
solutions to the Dirichlet Problem exist for all continuous boundary data.  Furthermore, uniqueness
holds on arbitrary domains $\O$.

In this paper the  results from [HL$_4$] will be generalized in several ways.
To begin we shall work on a general manifold $X$ and consider subequations given by closed subsets
$$
F\ss J^2(X)
$$
 of the 2-jet bundle which satisfy three conditions.  The first is the {\sl positivity condition} 
 $$
 F+\cp \ \ss\ F
 \eqno{(P)}
 $$ 
where $\cp_x$ is the set of 2-jets of non-negative functions with critical value zero at x. The second is the  {\sl negativity condition} 
$$
F+\cn \ \ss\ F,
 \eqno{(N)}
 $$  
where $\cn = $ the jets of non-positive constant functions. The third is
 a mild {\sl topological condition}, which is  satisfied in all interesting cases.  
(See \S 2.) We always assume our subequations $F$ satisfy these conditions.
For the pure second-order constant coefficient subequations considered in [HL$_4$]
the positivity conditions (P) and (1.1) are equivalent and automatically imply the other two conditions (N) and (T),
as long as  $F$ is a closed set.

Here the dual subequation $\ft$ is  defined by (1.2) exactly as in  [HL$_4$].  While duality continues
to be important, it is not used to define 
$F$-subharmonic  functions
in the  upper semi-continuous  setting. 
This is because subaffine (i.e., co-convex) functions do
not satisfy the maximum principle on general riemannian manifolds.
 So we use instead the equivalent
[HL$_4$,  Remark 4.9] viscosity
definition (cf. [CIL]).  The $F$-subharmonic functions have most of the important properties
of classical subharmonic functions, such as closure under taking maxima,
 decreasing limits, uniform limits, and upper envelopes.
 This is discussed in Section 2 and Appendices A and B.
 
 The notion of a strictly  $\Fa$-convex boundary can be extended  to this general
 context. This is discussed below in the introduction and treated in Section 11.
  
 The focus of this paper is an important class of subequations   constructed on
 riemannian manifolds, with an appropriate reduction of  structure group,
 from constant coefficient equations on $\rn$.
 It is a basic fact that on riemannian manifolds $X$, there is a canonical bundle
 splitting 
 $$
 J^2(X) \ =\ \bbr\oplus T^*X\oplus \Sym(T^*X)
 \eqno{(1.3)}
 $$
given by the riemannian hessian.  This enables us for example to carry over
any O$_n$-invariant constant coefficient subequation $F\ss \Symn$
(the so-called {\sl Hessian equations} in [CNS])
to all riemannian manifolds.  That is, any purely second order subequation on $\rn$  which depends
only on the eigenvalues of the hessian carries over to general riemannian manifolds.

Much more generally however, consider a constant coefficient subequation, i.e., a  subset
$$
\bbf \ \ss\  \bbr\oplus  \rn\oplus \Symn 
$$ 
 satisfying (P), (N) and (T), and its compact invariance group
$$
G\ \equiv\ \{g\in  {\rm O}_n : g(\bbf)=\bbf\}.
$$ 
Suppose $X$ is a riemannian manifold with a {\sl topological  $G$-structure}. 
This means that  $X$ is provided with an  open covering
$\{U_\a\}_\a$  and  a  tangent frame field $e^\a=(e^\a_1,...,e^\a_n)$ on each
$U_\a$ so that each change of framing $G_{\a\b}:U_\a\cap U_\b \to G \ss {\rm O}_n$
is valued in the subgroup $G$.  Then from the splitting (1.3), the subequation  $\bbf$ 
can be transplanted to a globally defined subequation $F$ 
 on $X$.  Such subequations are called {\sl $G$-universal} or {\sl riemannian $G$-subequations}
  and make sense
on any riemannian manifold with a topological $G$-structure.

An important new ingredient  is
the concept of a convex {\sl monotonicity cone} for $F$ (Definition 9.4).  This is a subset $M\ss J^2(X)$
which is a convex cone (with vertex 0) at each point and satisfies
$$
F+M\ \ss\ F.
$$
The subequation $P$ defined by requiring that the riemannian hessian
be $\geq 0$  is a monotonicity  cone for all 
purely second-order subequations.  The corresponding $P$-subharmonic functions are just the
riemannian convex functions. (Note that these exist globally on $\rn$ and, in fact, on
any complete simply-connected manifold of non-negative sectional curvature.)
Similarly,  on an almost complex hermitian manifold $X$ the subequation
$P^\bbc$,  defined by requiring  the hermitian symmetric part
of $\Hess \,u$ to be $\geq 0$, is a  monotonicity cone for 
any U$_n$-universal subequation depending only on the hermitian symmetric part. 
The corresponding $P^\bbc$-subharmonic functions are 
{\sl hermitian  plurisubhamonic functions} on $X$.

Our first main result is the following.

\Theorem{13.1} {\sl Suppose $F$ is a riemannian $G$-subequation on a riemannian manifold
$X$ provided with a topological $G$-structure as above.
Suppose there exists a 
$C^2$ strictly $M$-subharmonic function on $X$ where $M$ is a monotonicity cone for $F$.

Then for every domain $\O\ss\ss X$ whose boundary is strictly
$F$- and $\ft$-convex, both existence and uniqueness
hold for the Dirichlet problem.
That is, for every    $\vf\in  C(\bo)$ there exists a unique $F$-harmonic
function $u\in C(\overline{\O})$ with  $u \bigr|_{\bo}=\vf$.      }
 
\medskip

The simplest case of this, where $F$ is a constant coefficient equation
in euclidean $\rn$ (and $G=\{1\}$), already generalizes
the main  result in [HL$_4$]. Here comparison holds for any subequation
which is gradient independent since the squared distance to a point
can be used to construct the required $M$-subharmonic function. (See Theorem 13.4.)
A similar comment holds for any complete simply-connected 
manifold $X$ of non-positive sectional curvature.

In general, requiring the existence of a strictly $M$-subharmonic function (or something
similar) is an intuitively necessary global hypothesis;   for example in the
case where $G= {\rm SO}_n$, one is free to arbitrarily change the riemannian
geometry inside a domain $\O$ while preserving the  $F$-convexity
of $\bo$. However, even for quite regular metrics -- domains in $S^3\times S^3$ --
examples in Appendix D show that uniqueness fails without this hypothesis.

Theorem 13.1 applies in a quite broad context.  We point out some examples here.
They will be discussed in detail in the latter sections of the paper.

\medskip\noindent
{\bf Example A. (Parallelizable manifolds).} Even when $F$ admits absolutely no symmetries, i.e. $G=\{1\}$,
the theorem has broad applicability.  Suppose $X$ is a riemannian manifold which
is parallelizable, that is, on which there exist global vector fields $e_1,...,e_n$ which
form a basis of $T_xX$ at every point.  
The single open set $U=X$ with $e=(e_1,...,e_n)$ is a $G=\{1\}$-structure, and 
{\sl every constant coefficient subequation $F$ in $\rn$ can be carried over to a subequation $F$ on $X$}.

Note that many manifolds are indeed parallelizable.  For example, all orientable 3-manifolds
have this property.

\medskip\noindent
{\bf Example B.  (General riemannian manifolds).} Any constant coefficient
subequation $\bbf  \ss \bbr\oplus \rn\oplus\Symn$ which is invariant under
the action of O$_n$ carries over to all riemannian manifolds.  There are, for example,
many invariant functions of the riemannian hessian which yield
universal  ``purely second-order''  equations.
 Associated to each $A\in\Symn$
is its set of ordered eigenvalues $\l_1 \leq \cdots\leq\l_n$.  These eigenvlaues have the
property that $\l_q(A+P)\geq \l_q(A)$ for any $P\in \cp$. Thus for example,
 the set $\bbp_q= \{A: \l_q(A)\geq 0\}$ gives a subequation on any riemannian manifold.
This subequation $P_q$ is the $q$th branch of the Monge-Amp\`ere equation.
 
 There are  many other equations of this type.  In fact, let $\L \ss \rn$ be any subset
 which is invariant under permutations of the coordinates and satisfies the   positivity condition:
 $\L+\bbr_+^n \ss \L$ where $\bbr_+\equiv \{t\geq 0\}$. Then this set, considered as a relation on the eigenvalues
 of the riemannian hessian, determines a universal subequation on every riemannian
 manifold. 
 
A  classical case  is $\bbf(\s_k) = \{ \l : \s_1(\l)\geq 0,..., \s_k(\l)\geq 0\}$ (the principal branch of $\s_k(\l)= 0$
where $\s_k$ denotes the $k$th elementary symmetric function). 
One can compute that  the  convexity of $\bo$, for  this equation and its dual,
corresponds to the  $\bbf(\s_{k-1}) $-convexity of its second fundamental form.
For domains in $\rn$ this result was proved in [CNS].
It is generalized  here in two ways.  We establish the result for {\sl all branches} of the equation
 and we carry it over 
to general riemannian manifolds. (See Theorem 14.4.)

 We are also able to treat inhomogeneous subequations such as $\Hess \, u \geq 0$ and $\det \Hess\, u
 \geq f(x)$ for a positive function $f$ on $X$ (and the  analogous subequations for the other    $\s_k$ as above).
 This follows from Theorem 10.1 on  local affine jet-equivalence.

\medskip\noindent
{\bf Example C.  (Almost complex hermitian manifolds).}  Consider $\bbc^n = (\bbr^{2n}, J)$ where 
$J:  \bbr^{2n}\to  \bbr^{2n}$ represents multiplication by $\sqrt{-1}$. To any $A\in \Sym( \bbr^{2n})$ 
one can associate the hermitian symmetric part $A^\bbc \equiv \half (A-JAJ)$ which has
$n$ real eigenvalues $\l_1\leq\cdots \leq\l_n$ occuring with multiplicity 2 on $n$ complex lines.
The entire discussion in Example B now applies.  Any permutation-invariant subset
$\L\ss\rn$ satisfying the $\bbr_+^n$-positivity condition $\L+\bbr_+^n \ss \L$ gives a natural subequation
on any manifold with $U_n$-structure, i.e., any almost complex manifold with a compatible riemannian
metric. This includes for example all branches of the homogeneous complex Monge-Amp\`ere equation.
We note that almost complex hermitian manifolds play an important role in modern symplectic geometry.

One can also treat the Dirichlet problem for Calabi-Yau-type equations in this context of almost complex
hermitian manifolds.  See Example 6.15 for example.

One can also consider U$_n$-invariant functions of the skew-hermitian part
$A^{\rm sk}=\half(A+JAJ)$. An important case of this is   the Lagrangian equation discussed below.

\medskip\noindent
{\bf Example D.  (Almost  quaternionic hermitian manifolds).}   A discussion parallel to that in Example C holds
with the complex numbers $\bbc$ replaced by the quaternions $\bbh = (\bbr^{4n}, I,J,K)$.  
In particular Theorem 13.1 applies to
the Dirichlet problem for  all branches of the quaternionic Monge-Amp\`ere equation on 
almost  quaternionic hermitian manifolds.

\medskip\noindent
{\bf Example E.  (Grassmann structures).} 
Fix any closed subset $\GG\ss G(p,\rn)$ of the Grassmannian of $p$-planes in $\rn$, and consider the
subequation $F(\GG)$ defined for $(r,p,A) \in \bbr\oplus\rn\oplus\Symn$ by the condition
$$
\tr_\xi A \ \equiv\ \tr\left( A\bigr|_\x \right)\ \geq\ 0 \fa \x\in\GG.
$$
Here the $F(\GG)$-subharmonic functions are more appropriately 
called $\GG$-{\sl plurisubhamonic}.
Now $F(\GG)$ carries over to a subequation on any riemannian manifold with $G$-structure
where 
$$
G \equiv \{g\in{\rm O}_n : g(\GG)=\GG\}.
$$

For example, if $\GG = G(p,\rn)$, then $G={\rm O}_n$ and the corresponding subequation,
which makes sense on any riemannian manifold, states that the sum of any $p$ eigenvalues 
of $\Hess\,u$ must be $\geq0$.  This is called {\sl geometric $p$-convexity} or {\sl $p$-plurisubhamonicity}.

A particularly  interesting example is given by the
U$_n$-invariant  set $\LAG\ss G(n,\bbc^{n})$ of Lagrangian subspaces
of $\bbc^n= (\bbr^{2n},J)$.  
On any almost complex hermitian manifold, this gives rise to a notion of {\sl Lagrangian subharmonic}
and {\sl Lagrangian harmonic} functions with an associated Dirichlet problem which is
solvable on Lagrangian convex domains.
This is discussed in more detail below.

Another rich set of examples comes from {\sl almost calibrated manifolds} such
as  manifolds with topological G$_2$ and Spin$_7$ structures. These are also discussed at the
end of the introduction.

For all such structures Theorem 13.1 gives the following general result.
Call a domain $\O\ss X$  {\bf strictly $\GG$-convex} if it  has a strictly
$\GG$-psh defining function. This holds if $\bo$ is strictly
$\GG$-convex (cf. (1.5) below) and there exists $f\in C^2(\ob)$ which is strictly $\GG$-psh.

\vfill\eject

\Theorem {16.1}  {\sl Let $X$ be a riemannnian manifold with topological 
$G$-structure so that the riemannian $G$-subequation $F(\GG)$ is defined on $X$.
Then on any strictly $\GG$-convex domain $\O\ss\ss X$
the Dirichlet problem for $F(\GG)$-harmonic functions is uniquely solvable
 for all continuous boundary data.}

\vskip .3in
\centerline{\headfont Jet-Equivalence of Subequations}\bigskip

Although it is not our  main focus,  one of the important results of this paper is that
Theorem 13.1 actually holds for a much broader class of subequations, namely those
which are locally affinely jet-equivalent to constant coefficient subequations.  
The notion of (ordinary) jet-equivalence  is expressed  in terms of  automorphisms of the 2-jet bundle.
An (linear) automorphism of $J^2(X)$ is a smooth bundle isomorphism $\Phi: J^2(X) \to J^2(X)$
which has certain natural properties with respect to the short exact sequence
$0\to \Sym(T^*X) \to J^2(X) \to J^1(X)\to 0$. (See Definition 6.1 for details).
The automorphisms form a group, and two subequations $F, F'\ss J^2(X)$ are
called {\sl jet-equivalent} if there exists an automorphism $\Phi$ with $\Phi(F)=F'$.

A subequation on a general manifold $X$ is said to be {\sl locally jet-equivalent to a constant coefficient
subequation}    if each point $x$ has a coordinate neighborhood $U$ so that 
$F\bigr|_U$ is jet-equivalent to a constant coefficient subequation $U\times \bbf$ in those
coordinates.

Any riemannian $G$-subequation on a riemannian manifold with topological $G$-structure
is locally jet-equivalent to a constant coefficient subequation.  Theorem 13.1 remains true 
if one  drops the riemannian metric on $X$ and 
replaces the `` $G$-universal assumption'' with the assumption that $F$ is locally jet-equivalent to a constant coefficient subequation.  This is a strictly broader class of equations.  For example,  a generic 
subequation defined by a smoothly varying linear inequality in the fibres of $J^2(X)$ 
is locally jet-equivalent to constant coefficients.

The notion of jet-equivalence can be further substantially broadened by using the {\sl affine} automorphism group.
This is the fibrewise extension of the linear automorphisms by the full group
of translations in the fibres $J^2_x(X)$.   

\Theorem {13.1$'$}  {\sl 
Theorem 13.1 holds for subequations
which are  locally affinely jet-equivalent to  riemannian $G$-subequations.}
\medskip
\noindent
  This allows one to 
treat inhomogeneous equations with variable right hand side.  For example one can establish 
existence and uniqueness  of solutions to the 
Dirichlet problem for the Calabi-Yau type equation
$$
\left( {1\over i}\partial\overline{\partial} u  + \o \right)^n\ =\ e^u f(x) \o^n.
\eqno{(1.4)}
$$
or the same equation with $e^u$ replaced by any non-decreasing positive  function $F(u)$.
Here the results hold on any {\bf almost complex}  hermitian manifold (see Section 19).
One also gets existence and uniqueness results for equations of the type
$$
\s_1(\Hess\, u)\ \geq \ 0,\ ...\ , \s_{k}(\Hess\, u)\ \geq\ 0  \qquad{\rm and}
$$
$$
\s_k(\Hess\, u)\ =\ f(x)
$$
where $f>0$.  This is also true for equations such as 
$$
\l_q(\Hess\, u) \ =\ f(x)
$$  for any smooth function $f$.  Examples are given in \S 18.

\smallskip
\noindent
{\bf Note.} The notion of jet-equivalence is a quite weak relation.
A jet-equivalence of subequations $\Phi:F\to F'$ does {\sl not} induce a correspondence between
 $F$-subharmonic 
functions and  $F'$-subharmonic functions.  In fact, 
 for a $C^2$-function $u$, $\Phi(J^2 u)$ is almost never
the 2-jet of a function.

\vskip .3in
\centerline{\headfont Comparison Results}\bigskip

This paper contains a number of other possibly quite useful existence and uniqueness
results which lead to Theorem 13.1 above.

A central concept in this subject is that of comparison for upper semi-continuous
functions, which is treated in Section 8.

\medskip
\noindent
{\bf Definition.}  We say that {\sl  comparison  holds} for $F$ on X if for all compact
subsets $K\ss X$ and functions 
$$
u\in F(K) \and  v\in \wt F(K), 
$$
   the {\bf Zero Maximum Principle} holds for $u+v$ on $K$, that is,
$$
u+v\ \leq \ 0\quad {\rm on\ \ } \partial K\qquad\Rightarrow \qquad
u+v\ \leq \ 0\quad {\rm on\ \ } K.
\eqno{(ZMP)}
$$

\medskip

If comparison holds for $F$ on $X$, then one easily deduces the uniqueness of solutions
to the Dirichlet problem for $F$-harmonic functions on every compact subdomain.

Obviously comparison for small compact sets $K$ does not imply that comparison holds
for arbitrary compact sets.  However, this is true for a weakened form of comparison involving
a notion of strictly $F$-subharmonic functions.

Consider a subequation $F$ on a riemannian manifold $X$. For each constant $c>0$ we 
define $F^c\ss F$ to be the subequation whose fibre at $x$ is
$$
F_x^c \ \equiv\ \{J\in F_x : \dist(J, \sim F)\geq c\}
$$
where ``$\dist$'' denotes distance in the fibre $J_x(X)$.  We define an  upper semi-continuous  function $u$ 
on $X$ is to be  {\sl strictly $F$-subharmonic} if for each point $x\in X$ there is a neighborhood
$B$ of $x$ and a $c>0$ such that $u$ is $F^c$-subharmonic on $B$.
\medskip
\noindent
{\bf Definition.}  We say that {\sl weak comparison  holds} for $F$ on X if for all compact
subsets $K\ss X$ and functions 
$$
u\in F^c(K)\ {\rm (some\ }c>0) \and  v\in \wt F(K),
$$
   the  Zero Maximum Principle holds for $u+v$ on $K$.
 We say that {\sl local weak comparison holds} for $F$ on $X$ if every
$x\in X$ has  some neighborhood  on which weak comparison holds.

\medskip

It is proved (see  Theorem 8.3)  that:
 \medskip

\centerline {\sl Local weak comparison implies weak comparison.}
\medskip
\noindent
Then in Section 10 we prove that:
\medskip
\centerline{\sl Local weak comparison holds for any subequation which is  }

\centerline{\sl  locally affinely jet-equivalent to a constant coefficient subequation.}

\medskip\noindent
Taken together we have that global weak comparison holds for all such subequations.
This constitutes a very broad class.  In particular, we have that 
\medskip
\centerline{\sl Weak comparison holds for riemannian $G$-subequations }

\centerline{\sl  on a riemannian manifold.}
\medskip
\vfill\eject

We then establish, under certain global conditions, that weak comparison implies comparison.

\Theorem{9.7} {\sl Suppose $F$ is a   subequation for which local
weak comparison holds.
Suppose there exists a 
$C^2$ strictly $M_F$-subharmonic function on $X$ where $M_F$ is a monotonicity
cone for $F$.
Then  comparison holds for $F$ on $X$.}

\medskip

Combined with our notions of boundary convexity, we prove the following.

\Theorem{13.3} {\sl Suppose comparison holds for a  subequation   $F$ on $X$

Then for every domain $\O\ss\ss X$ whose boundary is strictly
$F$- and $\ft$-convex, both existence and uniqueness
hold for the Dirichlet problem.}\medskip

In Section 9 we introduce the concept of strict approximation for $F$ on $X$
and show that if $X$ admits a $C^2$ strictly $M_F$-subharmonic function
(where $M_F$ is a monotonicity cone for $F$ as above), then
strict approximation holds for $F$ on $X$.  Furthermore, we show that weak comparison
plus strict approximation implies comparison.

\vskip .3in
\centerline{\headfont Boundary Convexity}\medskip

An important part of this paper is the formulation and study of the notion
of boundary convexity for a general fully non-linear second-order equation.
This concept is presented in Section 11.  It strictly generalizes the boundary
convexity defined in [CNS] and in [HL$_{4}$]. 

 In  the Grassmann examples
this boundary convexity condition is particularly transparent and geometric.
Suppose $F(\GG)$ is a  
purely second-order equation, defined as in Example E by a subset $\GG\ss G(p,TX)$  of the Grassmann bundle. Then boundary convexity for 
a domain $\O$ becomes the
requirement that the second fundamental form $II_{\bo}$ of $\bo$  satisfy
$$
\tr_\xi II_{\bo}\ \geq\ 0 
\eqno{(1.5)}
$$ 
for all $\x\in\GG$ such that $\x\ss T(\bo)$. (When 
 there are no $\GG$-planes in $T_x\bo$,  $\GG$-convexity automatically holds
at $x$.) This convexity automatically implies convexity for the dual subequation.

 Domains with convex boundaries in the riemannian sense
are $F$-convex for 
any purely second-order subequation. 

For subequations which do  no involve the  dependent variable, $F$-convexity
is defined in terms of the {\sl asymptotic interior} of $F$.  This is an open, point-wise
conical set consisting of rays with conical neighborhoods which eventually lie in $F$.
For general subequations, we freeze the dependent variable to be a constant $\l$
(e.g. replace $\D u \geq e^u$ with $\D u \geq e^\l$) and use the associated asymptotic
interiors to define convexity. 

We then derive  a  simple, geometric criterion for $F$-convexity in
terms of the second fundamental form of $\bo$.  In almost all interesting cases
the boundary convexity condition is straightforward to compute.

Some subequations $F$ have the property that 
every boundary is $F$-convex. These  include
the $p$-Laplace-Beltrami subequation,
$$
\|\nabla u\|^2\Delta u + (p-2)(\nabla u)^t (\Hess\, u) (\nabla u) \geq\ 0
$$
for $1 <  p<\infty$, and the infinite Laplace-Beltrami subequation
$$
(\nabla u)^t (\Hess\, u) (\nabla u) \ \geq\ 0.
$$

 For the general minimal surface subequation 
$(1+\|\nabla u\|^2)\D u - (\nabla u)^t (\Hess\, u) (\nabla u) \geq 0$
strict boundary convexity is equivalent  to strictly positive mean curvature with respect to the interior normal.
For the Calabi-Yau equation, and also for the homogeneous complex Monge-Amp\`ere equation,
$F$-convexity of the boundary is simply standard pseudo-convexity (in the hermitian almost-complex case).
 
Finally we note that for a general $F$ {\sl  all sufficiently small balls  in any local coordinate system are strictly
$F$-convex} under the mild condition that $F$ contains critical jets (see Proposition 11.9).

\vskip .3in
\centerline{\headfont  Existence without Comparison}\medskip

Certain methods taken from complex analysis, which go back to Bremermann/Perron 
and then Walsh ([B], [W]), enable us to prove existence in the   constant coefficient setting
 with no assumptions about a monotonicity cone.  This includes many cases where uniqueness
 does not hold.  More generally, we prove in Section 12 that if $X = K/G$
  is a riemannian homogeneous space
and $F\ss J^2(X)$ is a $K$-invariant subequation, then existence holds for all
continuous boundary data on any domain which is strictly $F$- and 
$\ft$-convex. This applies in particular to the $p$-Laplace-Beltrami
and infinite Laplace-Beltrami subequations mentioned above, where
all domains are strictly $F$- and 
$\ft$-convex.

The euclidean version of this result is stated as Theorem 12.7.  
It establishes existence for any constant coefficient subequation
$\bbf$ on all strictly ${\bbf} $ and ${\wt \bbf}$-convex domains
(when they exist).
We follow this with 
an Example 12.8 of a second order equation where uniqueness does in fact fail.

\vskip .3in
\centerline{\headfont Further Examples.}\medskip

There are many geometrically interesting subequations which are covered by the
results above. We examine a few more examples here. We start with a general observation.

\medskip\noindent
{\bf Inhomogeneous Equations.}  The methods above apply to any subequation which is locally
{\sl affinely} jet-equivalent to a constant coefficient equation. This greatly extends the equations 
that one can treat.  For example, suppose that $F$ is a riemannian $G$-subequation
with a monotonicity cone $M$ on a manifold $X$, and let $J$ be any smooth section of 
the 2-jet bundle $J^2(X)$.  Then $F_J \equiv F +J$ (fibre-wise translation) is also a subequation
having the same asymptotic interior  $\overrightarrow{F_J} = \Fa$ and also having $M$ as a monotonicity
cone.

As a simple but interesting example, suppose that $F$ is one of the branches of the 
homogeneous Monge-Amp\`ere equation (in  the real,  complex or quaternionic case), 
and choose $J_x = f(x) I$ where $f$ is an arbitrary smooth function on $X$ (and $I$ is the identity
section of $\Sym(T^*X)$).  Then one can treat the inhomogeneous equation
$$
\l_k(\Hess\, u)\ =\ f(x)
$$
under the same conditions that one can treat the homogeneous equation 
 $\l_k(\Hess\, u)=0$.  For the principal branch $\cp$, $\cp^\bbc$ etc.,
 and $f\equiv -1$, this yields functions which are quasi-convex, quasi-psh, etc..
 The higher branches with variable $f$ are more interesting.

Similar remarks hold for any purely second-order subequation, such as those below.
However, many other interesting subequations arise from local affine jet-equivalence
(cf. Example 6.15).

\medskip\noindent
{\bf Example F.  (Almost calibrated manifolds).}  An important class of  Grassmann structures
(discussed in Example E) are given  by   calibrations. A  constant coefficient
calibration is a   $p$-form $\phi$ with the property that 
$\f(\x)\leq1$ for all $\x\in G(p,\rn)$. The associated set is
$$
\GG(\f)\ \equiv\ \{\x\in G(p,\rn) : \f(\x)=1\}
$$
and the associated group is 
$$
G(\f) = \{g\in {\rm O}_n : g^*(\f)=\f\}.
$$
Any riemannian manifold with a topological $G(\f)$-structure will
carry a global $p$-form $\wt\phi$, called an {\sl almost calibration}, which is of type $\f$ at every point
but is not necessarily closed under exterior differentiation.   
Some of the subequations already referred to can be defined in this way.
For example, almost complex hermitian geometry arises from 
$\phi=\o$, the standard  (not necessarily closed)  K\"ahler form.  There are however many others.
Several interesting examples are given next.

\medskip\noindent
{\bf Example G.  (Almost Hyperkahler Manifolds).}
Here we consider a $4n$-dimensional manifold $X$  
equipped with a subbundle ${\cal Q}\ss \Hom(TX,TX)$
generated by {\bf global} sections $I,J,K$ satisfying the standard quaternion relations:
$I^2=J^2=K^2=-1$, and $ IJ=K$, etc. and equipped with a compatible
 riemannian metric $g$. The topological structure group is  $ {\rm Sp}_n$.
This is a special case of Example D above so the quaternionic 
plurisubharmonic functions are defined, and when the Hyperkahler structure is integrable, they 
coincide with the ones  used by Alesker and Verbitsky  to study the 
Quaternionic Monge-Ampere equation in the principal-branch  case. (See  [A$_{1,2}$], [AV])

However a  manifold with a topological  Sp$_n$-structure carries other almost calibrations  
such as the generalized Cayley calibration:
$
\wt\Omega \ =\  \smfrac12\left\{ \o_I^2+\o_J^2-\o_K^2\right\}
$ 
introduced in [BH].

\medskip\noindent
{\bf Example H.  (Almost Calabi-Yau Manifolds).} 
 This is an almost complex hermitian
manifold with a global section of $\Lambda^{n,0} (X)$ whose real part $\Phi$ has  
comass $\equiv 1$. This is equivalent to having topological structure group SU$_n$.
In addition to the structures discussed in Case 2, these manifolds carry
functions associated with the {\sl Special Lagrangian calibration} $\Phi$.
The $\Phi$-submanifolds are called  {\sl Special Lagrangian submanifolds} 
and the $\Phi$-subharmonic functions are said to be  {\sl Special Lagrangian subharmonic.}

\medskip\noindent
{\bf Example I.  (Almost G$_2$  manifolds).} 
 Let  Im$\bbo = \bbr^7$ denote the imaginary
Cayley numbers.  This space is acted on by the group G$_2$ of automorphisms of $\bbo$
which preserves the 3-form
$$
\vf(x,y,z) \ =\ \langle x\cdot y, z\rangle
$$
called the {\sl associative calibration}. Any 7-dimensional manifold $X$ with a 
topological G$_2$-structure carries an associated riemannian metric and
a global (non-closed) calibration $\vf$. There also exists a degree 4 calibration
$\psi = *\vf$ on $X$.  One then has $\vf$ and $\psi$ subharmonic and harmonic functions
on $X$ and one can consider the Dirichlet Problem on bounded domains.
We note that these structures exists in abundance. In [LM, page 348] it is shown that
for any 7-manifold $X$
$$
X\ \ {\rm has \  a \ topological \ } {\rm G}_2\ {\rm structure\ } \qquad \iff \qquad X\ \ {\rm is\ spin}.
$$

\medskip\noindent
{\bf Example J.  (Almost Spin$_7$ manifolds).} 
  On the Cayley numbers $\bbo=\bbr^8$ there
is a 4-form of comass one defined by
$$
\Phi(x,y,z,w) \ =\ \langle( x\cdot y)\cdot z - x\cdot (y\cdot z), w\rangle
$$
and preserved by the subgoup Spin$_7 \ss {\rm SO}_8$ (cf. [HL$_1$], [H], [LM]).
This determines a non-closed calibration $\Phi$ on any 8-manifold $X$ with
a topological Spin$_7$-structure.  In [LM, page 349] it is shown that
for any 8-manifold $X$
$$
X\ \ {\rm has \  a \ topological \ } {\rm Spin}_7\ {\rm structure\ } \quad \iff \quad X\ \ {\rm is\ spin\ and\ }
p_1(X)^2 -4p_2(X) +8\chi(X)=0
$$
for an appropriate choice of orientation on $X$. Here $p_k(X)$ is the $k$th Pontrjagin class
and $\chi(X)$ denotes the Euler class of $X$.

\medskip\noindent
{\bf Example K.  (Lagrangian subhamonicity).}  
Suppose $(X,J)$ is an almost complex hermitian manifold of real dimension $2n$.
  Then, as mentioned in Example E,
there is a natural Grassmann structure
 $\LAG \ss G(n, TX)$ consisting of the Lagrangian
 $n$-planes. This gives rise to the Lagrangian subequation defined
in terms of the Riemannian hessian by the 
condition that
$$
\tr \left\{\Hess\, u\bigr|_\xi\right\} \ \geq \ 0 \fa \x\in \LAG
$$
Interestingly, there is a beautiful polynomial operator 
which vanishes on the Lagrangian harmonic functions, and so the Dirichlet problem here
can be considered to be for this operator.  Furthermore, just as in the Monge-Amp\`ere case,
this operator has many branches, each of which is another U$_n$-invariant subequation
on the manifold.  This follows from G\aa rding's theory of hyperbolic polynomials [G], [HL$_4$], [HL$_7$].

\medskip\noindent
{\bf Example L.  (Equations involving the dependent variable but independent of the gradient).}  
Many  purely second-order
equations can be enhanced to ones which involve the dependent variable $u$ and
our theory continues to apply. For example, consider the O$_n$-universal subequation
$F$ defined by requiring
$$
f( u, \l_q(A))\ \geq\ 0
$$
where $f(x,y)$ is non-increasing in $x$ and non-decreasing in $y$.  As a special
case,  consider
$$
\l_q(A) - \vf(u) \ \geq\ 0
$$
where $\vf $ is monotone non-decreasing.  The dual subequation $\ft$ is given
by
$$
\l_{n-q+1} + \vf(-u) \ \geq\ 0.
$$
If $\vf(0)=0$, then the required convexity of the boundary is that
$$
\min\{\l_q(II_{\bo}), \l_{n-q}(II_{\bo})\} \   > \ 0.
$$

Another interesting O$_n$-universal subequation  is
$$
F \ \equiv\ \{(r,p, A) : A\geq 0 \ {\rm and\ }\det A - e^r \geq 0\}
$$
which is discussed in Remark 12.9.  Its dual equation is 
$$
\ft\ =\ \{(r,p,A) : -A\not< 0 \ {\rm or} \  -|\det A|+e^{-r} \geq 0\}.
$$

\vskip .3in
\centerline{\headfont  G\aa rding hyperbolic polynomials}\medskip

G\aa rding's  beautiful  theory of hyperbolic polynomials [G], when applied to homogeneous 
polynomials  $M$   on $\Symn$, fits perfectly into this paper.
It unifies and generalizes many of our constructions.  We give here a brief sketch of
how this works, and refer to  [HL$_7$] for full details.

By definition a homogeneous polynomial $M$ of degree $m$ on $\Symn$ is {\sl hyperbolic
with respect  to the identity $I$}   if for all $A\in \Symn$ the polynomial
$s\mapsto M(sI+A)$ has exactly $m$ real roots. The negatives of these roots are called the
{\sl $M$-eigenvalues of $A$}.  It is useful to order these eigenvalues
$$
\l_1^M(A)\ \leq\ \cdots\ \leq \l_m^M(A)
$$
and normalize so that $M(I)=1$. Then the polynomial factors as
$
M(sI+A)\ =\ \prod_{k=1}^m\left(s+\l_k^M(A)\right)
$.
Using the ordered eigenvalues we can define branches 
$$
\bbf^M_k\ \equiv\ \{A\in \Symn : \l_k^M(A)\geq0\}
$$
which satisfy $\bbf^M_1\ss\bbf^M_2\ss\cdots$.  
The {\sl principal branch} $\bbf^M \equiv\bbf_1^M$ is the connected component 
of $\{M\neq 0\}$ containing $I$.
G\aa rding proves that:
\medskip

(1)\ \ The principal branch $\bbf^M$ is a convex cone.

\medskip

(2)\ \ $\bbf^M_k+\bbf^M \ \ss\ \bbf^M_k$ for all $k$.\medskip

\noindent
Under the positivity assumption: $\bbf^M + \cp \ss \bbf^M$ we then have that 
\medskip

\centerline{\sl  Each $\bbf^M_k$
is a constant coefficient subequation for which $\bbf^M$ is a monotonicity cone. }

\medskip\noindent
When $M=\det$ we get the
branches of the Monge-Amp\`ere equation (Example B).
However, this applies to many other interesting
cases (such as  Examples C, D and K). Furthermore, for each subset
$\L\ss\rn$ as in 14.1 below, one can construct  a subequation $\bbf^{M,\L}$
using the $\l_k^M$ (see [HL$_7$]).

If the polynomial $M$ is invariant under a subgroup
$G\ss {\rm O}_n$, then  these subequations carry over to any manifold with a 
topological $G$-structure.

\vskip .3in
\centerline{\headfont  Parabolic Subequations}\medskip

The methods and results of this paper carry over effectively to parabolic equations.
Suppose $X$ is a riemannian manifold  equipped with a riemannian
$G$-subequation $F$ for $G\ss {\rm O}_n$. We assume $F$ is induced from a universal
model 
$$
\bbf \ =\ \{J\in \bbj^2 : f(J)\geq0\}
$$
where $f: J^2(X) \to \bbr$  is $G$-invariant,  $\cp$- and $\cn$-monotone,
and Lipschitz in the reduced variables $(p,A)$.
Then on the riemannian product $X\times \bbr$  we have the associated $G$-universal parabolic
subequation $H$ defined by
$$
f(J)-p_0\ \geq\ 0
$$
where $p_0$ denotes the $u_t$ component of the 2-jet of $u$. The $H$-harmonic functions
are solutions of the equation $u_t = f(u, Du, D^2 u)$. Interesting examples which can be
treated include:

\smallskip
 (i) \  $f = \tr A$, the standard heat equation $u_t = \D u$ for  the Laplace-Beltrami operator on $X$.

\smallskip
 (ii)     $f = \l_q(A)$, the $q$th ordered eigenvalue of $A$.  This is 
 the natural parabolic equation associated to the $q$th branch of the Monge-Amp\`ere equation.

\smallskip
 (iii)  \  $f = \tr A +  { k\over  |p|^2+\e^2}  p^t A p$ for $k\geq - 1$ and $\e > 0$.
When $X=\rn$ and $k=-1$, 
the solutions $u(x,t)$ of the associated parabolic equation, in the limit as
$\e\to0$, 
have the property that the associated level sets
$
\Sigma_t\ \equiv\ \{x\in\rn : u(x,t)=0\}
$
are evolving by mean curvature flow. (See the classical papers
of Evans-Spruck [ES$_*$] and Chen-Giga-Goto [CGG$_*$] , and the very nice accounts in [E] and [Gi].)

\smallskip
 (iv)  $f = \tr\{\arctan A\}$.  
 When $X=\rn$,  solutions $u(x,t)$ 
 have the property that the graphs of the gradients:
 $
 \G_t\ \equiv\ \{(x,y)\in \rn\times \rn =\bbc^n : y = D_xu(x,t)\}
 $
are Lagrangian submanifolds which evolve the initial data by 
mean curvature flow.  (See [CCH] and references therein.)
\smallskip

Straightforward application of the  techniques in this paper shows that:\medskip
\centerline
{\sl Comparison holds for all riemannian $G$-subequations $H$ on $X\times \bbr$.
}
\medskip

By standard techniques one can prove more. 
Consider a compact subset  $K \ss \{t\leq T\}\ss X\times \bbr$ and 
let $K_T \equiv K\cap \{ t=T\}$ denote the terminal time slice of $K$.  Let
$
\partial_0 K \  \equiv \ \partial K  - \Int K_T
$
denote the {\sl parabolic boundary of $K$}.  Here $\Int K$ denotes the relative interior in 
$\{t=T\}\ss X\times\bbr$.
 We say that {\sl parabolic comparison holds for $H$} if for all such $K$ (and $T$)
$$
u+v\ \leq\ p \quad{\rm on}\ \ \partial_0 K\qquad\Rightarrow\qquad u+v\ \leq\ p \quad{\rm on}\ \ \Int K
$$
for all $u\in H(K)$ and $v\in \wt H (K)$.  Then one has that:
\medskip
\centerline{\sl Parabolic comparison holds for all riemannian $G$-subequations $H$ on $X\times \bbr$.
}
\medskip

Under further mild assumptions on $f$ which are satisfied in the examples above, 
one also has existence results.   Consider a domain $\O \ss X$ whose boundary is
strictly $F$- and ${\wt F}$-convex.  Set $K= \overline\O \times [0,T]$.
Then 
\medskip
\centerline{\sl For each $\vf\in C(\partial_0 K)$ there exists a unique function $u\in C(K)$ such
that
}
\centerline{\sl $u\bigr|_{\Int K}$ is $H$-harmonic\ \  and\ \  $u\bigr|_{\partial_0 K} = \vf$.
}
\medskip
\noindent
One also obtains corresponding long-time existence results.  Details will appear 
elsewhere.

\medskip\noindent
{\bf A brief outline of the paper.}  Section 2 along with Appendices A and B provide a self-contained treatment 
of general $F$-subharmonic functions and their properties.  Section 3 introduces the concept
of a  {\sl subequation} and discusses the natural duality among subequations.
 Riemannian manifolds are considered in Section  4 where it shown that
there is a natural splitting of the 2-jet bundle induced by the riemannian hessian.
This splitting gives many geometric examples of subequations which are purely
second-order.   In Section 5 ``universal''  subequations are constructed.
Suppose $\bbf$ is a (euclidean) constant coefficient subequation with compact invariance group
$G\ss {\rm O}_n$, and $X$ is a riemannian manifold   with a 
topological $G$-structure.  Then it is shown that the euclidean model $\bbf$ induces a {\sl riemannian $G$-subequation}  $F$ on $X$ with 
the property that $F$ is locally jet-equivalent to  $\bbf$. 
Section 6 discusses  automorphisms of the 2-jet bundle and the general  notion of jet-equivalence
for subequations.
 In Section 7 a notion of {\sl strictly}  $F$-subharmonic functions is introduced
for upper semi-continuous functions.
In Section 8 a weak form of comparison  is defined using this  notion 
of strictness.  It  is shown 
that if weak comparison holds locally, then it is true globally.
Section 9 addresses the question of when weak comparison implies comparison.
We first note that this holds whenever approximation by strictly $F$-subharmonic functions
is possible. The main discussion concerns how this ``strict approximation'' can be 
deduced from a form of monotonicity.  This yields many geometric examples.
In Section 10 we prove that weak comparison holds for any subequation
which is locally (weakly) jet-equivalent to a constant coefficient subequation, and 
in particular for the riemannian $G$-subequations constructed in Section 5.
The proof relies on the Theorem on Sums stated in Appendix C.
In Section 11 we introduce the notion of the {\sl asymptotic interior} $\Fa$
of a reduced  subequation $F$.  This leads to the notion of strict $F$-boundary-convexity
which implies the existence of  barriers essential to existence proofs.
For riemannian $G$-subequations the existence of these barriers is actually equivalent
to strict $F$-convexity.  Section 12 addresses the existence question for the
Dirichlet problem.  Assuming strict  boundary convexity for both the subequation $F$
and its dual $\ft$, several existence theorems are proved.
Section 13 compiles and summarizes the results established for the Dirichlet Problem.

The remaining sections   are devoted to examples.  Section 14 examines O$_n$-universal subequations.
These are subequations that makes sense on any riemannian manifold. 
Particular attention is paid to subequations which are purely second-order.
Analogous results in the complex and quaternionic case are examined in Section 15.
Section 16 discusses the subequations $F(\GG)$ geometrically defined by a closed subset
$\GG$ of the grassmannian. Section 17 applies the theory to equations  defined in terms of
the principal curvatures of the graph of $u$. Section 18 applies affine jet-equivalence to 
give results for inhomogeneous and other equations.  Section 19 treats the Calabi-Yau type
equations on almost complex hermitian manifolds.

\medskip
We note that  in  [AFS]  and [PZ] standard viscosity theory has been  retrofitted
to  riemannian manifolds by using the distance function, parallel translation,
Jacobi fields, etc..   In our approach  this machinery
in not necessary. We get by with  the standard
viscosity techniques (cf.  [CIL], [C]).

\vskip .5in


\centerline{\headfont \ 2.\   F-Subharmonic  Functions.}
\medskip
Let $X$ be a smooth $n$-dimensional manifold. Denote by $J^2(X) \to X$ the bundle
of 2-jets whose fibre at a point $x$ is the quotient
$$
J^2_x(X)\ =\ C^\infty_x/ C^\infty_{x,3}
$$
where $C^\infty_x$ denotes the germs of smooth functions at $x$ and 
$C^\infty_{x,3}$ the subspace of germs which vanish to order three at $x$.
Given a smooth function $u$ on $X$, let  $\J^2_x u\in J^2_x(X)$ denote its 2-jet at $x$. 
and note that  $J^2u$ is a smooth section of the bundle $J^2(X)$.

The bundle on 1-jets $J^1(X)$ is defined similarly and has a natural
splitting $J^1(X) =  \bbr \oplus T^*X$ with $\J^1_x u = (u(x), (du)_x)$.
  There is a short exact sequence of bundles
$$
0 \to \Sym(T^*X) \to J^2(X)  \to J^1(X)  \to 0 
\eqno{(2.1)}
$$
Here the space  $\Sym(T^*_xX)$ of symmetric bilinear forms on $T_xX$
 is embedded as the space of 2-jets of functions with critical
value zero at the point $x$, i.e.,  

\medskip
\centerline
{$\Sym(T^*_xX)\cong    \{\J^2_x u :u(x)=0, (du)_x=0\}.$}
\medskip
  Note that if $u$ is such a function, and $V,W$ are
vector fields near $x$, then $(\Hess_x u)(V,W) = V\cdot W\cdot u
= W\cdot V\cdot u +[V,W]\cdot u = W\cdot V\cdot u= 
(\Hess_x u)(W,V)$ 
is a well-defined symmetric form on $T_xX$.
However, for functions $u$ with $(du)_x\neq0$, there is no natural definition
of $\Hess_x u$, i.e., the sequence (2.1) has no natural splitting.
Choices of splittings correspond to definitions of a hessian, and there is a canonical one
for each riemannian metric as we shall see in \S 4.

At a minimum point $x$ for a smooth function $u$, we have $(du)_x=0$
(so that $\Hess_x u\in\Sym(T^*X)$ is well defined), and $\Hess_x u\geq0$. 
The isomorphism 
$$
\{H\in \Sym(T^*X) : H\geq 0\}\ \cong\ \{J^2_x u : u \geq0  \ \ {\rm near\ } x \ {\rm and \ } u(x)=0  \}
\eqno{(2.2)}
$$
defines a cone bundle
$$
\cp\ \ss\  \Sym(T^*X)\ \ss\ J^2(X)
$$
with $\cp_x$ defined by (2.2).

Given an arbitrary subset $F\ss J^2(X)$ a function $u\in C^2(X)$ will be called
{\bf $F$-subharmonic} if its 2-jet satisfies
$$
J^2_x u\ \in\ F_x \fa x\in X,
$$
and {\bf strictly  $F$-subharmonic} if its 2-jet satisfies
$$
J^2_x u\ \in\ (\Int F)_x \fa x\in X,
$$
These notions are of limited interest for general sets $F$

\Def{2.1} A subset $F\ss J^2(X)$ satisfies the {\bf   Positivity Condition}  if
$$
F+\cp\ \subseteq \ F
\eqno{(P)}
$$
A  subset $F\ss J^2(X)$ which satisfies (P) will be called 
{\bf $\cp$-monotone}.
\medskip

Note that condition (P) implies that
$$
\Int F + \cp \ \ss\ \Int F.
$$

Monotonicity is a key concept in this paper.  For arbitrary subsets $M,F\ss J^2(X)$
we say that $F$  is {\bf $M$-monotone} or that $M$ is a {\bf montonicity set for $F$} if 
 $$F+M\ss F$$.

It is necessary and quite useful to extend the  definition of $F$-subharmonic 
to non-differentiable functions $u$.
Let $\USC(X)$ denote the set of $[-\infty, \infty)$-valued, upper semicontinuous 
functions on $X$.

\Def{2.2} A function $u\in \USC(X)$ is said to  be  $F$-{\bf subharmonic} 
if for each $x\in X$ and each function $\vf$ which is $C^2$ near $x$, one has that
$$
\left.
\cases
{ u-\vf  \ &$\leq$ \ 0 \    \quad {\rm near}\ $x_0$  and \cr 
  \ &= \ 0\ \qquad {\rm at}\ $x_0$  
 } 
 \right\} 
 \qquad \Rightarrow\qquad
 \J^2_x \vf \ \i \ F_x.
\eqno{(2.3)}
$$
We denote by $F(X)$ the set of all such functions.
\medskip

Note that if $u\in C^2(X)$, then
$$
u\in F(X)\quad\Rightarrow\quad \jt_{x} u \in F_{x} \qquad \forall \ x\in X
\eqno{(2.4)}
$$
since the test function $\vf$ may be chosen equal to $u$ in (2.3).  The converse   is not
  true for general subsets $F$.  However, we have the following.

\Prop{2.3} {\sl Suppose  $F$ satisfies the  Positivity Condition (P) and $u\in C^2(X)$.
Then}
$$
\jtx u \ \in\ F_x \fa x\in X \qquad\Rightarrow\qquad u\ \in\ F(X).
\eqno{(2.5)}
$$

\pf
Assume $\jt_{x_0} u \in F_{x_0}$ and  and $\vf$ is a $C^2$-function such that  
$$
\left.
\cases
{ u-\vf  \ &$\leq$ \ 0 \    \quad {\rm near}\ $x_0$   \cr 
  \ &= \ 0\ \qquad {\rm at}\ $x_0$  
 } 
 \right\} 
$$
Since $(\vf - u)(x_0) = 0$,
$d(\vf-u)_{x_0}  = 0$, and $\vf-u \geq 0$ near $x_0$,  we have $J_{x_0}^2(\vf-u) \in \cp_{x_0}$ by definition.
 Now the Positivity Condition implies that  
 $\jt_{x_0}\vf \in  J^2_{x_0} u +\cp_{x_0} \ss F_{x_0}$.  This proves
that $u\in F(X)$.   \qed  \medskip

 The  proof  above  shows that  {\sl we must assume 
$F$ satisfies the Positivity Condition}. Otherwise the definition of  $F$-subharmonicity would  not
extend the natural one for smooth functions.

There is an equivalent definition of $F$-subharmonic functions which is quite useful.
We record it here and prove it in Appendix A (Proposition A.1 (IV)).

\Lemma{2.4}  {\sl  Fix $u\in \USC(X)$, and suppose that $F\ss J^2(X)$ is closed. 
 Then  $u\notin F(X)$ if and only if 
there exists a point $ x_0\in X$, local coordinates $x$ at $x_0$, $\a >0$ and 
a quadratic function 
$$
q(x) \ =\ r+\langle p, x-x_0\rangle +\half \langle A(x-x_0), x-x_0 \rangle
$$
with $\J^2_{x_0} (q)\notin F_{x_0}$ so that}
$$
\eqalign
{
u(x)-q(x) \ &\leq \ -\a|x-x_0|^2 \qquad
{\rm near\ } x_0 \qquad{\sl and}  \cr
&=\ \ 0\qquad\qquad\qquad\ \ \ {\sl at}\ \ x_0
}
$$

\Remark{2.5} 
(a).  The hypothesis that $F$ be closed in Lemma 2.4 can be substantially weakened.
It suffices to assume only  that the fibres of $F$ under the map 
$J^2(X)\to J^1(X)$ are closed. In fact, Proposition A.1
is proved under this  weaker hypothesis.

(b). The positivity condition is rarely used in proofs.  This is because without it $F(X)$ is empty and the results are trivial. For example, the positivity condition is not required in the following theorem. ($F$ need only be closed.)

\medskip

It is remarkable, at this level of generality,   that $F$-subharmonic
 functions  share many of the important properties
of classical subharmonic functions.   

\Theorem{2.6. Elementary  Properties of  F-Subharmonic Functions}
{\sl
Let $F$ be an arbitrary  closed subset of $J^2(X)$.
\medskip

\item{(A)}  (Maximum Property)  If $u,v \in F(X)$, then $w=\max\{u,v\}\in F(X)$.

\medskip

\item{(B)}     (Coherence Property) If $u \in F(X)$ is twice differentiable at $x\in X$, then $\jtx u\in F_x$.

\medskip

\item{(C)}  (Decreasing Sequence Property)  If $\{ u_j \}$ is a 
decreasing ($u_j\geq u_{j+1}$) sequence of \ \ functions with all $u_j \in F(X)$,
then the limit $u=\lim_{j\to\infty}u_j \in F(X)$.

\medskip

\item{(D)}  (Uniform Limit Property) Suppose  $\{ u_j \} \ss F(X)$ is a 
sequence which converges to $u$  uniformly on compact subsets to $X$, then $u \in F(X)$.

\medskip

\item{(E)}  (Families Locally Bounded Above)  Suppose $\cf\subset F(X)$ is a family of 
functions which are locally uniformly bounded above.  Then the upper semicontinuous
regularization $v^*$ of the upper envelope 
$$
v(x)\ =\ \sup_{f\in \cf} f(x)
$$
belongs to $F(X)$.

}

\pf
 See Appendix B.
 \medskip
 
   \medskip
 \noindent{\bf Cautionary Note 2.7.}  Despite the elementary proofs of the properties in Theorem
 2.6, illustrating  how well adapted Definition 2.2 is to nonlinear theory, 
 there are difficulties with the linear theory.  If $F_1$ and $F_2$ are $\cp$-monotone subsets,
 then the  (fibrewise) sum $F_1+F_2$ is also obviously a $\cp$-monotone subset.
 However, the property 
 $$
 u\in F_1(X), v\in F_2(X) \quad\Rightarrow \quad u+v \in (F_1+F_2)(X)
 $$
 is difficult to deduce   from Definition 2.2 even in the basic case where
 $F_1 = F_2 = F_1+F_2$ is the linear ``subequation'' on $\rn$ defined by $\D u\geq 0$.

 \vskip .3in


\centerline{\headfont \ 3.\  Dirichlet Duality  and the Notion of a Subequation.}
\medskip

The following concept is the lynchpin for the Dirichlet Problem.

\Def{3.1}  Given a subset  $F\ss J^2(X)$
the {\bf Dirichlet dual} $\wt F$ of $F$ is     defined by
$$
\wt F\ =\ \sim(-\Int F)\ =\ -(\sim \Int F)
\eqno{(3.1)}
$$

Note the obvious properties:
 
\medskip

 (1)  \ \ $F_1 \ \ss \ F_2$     \quad$\Rightarrow$\quad  $\wt F_2 \ \ss \ \wt F_1$. \smallskip

\smallskip

 (2)  \ \  $\wt {F_1 \cap F_2} \ =\ \wt F_1 \cup \wt F_2$   \smallskip
 
\smallskip

 (3)  \ \  $\wt { \wt F} \ =\  F$ \qquad {\rm if \ and\ only\ if \  \  \ } $F \ =\ \overline{\Int F}$.
 
\smallskip

Thus to have a  true duality with $\wt{\wt F}\ =\ F$ we must assume that $F= \overline{\Int F}$.
For simplicity we also want to compute the dual fibrewise in the jet bundle.
Consequently we will assume the following  three topological conditions on $F$, combined as condition (T).

$$
 (i)\ \ F\ =\ \overline{\Int F}, \qquad (ii)\ \ F_x\ =\ \overline{\Int F_x}, \qquad 
 (iii)\ \ \Int F_x\ =\     (\Int F)_x
 \eqno{(T)} 
 $$
\smallskip
 \noindent
 (where $\Int F_x$ denotes interior with respect to the fibre).
 It is then easy to see that the fibre of $\ft$ at $x$ is given by 
$-(\sim \Int F_x) = \sim(-\Int F_x) = \wt{(F_x)}$, so there is no ambiguity in the
notation $\ft_x$.
 
 \Def{3.2}  A subset $F   \ss J^2(X)$  satisfying  (T) will be called a {\bf T-subset}.
 
 \Lemma{3.3} {\sl Suppose  $F\ \ss\ J^2(X)$  has the property that $F=\overline{\Int F}$. Then
\medskip

 (a)  \ \ $F$ satisfies Condition (P)   \quad$\iff$\quad $\wt F$ satisfies Condition (P). \smallskip

\smallskip

 (b)  \ \ $F$ satisfies Condition (T)    \quad$\iff$\quad $\wt F$ satisfies Condition (T). \smallskip
 
 \smallskip
   }
 \pf
Assertion (b)  is straightforward.
To prove (a) we use another property.

 \Lemma{3.4} {\sl Suppose that $F$ is a $T$-subset.  Then
  \medskip

 (4)  \ \  $\wt {F_x+J} \ =\ \wt F_x - J $\qquad for all $J \in J^2_x(X)$
 }
 
 \pf  Fix  $J \in J^2_x(X)$. Then  
 $J'\in \wt{F_x+J} \iff -J' \notin \Int(F_x+J) = \Int F_x +J
\iff -(J'+J)\notin \Int F_x \iff  J'+J\in \ft_x \iff J'  \in\ft_x-J$.
\qed

\Cor{3.5}  {\sl  Suppose  $F$ is a $T$-subset of $J^2(X)$ and $M$ is an
arbitrary subset of $J^2(X)$.  Then\medskip
\centerline
{ $F$ is $M$-monotone \quad$\iff$\quad $\wt F$ is $M$-monotone.}
}
\pf
Fix $J\in M_x$ and assume $F_x + J \ss F_x$, or equivalently 
$F_x\ss F_x-J$. By (1) this implies that $\wt{F_x-J}
\ss \ft_x$.  By (4) we have $\wt{F_x-J} = \ft_x+J$ so that $\ft_x +J \ss \ft_x$.
The converse follows from (3) and (b).
\qed
 \medskip

\noindent
{\bf Proof of (a).}  Take $M=\cp$ in Corollary 3.5.\qed\medskip
 
  Given a closed subset $F\ss J^2(X)$, note that 
 $$
 \partial F\ =\ F\cap (\sim \Int F) \ =\ F\cap (-\ft)
 \eqno{(3.2)}
 $$

 \Def{3.6}  A function $u$ is said to be 
 {\bf $F$-harmonic} if 
 $$
 u\ \in\ F(X) \and -u\ \in\ \wt F(X).
 $$
 
 Thus, under Condition (T) a function $u\in C^2(X)$ is $F$-harmonic if and only if
 $$
J^2_x u \ \in\ \partial F_x \fa x\in X.
 \eqno{(3.3)}
 $$
 Note also that an $F$-harmonic function is automatically continuous (since both $u$ and $-u$ are
 upper semi-contiuous).  The topological condition
 (T) implies that $\wt{\wt F} = F$ and hence that 
 \medskip
 \centerline{$u$ is $F$-harmonic \ \ $\iff$ \ \  $-u$ is $\ft$-harmonic.}
 \medskip
 
 The focal point of this paper is the {\bf Dirichlet Problem}, abbreviated (DP).
 Given a compact subset $K\ss X$ and a function $\vf\in C(\partial K)$, a function $u\in C(K)$ 
 is a solution to (DP) if 
 \medskip
 \centerline{ (1)  \ \ $u$ is $F$-harmonic on $\Int K$ \and (2) \ \ $u\bigr|_{\partial K} \ =\ \vf$.}

\Ex{3.7}  Consider the  one variable first order  subset $F$, defined by $|p|\leq 2|x|$,
and its dual $\ft$ defined by $|p|\geq 2|x|$. Of the conditions (P)  and (T), the subset
$F$ satisfies all but condition (ii) of (T) and its dual $\ft$  satisfies all but (iii) of (T).
 Both functions $x^2-1 $ and $-x^2+1$ are $F$-harmonic
 and have the same boundary values on $[-1,1]$.
 Thus without the full condition (T),  solutions of the  Dirichlet Problem may not be unique. 
\vskip .3in

\centerline{\bf Subequations}
\bigskip

Uniqueness for the Dirichlet problem requires a third hypothesis in addition to (P) and (T).
Consider the elementary   canonical splitting
$$
J^2(X) \ =\ \bbr\oplus J^2_{\rm red}(X)
$$
where $\bbr$ denotes the 2-jets of (locally) constant functions and 
$
J^2_{\rm red}(X)_{x} \equiv \{J^2_x u : u(x)=0\}
$
is the space of {\sl reduced} 2-jets at $x$, with the short exact sequence of bundles
$$
0\ \to\ \Sym(T^*X) \ \to\ J^2_{\rm red}(X)\ \to\ T^*X\ \arr\ 0
$$
which is the same as (2.1) except for the trivial factor $\bbr$.

 We define
$$
\cn\ \ss\ \bbr\ \ss\ J^2(X)
$$
to have fibres $\cn_x = \bbr^-= \{c\in \bbr : c\leq 0\}$.

\Def{3.8}  A subset $F\ss J^2(X)$
satisfies the {\bf Negativity Condition} if 
$$
F+\cn\ \subseteq \ F.
\eqno{(N)}
$$
A   subset $F\ss J^2(X)$ satisfying (N) will also be called  {\bf $\cn$-monotone}.
\medskip

Now we can add to Lemma 3.2 a third conclusion (assuming $F$ satisfies (T)):
\medskip

 (c)  \ \ $F$ satisfies Condition (N)    \quad$\iff$\quad $\wt F$ satisfies Condition (N) \smallskip
 
 \smallskip
   \noindent
   by taking $M=\cn$ in Corollary 3.5.
   
   In this paper the main results concerning existence and uniqueness for the Dirichlet problem
   assume that $F$ staisfies (P), (T) and (N).  This is formalized as follows.

\Def{3.9}  By a {\bf subequation} $F$ on a manifold $X$ we mean a subset
$$
F\ \ss\ J^2(X)
$$
which satisfies the three conditions (P), (T) and (N).

\Prop{3.10}  
$$
F\ \ {\sl is \  a\  subequation }\qquad \iff \qquad \ft\ \ {\sl is \  a\  subequation }
$$
\pf
See (a), (b) and (c) above.\qed
\medskip

Our investigation of the (DP) for a  subequation involves two additional subequations,
which are constructed from the original and have two additional properties. One is the 
{\bf cone property}:
$$
J\in F \quad\Rightarrow\quad tJ\in F \fa t>0,
$$
i.e.,  each fibre $F_x$ is a cone with vertex at the origin. 
Stronger yet is the 
{\bf convex cone property}:
\medskip
\centerline{each fibre $F_x$ is a convex cone with vertex at the origin.}
\medskip
\Def{3.11} A closed subset $F\ss J^2(X)$ having properties (P), (T), (N) and  the  
cone property will be called a {\bf    cone subequation}.   If is also has the convex
cone property it  will be called a {\bf  convex  cone subequation}.

\vskip .3in

\centerline{\bf Euclidean (or Constant Coefficient)  Subequations}
\bigskip

 The 2-jet bundle on $X=\rn$ is canonically trivialized by
 $$
 J^2_x u \ =\ \left(u(x), D_xu, D_x^2u  \right)
 \eqno{(3.4)}
 $$
 where
 $$
 D_x u \ =\ \left(\smfrac{\partial u}{\partial x_1}(x), ... ,  \smfrac{\partial u}{\partial x_n}(x)   \right)
 \and
  D_x^2 u \ =\  \left( \left(\smfrac{\partial^2 u}{\partial x_i\partial x_j}(x)\right)  \right)
 $$
 are the first and second derivatives of $u$ at $x$.  That is, for any open subset $X\ss\rn$
 there is a canonical trivialization
 $$
 J^2(X)\ \cong\ X\times \bbr\times\rn\times\Symn
  \eqno{(3.4)'}
 $$
 with fibre
 $$
 \bbj^2\ =\ \bbr\times\rn\times\Symn
 $$
 The standard notation $J=(r,p,A) \in\bbj^2$ will be used for the coordinates on $\bbj^2$.
 
 \Def{3.12} Any  subset $\bbf\ss \bbj^2$ which 
 satisfies  (P), (N) and for which  $\bbf=\overline{\Int \bbf}$ determines a 
   {\bf euclidean} or {\bf constant coefficient subequation} on any  
   open subset $X\ss \rn$  by setting  $F=X\times \bbf\ss J^2(X)$ (with constant fibre $F_x\equiv \bbf$).
 For simplicity we shall often denote this subequation simply by $\bbf$
 (with the euclidean coordinates implied), 
 and refer to the $F$-subharmonic functions as $\bbf$-subharmonic
 functions.
 
  \medskip
  Note that our assumptions on $\bbf$ imply the full condition  (T) for $F$.

In this paper our first concern will be to start with an arbitrary euclidean subequation $\bbf$ and then to 
describe those riemannian manifolds which have a variable coefficient subequation $F$
modeled on  $\bbf$. We will give an explicit construction of  this $F$. 
When, for example,   $\bbf$ is the euclidean Laplace subequation, this analogue can be
constructed on any riemannian manifold $X$ and is simply the Laplace-Beltrami subequation on $X$.

The general linear group $\GL_n(\bbr)$ has a natural action on $\bbj^2$ given by
$$
h(r,p,A) \ =\ (r, hp, hAh^t)    \qquad{\rm for}\ h\in \GL_n(\bbr).
\eqno{(3.5)}
$$
Each euclidean subequation $\bbf\ss\bbj$ has a {\bf  compact invariance group}
$$
G(\bbf)  \ =\  \{h\in {\rm O}_n : h(\bbf)=\bbf\}
\eqno{(3.6)}
$$
which will be instrumental in determining which riemannian manifolds can be fitted with an
analogous or companion subequation $F$.  When $G(\bbf) = {\rm O}_n$, no restriction on $X$ will be required.  However, even the extreme  case $G(\bbf) = \{I\}$ is not vacuous and will prove interesting.

\vskip .3in


\centerline{\headfont \ 4.\  The Riemannian Hessian -- A Canonical Splitting of $J^2(X)$.}
\medskip

\noindent
{\bf 4.1. The riemannian hessian.} 
Assume now that $X$ is equipped with a riemannian metric. 
Then, using the riemannian connection, any $C^2$-function $u$ on $X$ 
has a canonically defined {\bf riemannian hessian}
at every point given as follows.  Fix $x\in X$ and vector fields $V$ and $W$ defined
near $x$. Define
$$
(\Hess_x u) (V,W) \ \equiv\ V\cdot W\cdot u - \left(\nabla_V W\right)\cdot u
\eqno{(4.1)}
$$
where the RHS is evaluated at $x$. Since $(\nabla_V W)\cdot u - (\nabla_W V)\cdot u
= [V,W] \cdot u$, we see that  $(\Hess_x u) (V,W)$ is symmetric and depends only on the
values of $V$ and $W$ at the point $x$, i.e., $\Hess\,u$ is a section of $\Sym(T^*X)$.
Of course at a critical point the riemannian hessian always agrees with the hessian defined
at the beginning of  Section 2.

The subequations we shall consider will be obtained by putting constraints on the 
{\bf riemannian 2-jet} 
$$
(u, du, \Hess\,u)
\eqno{(4.2)}
$$
which combines the riemannian hessian with the exterior derivative.

\medskip
\noindent
{\bf 4.2. The canonical splitting.}
The riemannian hessian provides a bundle isomorphism
$$
 J^2(X) \ \harr{\cong}{}\ \bbr\oplus T^*X \oplus \Sym(T^*X)
\eqno{(4.3)}
$$
by mapping  
$$
J_x^2 u\ \mapsto\ \left( u(x), (du)_x, \Hess_x u\right)
$$
 for a $C^2$-function $u$ at $x$. Note that the right hand side depends only on the 2-jet 
of $u$ at $x$, and hence this is a well defined bundle map. It provides a canonical
splitting of the short exact sequence (2.1)

\medskip
\noindent
{\bf 4.3. Local trivializations associated to framings.}
The picture can be trivialized by a choice of local framing  $e=(e_1,...,e_n)$ 
of the tangent bundle $TX$ on some neighborhood  $U$.  
(This framing is {\sl not} required to be orthonormal.)
The canonical splitting (4.3) determines a   trivialization of $J^2(U)$
given at $x\in U$ by
$$
\Phi^e:  J^2_x(U) \ \arr\   \bbr \oplus \rn \oplus \Symn \qquad{\rm defined\  by}
\eqno{(4.4a)}
$$
$$
\Phi^e\left(J^2_x u\right)\ \equiv \ (u, e(u), (\Hess\, u)(e,e)) 
\eqno{(4.4b)}
$$
where $e(u)= (e_1 u,...,e_n u)$ and  $(\Hess\, u)(e,e)$
 is the $n\times n$-matrix with entries $\left(\Hess\, u \right) (e_i,e_j)$.

\Lemma{4.1 (Change of Frame)} {\sl 
Under a change of frame $e'=he$ (that is, $e_i' = \sum_j h_{ij} e_j$) where $h$ is a smooth
$\GL_n(\bbr)$-valued function, one has
$$
\left(u, e'(u),  (\Hess\, u)(e',e') \right)\ =\ \left(u, he(u), \,  h(\Hess\, u)(e,e) h^t \right).
\eqno{(4.5)}
$$
Said differently, for a 2-jet $J\in J^2_xX$}, 
$$
  \Phi^e(J) \ =\ (r,p,A)    \qquad  
 \Rightarrow\qquad  \Phi^{e'}(J) \ =\ (r,hp, hAh^t).
\eqno{(4.6)}
$$
The elementary proof is omitted.
\medskip

The rest  of this section presents some general constructions of subequations
using the riemannian hessian. For the sake of continuity the reader may want 
  to pass directly to  section 5.
  

\vskip .3in
\centerline{\bf Some Classes of Subequations Based on the Riemannian Hessian}
\medskip
\noindent 
{\bf 4.4. Pure Second-Order Subequations.}
These are the subequations which involve only
the riemannian hessian and not the value of the function or its gradient.  That is,
with respect to the riemannian splitting (4.3), a subequation of the form
$$
F\ =\ \bbr\oplus T^*X\oplus F',
$$
for a closed subset $F'\ss \Sym(T^*X)$ will be called {\bf pure (riemannian-)second-order}.
In other words, $J^2_xu\in F$ if and only if $\Hess_x u\in F'$.

We now observe that a closed subset $F'\ss \Sym(T^*X)$ determines a pure second-order  subequation 
 if and only if positivity  $ F'+\cp\ss F'$ holds and
 $$
  \Int F_x' \ = \ (\Int F')_x \fa x.
  \eqno{(4.7)}
 $$
 where $\Int F_x'$ denotes interior with respect to the fibre $J_x^2(X)$. 
Condition (N) is obvious.  It  is clear that condition (P) for $F$ is equivalent to
$F'+\cp \ss F'$ which implies that  
$\Int F'+ \Int \cp \ss \Int  F'$, because any $A\in \Int F_x'$ can be extended
  to a local  section of $\Int F'$.
Since  $A+\e I$ approximates any $A\in F'$, it follows that Condition (P) implies 
$$
(i)\ \ F\ =\ \overline{\Int F}\and  (ii)\ \ F_x\ =\ \overline{\Int F_x}.
$$     
Thus  with (4.7) condition (T) holds.

\medskip
\noindent {\bf   4.5. The complex and quaternionic hessians.} 
 Let $X$ be a riemannian manifold
equipped with a pointwise orthogonal almost complex structure $J:TX\to TX$. On this hermitian
almost complex manifold one can define the {\bf complex hessian } of a $C^2$-function $u$
by
$$
\Hess^\bbc u\ \equiv \ \half \left\{ \Hess\, u - J(\Hess\, u) J\right\}.
\eqno{(4.8)}
$$
This is a hermitian symmetric quadratic form on the complex tangent spaces of $X$. 
In particular its eigenvalues are real with even multiplicity  and its eigenspaces are $J$-invariant.
(See \S 15 for more details.)

Analogously suppose $X$ is equipped with a hermitian almost quaternionic structure, i.e., 
orthogonal bundle maps $I,J,K:TX\to  TX$ satisfying the standard quaternionic 
identities: $I^2=J^2=K^2=-1, \ IJ=-JI=K$, etc.  Then one can define the {\bf quaternionic hessian }
$$
\Hess^\bbh u\ \equiv \ \smfrac 1 4  \left\{ \Hess\, u - I(\Hess\, u) I- J(\Hess\, u) J- K(\Hess\, u) K\right\}.
\eqno{(4.9)}
$$
A number of basic subequations are defined in terms of these hessians.
Primary among them are the following.

\medskip
\noindent {\bf   4.6. The Monge-Amp\`ere subequations.}
  A classical subequation, which
is defined on any riemannian manifold $X$ and 
will play an important role in this  paper, is the (real) {\bf Monge-Amp\`ere subequation}
$$
P\ \equiv\  P^\bbr \ \equiv\ \{\J^2 (u) : \Hess\, u \geq 0\}
\eqno{(4.10)}
$$

When $X$ carries an almost complex or quaternionic structure as above, we also have the
associated {\bf complex  Monge-Amp\`ere subequation}
$$
 P^\bbc \ \equiv\ \{\J^2 (u) : \Hess^\bbc u \geq 0\}
\eqno{(4.11)}
$$
and   {\bf quaternionic Monge-Amp\`ere subequation}
$$
  P^\bbh \ \equiv\ \{\J^2 (u) : \Hess^\bbh u \geq 0\}
\eqno{(4.12)}
$$
One easily checks that these sets are in fact  subequations and have the convex cone property.  
They provide important examples of {\sl monotonicity cones} for many other subequations
and play an important role in the study of those subequations (See Section 8).
  In addition they typify  an important general construction which we  now  present.

\medskip
\noindent {\bf  4.7. Geometrically defined subequations -- Grassmann structures.}    
  Consider the   Grassmann bundle $\pi:G(p,TX)\to X$ with  fibres  $G(p,T_xX)$,   the set
 of unoriented $p$-planes $\x$
  through the origin in $T_xX$.  If $X$ is a riemannian manifold, then we can identify
   $\x\in G(p,T_xX)$ with orthogonal projection $P_\x \in \Sym(T^*_xX)$ onto $\x$.
   
   Given $A\in \Sym(T^*_xX)$  and  $\x\in G(p,T_xX)$, the {\bf $\x$-trace of } $A$ is defined by 
 $$
   \tr_\x A\ =\ \langle A, P_\x\rangle
   \eqno{(4.13)}
$$
using the natural inner product on $\Sym(T^*_xX)$.  Equivalently, 
 $$
   \tr_\x A\ =\  {\rm trace} \left( A\bigr|_\x\right)
   \eqno{(4.14)}
$$
where $A\bigr|_\x \in \Sym(\x)$ is the restriction of $A$ to $\x$.

Note that for any closed subset $\GG_x \ss G(p,T_xX)$ the set 
$$
F(\GG_x)\ =\ \{ A\in \Sym(T^*_xX) : \tr_\x A \geq 0 \ \ \forall\, \x\in \GG_x\}
$$
is a closed convex cone with vertex at the origin in $\Sym(T^*_xX)$. 
Moreover, $F_x$ automatically satisfies positivity since if $P\geq0$,
 $\tr_\x P\geq 0$ for all $\x\in G(p,T_xX)$.
 
Given  $\GG \ss G(p, TX)$,  define $F(\GG)\ss \Sym(T^*X)$ to be the subset whose 
fibre at $x$ is $F(\GG_x)$.  It is left to the reader to verify that:
$F(\GG)$ is a closed subset if and only if $\pi: \GG\to X$ is a local surjection.
Furthermore, the topological condition (T)(iii), i.e., (4.7), is satisfied in this case. 
It now follows from the discussion in 4.4 (Pure Second Order Subequations) that the 
following holds.

\Prop{4.2}  {\sl Given a closed subset $\GG\ss \Sym(T^*X)$ with $\pi:\GG \to X$ a local surjection,
the subset
$$
\bbr\oplus T^*X \oplus F(\GG)
$$
is a pure second order convex cone subequation.
}
\medskip
The set $F(\GG)$ is said to be {\sl geometrically defined by $\GG$},
and the  $F(\GG)$-subharmonic functions will be referred to as 
{\sl $\GG$-plurisubharmonic functions}.

\medskip
\noindent {\bf 4.8. Gradient independent  subequations.}  
  On a riemannian manifold $X$ there are the subequations which only involve the 
  riemannian hessian and the value of the function.  That is,
with respect to the riemannian splitting (4.3) a subset of the form
$$
F\ =\     T^*X\oplus F',
$$
for a closed subset $F'\ss \bbr\oplus  \Sym(T^*X)$ is said to be {\bf gradient-independent}.
Note that  $F'  \equiv \cn+\cp \ss\bbr\oplus \Sym(T^*X)$ provides an example 
of a gradient-independent subequation.  It is also a convex cone subequation.
 
 Note that a closed subset $F = T^* X\oplus F'$ satisfies (P) and (N) if and only if 
 $$
 F' +\cn+\cp\ \ss\ F'.
 $$
  The topological conditions (T)(i) and (T)(ii)  then  follow from (P) and (N)  by
  using $(r,A) = \lim_{\e\to0}(r-\e, A+\e I)$. This proves that $F\ss J^2(X)$ is a 
 gradient-independent  subequation if and only if  $F$ is of the form $F = T^* X\oplus F'$
 and $F'$ satisfies:
  $$
  F'\ \ {\rm is \ closed},\quad F'+\cn+\cp\ \ss\ F', \ \ \ {\rm and}\ \ \  \Int F_x' \ =\ (\Int F')_x \fa x.
  $$
  Said differently, a closed subset $F'\ss \bbr\oplus \Sym(T^*X)$ with $\Int F_x' = (\Int F')_x$ for all
  $x$,  is a (gradient-independent) subequation if and only if $F'$ is $(\cn+\cp)$-monotone.

\medskip
\noindent {\bf 4.9. Subequations of Reduced Type.}  
  These are the subequations $F\ss J^2(X)$ that do not dependent variable $r$
  in the fibres of $J^2(X)$.  That is, $F=\bbr\times F'$ where $F'$ is a subset
  of the reduced 2-jet bundle $J^2_{\rm red}(X)$,  and we are using the 
  natural bundle  splitting $J^2(X)  =\bbr\times J^2_{\rm red}(X)$.

\vfill\eject
 

\centerline{\headfont \ 5.\   Universal Subequations on Manifolds with Topological G-Structure.}
\medskip

In this section we construct the subequations of principal interest in the paper.
The construction starts with a   ``universal model'',  which is  a euclidean (constant coefficient)
subequation 
 $
 \bbf \ss \bbj^2
 $.
  The idea is to find subequations   $F\ss J^2(X)$  which are locally  modeled on  $\bbf$.
 To accomplish this we recall from (3.6)  the compact  invariance group
 $$
 G \ =\ G( \bbf )\ \equiv\ \{g\in {\rm O}_n : g( \bbf )= \bbf \}
 \eqno{(5.1)}
 $$ 
 where $O_n$ acts on $\bbj^2$ by
 $$
 g(r,p, A) \ \equiv\ (r, gp, gAg^t).
 \eqno{(5.2)}
 $$ 
 The main point is that whenever  a riemannian manifold $X$ 
 is given a topological $G$-structure, we can construct the desired subequation
  by using this structure and the canonical splitting of $J^2(X)$ in the previous section.
Subequations $F$ constructed from an $\bbf$ in this way will be called {\sl riemannian $G$-subequations}.

\medskip
\noindent
{\bf 5.1. Riemannian $G$-manifolds.}
Recall that for a fixed subgroup 
$$
 G\subseteq   {\rm  O}_n,
$$
  a   {\bf  topological $G$-structure on}  $X$ 
is  a family of $C^\infty$ local trivializations
of $TX$ over open sets in a covering  $\{U_\a\}$  of $X$ whose transition functions
have values in $G$, i.e.,  the transition function
from the $\a$-trivialization to the $\b$-trivialization is just a  map 
$
g^{\b,\a}:U_\a\cap U_\b \ \arr\ G.
$
A local trivialization on $U_\a$ is simply a choice of framing
$e^\a_1,...,e^\a_n$ for $TX$ over $U_\a$, and one can think
of these trivializations as a family of {\bf admissible framings}  or {\bf admissible $G$-framings} 
for $TX$.  At each point $x\in X$ the structure determines a family of admissible $G$-frames, in the 
subset of all tangent frames, on which $G$ acts simply and transitively.  This gives a principal
$G$-bundle in the bundle of all tangent frames.

  \Def{5.1}  A {\bf riemannian $G$-manifold} is a riemannian manifold equipped with a 
  topological $G$-structure.
 
\medskip
\noindent
{\bf Important Notes.} (1) A topological $G$-structure  on $X$ (for $G\ss {\rm O}_n$)
determines a second riemannian metric on $X$ by declaring the admissible $G$ frames
to be orthonormal.  These two metrics coincide if the admissible frame fields are orthonormal in
 the original metric, and this will  typically be the case.  However,  the results will hold without this assumption.

(2) Furthermore, one can as well consider topological $G$-structures on $X$ where $G$ is a
subgroup of the {\sl full invariance group}
$$
 \wt G \ =\ \wt G( \bbf )\ \equiv\ \{g\in {\rm GL}_n : g( \bbf )= \bbf \}.
 $$ 
{\sl  The arguments presented here carry through in this more general case}, as the reader can easily verify.
This has genuine applicability since there are interesting subequations which are invariant
under groups like the conformal group, or GL$_n(\bbr)$   or GL$_n(\bbc)$. In these cases
one can do local calculations  with more general frame fields. However, no new subequations
arise from this generalization, and so our exposition is restricted to compact $G$.

\medskip
\noindent
{\bf Important Point:} The compact invariance group $G$ is not the riemannian holonomy group (for either metric).  
No integrability assumption is made. Topological $G$-structures  are soft and easily constructed.

\medskip
\noindent
{\bf 5.2. Riemannian $G$-subequations.}
The idea now is the following.  Suppose $X$ is a riemannian $G$-manifold with admissible 
framings $\{(U_\a, e^\a)\}$. Then each framing $e^\a$ determines a local trivialization of
$J^2(U_\a)$ given by (4.4b).  Now if $G$ lies in the invariance group of a euclidean 
subequation $\bbf$, the local subequations $U_\a\times \bbf$ will be preserved
under changes of framing and thereby determine a global subequation $F$ on $X$.
  We formalize this in the following lemma.

\Lemma{5.2}  {\sl
 Suppose  $\bbf$  is a euclidean subequation with compact invariance group $G$
  and $X$ is a riemannian $G$-manifold.    For $x\in X$, the condition on a 2-jet
  $J\equiv J^2_x u$ that
  $$
\Phi^e(J) \equiv  \left(u(x), e_x(u), (\Hess_x u)(e,e)   \right) \ \in \  \bbf
  \eqno{(5.3)}
  $$ 
  is independent of the choice of $G$-frame $e$ at $x$.
Hence there is a well defined subset $F\ss J^2(X)$ given by}
  $$
  J\in F_x \quad \Leftrightarrow \quad \Phi^e\left(J \right)(x) \in \bbf
  \eqno{(5.4)}
  $$

  \pf
Note that  (5.4) is independent of the $G$-frame $e$, since if  $e'=he$ is another admissible $G$-frame
on a neighborhood of $x$, then 
  $
  \Phi^{e'}(J) = h\left(\Phi^e(J)\right) \in h(\bbf)=\bbf
  $
  by (5.1), (5.2) and the change of frame formula (4.6).  Hence, $F$ is well defined by (5.4).  \qed\medskip

  It is easy to see that 
   \smallskip

 $F$ satisfies (P) \quad $\iff$ \quad $\bbf$ satisfies (P),
 
 $F$ satisfies (N) \quad $\iff$ \quad $\bbf$ satisfies (N), and

 $F$ satisfies (T) \quad $\iff$ \quad $\bbf$ satisfies (T).  
 \smallskip
 \noindent
 so that $F$ is a subequation on $X$.
 
 \Def{5.3}  The subequation $F$ defined by  (5.4) above will be called the {\bf riemannian
 $G$-subequation on $X$ with euclidean model $\bbf$}.
 For simplicity the $F$-subharmonic functions on $X$ will be called $\bbf$-subharmonic.

\medskip
Note that $ \bbf $ is $G$-invariant if and only if $\wt \bbf $ is $G$-invariant.
 It is also easy to show
 \Lemma{5.4} {\sl  The subequation $\ft$ is the  riemannian
 $G$-subequation on $X$ with euclidean model   $\wt \bbf$.}
\medskip

  The extreme cases are where $G=\{I\}$ or $G={\rm O}_n$.  Both are interesting.
  
  \Ex{$G=\{I\}$} Here $\bbf$ can be any euclidean subequation.  A riemannian $\{I\}$-manifold
  is one which is equipped with  $n$ vector fields which are  linearly independent at every point. 
  For example, such fields can be found  on any orientable riemannian 3-manifold.  
  They also exist on any Lie group with a metric
 invariant under left translations by elements of the group. Thus, for an arbitrary euclidean
 subequation $\bbf$ there are plenty of riemannian manifolds which support a  riemannian
 subequation $F$ modeled on $\bbf$.

  \Ex{$G= {\rm O}_n$}  Given a  euclidean subequation $\bbf$ which is invariant under the full orthogonal
  group, there exists a  companion  riemannian O$_n$-subequation $F$ on {\sl every}  riemannian manifold.

    \Ex{$G= {\rm U}_n$}   If $X$ is an almost complex manifold with a compatible (hermitian)
 metric, then one can consider $ \bbf  = \bbr\times\bbr^{2n}\times \bbp_\bbc$
 where $ \bbp_\bbc\ss\Symn$ are the matrices whose hermitian symmetric part is $\geq 0$
 (cf. (4.8), (4.11)  and Section 15). This leads   to solving the the Dirichlet problem for the 
 homogeneous complex Monge-Ampere equation
 on domains in $X$
    
      \Ex{$G=  {\rm Sp}_1\cdot  {\rm Sp}_n$ or ${\rm Sp}_n$} There are analogues of the preceding 
      example on any almost quaternionic or almost
      hyperkahler manifold with a compatible metric, .

\medskip
\noindent {\bf 5.3.  Riemannian $G$-subequations in local coordinates.}  
  Finally we derive the key formula needed to prove a comparison result in Section 10.
Recall that the riemannian  hessian has a simple expression in terms of local coordinates
 $x=(x_1,...,x_n)$ on $X$, namely
 $$
(\Hess u) \left({\partial\over \partial x_i} ,{\partial\over \partial x_j} \right)\ =\ 
{\partial^2 u\over \partial x_i \partial x_j}  - \sum_{k=1}^n \Gamma_{ij}^k (x) {\partial u\over \partial x_k} .
\eqno{(5.5)}
$$
where $\Gamma_{ij}^k$ denote the {\sl Christoffel Symbols}  
of the metric connection defined by the relation
 $
\nabla_{\partial\over \partial x_i} {\partial\over \partial x_j} \ 
=\ \sum_{k=1}^n \Gamma_{ij}^k {\partial\over \partial x_k}
$.
The equation (5.5) can be written more succinctly   as 
$$
(\Hess\, u) \left ({\partial\over \partial x},{\partial\over \partial x} \right)\ =\ D^2 u - \Gamma_x  (Du)
\eqno{(5.5)'}
$$
where $Du$ and $D^2u$ are the first and second derivatives
of $u$ in the coordinates $x$ and
$$
\Gamma_x : \rn \to \Symn 
$$
denotes the linear Christoffel map defined above.

 \Prop {5.5}  {\sl  Let $F$ be a riemannian
 $G$-subequation on $X$ with euclidean model $\bbf$
  on a riemannian $G$-manifold $X$. Suppose $x=(x_1,...,x_n)$ is a local coordinate system on $U$
 and that $e_1,...,e_n$ is an admissible $G$-frame on $U$.  Let $h$ denote the GL$_n$-valued
 function on $U$ defined by $e= h{\partial\over \partial x}$.
 Then a $C^2$-function $u$ is $F$-subharmonic on $U$ if and only if}
 $$
 \left(u,\ hDu, h\left(D^2u - \Gamma(Du)\right)h^t\right) \ \in\  \bbf \qquad{\sl on\ \ } U.
 \eqno{(5.6)}
$$

 \Remark{5.6}  Fix  $J\in J^2_x$ with $x\in U$. \ Then by  using the standard isomorphism
 $
 J_x^2 \cong \bbr\times\rn\times\Symn
 $
 induced by the coordinate system $x=(x_1,...,x_n)$ to represent $J$ as $J=(r,p,A)=(u, Du, D^2u)$,
condition (5.4) can be restated as saying that 
 $$
J\ \in\ F_x \qquad\iff\qquad 
(r,hp, h(A-\Gamma(p))h^t)\ \in\  \bbf 
$$
 
 \pf
 By (5.3) we must show that 
 $$
 \left(u, e(u),  (\Hess\, u)(e.e)   \right)\ =\  \left(u, h Du,  h(D^2 u-\Gamma(Du))h^t    \right).
 \eqno{(5.7)}
$$
The fact that  $e(u) = hDu$, i.e., that
 $
 e_i(u) \ =\ \sum_{j=1}^n h_{ij} {\partial u\over \partial x_j},
 $
is obvious. For the proof that the matrix $(\Hess \, u)(e_i,e_j)$  equals $h(D^2u - \Gamma(Du))h^t$,
note that 
 $$\eqalign
 {
 (\Hess \, u)(e_i,e_j) \ &=\ (\Hess \, u) \left(\sum_{k=1}^n h_{ik} {\partial \over \partial x_k},
 \sum_{\ell=1}^n h_{j\ell} {\partial \over \partial x_\ell} \right) = \ \sum_{k, \ell}^n h_{ik}    (\Hess \, u)   \left( {\partial \over \partial x_k}, {\partial \over \partial x_\ell}\right) h_{j\ell} \cr
  }
 $$
and then apply (5.5).\qed

 \vfill\eject
 

\centerline{\headfont \ 6.\   Jet-Equivalence of Subequations.}
\medskip

In this section we introduce the  important  notion of jet-equivalence for  subequations.
This concept enables us to solve the Dirichlet problem for a very broad spectrum of
subequations on manifolds -- namely those which are  locally   jet-equivalent to a euclidean
one.  They include the subequations considered so far.   In fact, Proposition 5.5 just states
that $F$ is  jet-equivalent to $\bbf$.  However,   jet-equivalence goes far beyond this realm, to
quite general variable coefficient inhomogeneous  equations.  These include the 
basic Calabi-Yau equation (see Example 6.15 below).   It overcomes an important obstacle 
in understanding, on an almost complex manifold,  the {\sl intrinsic homogeneous complex Monge-Amp\`ere subequation} (i.e., intrinsic plurisubharmonic functions) 
where no hermitian metric is used (cf. [HL$_{8}$]).\footnote {*} {Examples show that even on $\bbc^n$
there are hermitian metrics for which not all the classical psh-functions are hermitian psh.}

\medskip
\noindent
{\bf  6.1.  Automorphisms.}  To define jet-equivalence
 we first need to understand  the bundle automorphisms of $J^2(X)$.

 \Def{6.1} (a)\  An {\bf automorphism} of the reduced jet bundle  $J^2_{\rm red}(X)$ is a bundle isomorphism
 $\Phi:J^2_{\rm red}(X)\to  J^2_{\rm red}(X)$ such that with respect to the
 short exact sequence
 $$
 0\ \arr\ \Sym(T^*X)\ \arr\ J^2_{\rm red}(X)\ \arr\ T^*X\ \arr\ 0
 \eqno{(6.1)}
$$
we have 
$$
\Phi(\Sym(T^*X)) = \Sym(T^*X)
\eqno{(6.2)}
$$
 so there is an induced bundle automorphism
$$
g=g_{\Phi}:T^*X \ \arr\ T^*X
\eqno{(6.3)}
$$
and we further require that there exist a second bundle automorphism 
$$
h=h_{\Phi}:T^*X \ \arr\ T^*X
\eqno{(6.3)'}
$$
such that on $\Sym(T^*X)$, $\Phi$ has the form $\Phi(A) = hAh^t$, i.e., 
$$
\Phi(A)(v,w) \ =\ A(h^tv,h^tw)    \qquad {\rm for\ } v,w\in TX.
\eqno{(6.4)}
$$
(b)\ An {\bf automorphism} of the full jet bundle $J^2(X) = \bbr\oplus J^2_{\rm red}(X)$ is a  bundle  isomorphism
 $\Phi:J^2 (X)\to  J^2(X)$ which is  the direct sum of  the identity on $\bbr$ and an automorphism 
 of $J^2_{\rm red}(X)$.

\Lemma{6.2} {\sl  The automorphisms of $J^2(X)$ form a group. They are the sections 
of the bundle of groups whose fibre at $x\in X$ is the group of automorphisms
of $J_x^2(X)$ defined by (6.2), (6.3) and (6.4) above.}

\pf
 We first show that the
 composition of two automorphisms of $J^2(X)$ is again an automorphism.
Suppose $\Psi$ and $\Phi$ are bundle automorphisms. Then $\Psi\circ\Phi$ clearly 
satisfies condition (6.2) and one sees easily that $g_{\Psi\circ\Phi} = g_{\Psi}\circ g_{\Phi}$.
Finally,  $({\Psi\circ\Phi})(A) = \Psi(\Phi(A)) =   h_{\Psi} \Phi(A)h_{\Psi}^t =  
 h_{\Psi}  h_{\Phi} A  h_{\Phi}^th_{\Psi}^t = 
 (h_{\Psi} \circ h_{\Phi} )A (  h_{\Psi} \circ h_{\Phi} )^t =  h_{{\Psi}\circ {\Phi}}  A  h_{{\Psi}\circ {\Phi}}^t$.
The proof that  the inverse of an automorphism is an automorphism  is similar.
\qed

\vfill \eject

\Prop{6.3} {\sl With respect  to any  splitting 
$$
J^2(X) \ =\ \bbr\oplus T^*X \oplus \Sym(T^*X)
$$
of the short exact sequence (6.1), a bundle automorphism has the form
$$
\Phi(r, p,  A) \ =\ (r, gp, hAh^t + L(p))
\eqno{(6.5)}
$$
where  $g, h:T^*X\to T^*X$ are bundle isomorphisms
and $L$ is a smooth section of the bundle 
$\Hom (T^*X, \Sym(T^*X))$. }

\pf With respect to a splitting $J^2(X) \cong T^*X \oplus \Sym(T^* X)$ any bundle isomorphism
is of the form $\Phi(p,A) = (gp+M(A), H(A) + L(p))$. The requirement that $\Phi$ leaves $\Sym(T^*X)$ 
invariant implies that the point-wise linear map $M$ is zero,
and the property that $\Phi(0,A) = (0, hAh^t)$ implies that $H(A)=hAh^t$.
\qed

\Ex{6.4}  Given a local coordinate system $(x_1,...,x_n)$ on an open set
$U\ss X$, the  canonical trivialization 
$$
J^2(U)\ =\ U\times\bbr\times\rn\times\Symn
\eqno{(6.6)}
$$
is determined by $J^2_xu  = (u,Du,D^2u)$ evaluated at $x$, where $Du=(u_1,...,u_n)$ and
$D^2u = ((u_{ij}))$ (cf. Remark 2.8).
With respect to this splitting, every automorphism is of the form
$$
\Phi(u, Du, D^2u) \ =\ (u, \ g Du, \ h\cdot D^2u\cdot h^t + L(Du))
\eqno{(6.7)}
$$
where $g_x, h_x \in \GL_n$ and $L_x : \rn\to \Symn$ is linear for each point $x\in U$.

\Ex{6.5} The trivial 2-jet bundle on $\rn$ has fibre
$$
 \bbj^2  = \bbr\times\rn\times \Symn.
 $$
with automorphism group
 $$
{\rm Aut}(\bbj^2) \ \equiv\  {\rm GL}_n\times {\rm GL}_n\times \Hom(\rn, \Symn)
 $$
 where the  action is  given by
 $$
 \Phi_{(g,h,L)}(r,p,A) \ =\ (r, \ gp, \ hAh^t +L(p)).
 $$
Note that  the group law is
 $$
(\bar g, \bar h, \bar L)\cdot  (g, h,L)\ =\ (\bar g g,\  \bar h L \bar h^t + \bar L\circ g)
 $$

\Remark{6.6}  Automorphisms at a point, with $g=h$,  appear naturally when one considers the action of diffeomorphisms.
 Namely, if $\vf$ is a diffeomorphism
fixing a point $x_0$, then in local coordinates (as in Example 6.4 above) the right action on $J_{x_0}^2$, induced
by the pull-back $\vf^*$ on 2-jets,  is given by (6.7) where $g_{x_0}= h_{x_0}$ is the transpose on the
Jacobian matrix $(({\partial \vf^i\over \partial x_j}))$
and $L_{x_0}(Du) = \sum_{k=1}^n u_k {\partial^2 \vf^k\over \partial x_i\partial x_j}(x_0)$.

\medskip
\noindent
{\bf Cautionary Note 1.}  Despite the remark above, automorphisms of the 2-jet bundle $J^2(X)$,
even those with $g=h$,  have little
to do with global diffeomorphisms or global changes of coordinates. In fact an automorphism radically restructures
$J^2(X)$ in that the image of an integrable section (one obtained by taking $J^2u$ for a fixed smooth function $u$ on $X$) is essentially never integrable.

\bigskip
\noindent
{\bf 6.2.  Jet-Equivalence of subequations.}

\Def{6.7}  Two subequations $F, F'\ss J^2(X)$ are said to be   {\bf jet-equivalent}  if 
there exists an automorphism  $\Phi:J^2 (X)\to  J^2(X)$ with 
$\Phi(F)=F'$.

\medskip\noindent
{\bf Cautionary Note 2.}  A jet-equivalence  $\Phi: F\to F'$ does not take $F$-subharmonic functions
to $F'$-subharmonic functions.   In fact as mentioned above, for  $u \in C^2$,  
$\Phi(J^2 u)$ is almost never the 2-jet of a 
function. It happens if and only if  $\Phi(J^2 u) = J^2u$.

\medskip
Nevertheless, it is easily checked that if $\Phi$ is an automorphism and $F$  is a subequation, then 
$\Phi(F)$ is also a subequation.

\Def{6.8}  We say that a  subequation $F\ss J^2(X)$  {\bf is locally  jet-equivalent to a euclidean 
subequation} if each point $x$ has a coordinate neighborhood $U$ such that 
$F\bigr|_U$ is jet-equivalent to a euclidean (constant coefficient) subequation $U\times \bbf$ in those
 coordinates.  Such a coordinate chart will be called {\bf distinguished}.

\medskip

Using  Example 6.4  we see that Proposition 5.5 can be restated as follows.
Let $\bbf$ be a euclidean subequation with compact invariance group $G=G(\bbf)$
(see (3.6)).

\Prop {6.9} {\sl
Suppose that $F$ is a riemannian $G$-subequation with euclidean model  
 $\bbf$ on a riemannian $G$-manifold $X$.  Then $F$ is locally jet-equivalent to $\bbf$ on $X$.
}

\Lemma{6.10} {\sl  Suppose $X$ is connected and $F\ss J^2(X)$   is locally jet-equivalent  to a euclidean
subequation on $X$. Then there is a euclidean subequation $\bbf\ss \bbj^2$, unique up to automorphisms,
such that $F$ is jet-equivalent to $U\times \bbf$ on every distinguished coordinate chart.}

\pf 
In the overlap of any two distinguished charts $U_1\cap U_2$ choose a point $x$.  Then the 
local automophisms $\Phi_1$ and $\Phi_2$, restricted to $F_x$, determine an automorphism 
taking  $\bbf_1$ to $\bbf_2$.  Thus the local euclidean subequations on these
charts are all equivalent under automorphisms, and they can be made equal by applying the appropriate
constant automorphism on each chart.\qed

\medskip
\noindent{\bf 6.3.  Affine automorphisms and affine jet-equivalence.}
 The automorphism group ${\rm Aut}(\bbj^2)$ can be naturally extended
by the translations of $\bbj^2$.  Recall that the group
of affine transformations of a vector space $V$ is the product ${\rm Aff}(V) = \GL(V) \times V$
acting on $V$ by $(g,v)(u) = g(v)+u$.  The  group law is $(g,v)\cdot (h,w) = (gh, v+g(w))$.
There is a short exact sequence 
$
0\to V\to {\rm Aff}(V) \harr{\pi} {} \GL(V) \to \{I\}.
$

\Def{6.11}  
The {\bf affine automorphism group} is   the 
inverse image 
$$
{\rm Aut}_{\rm aff}(\bbj^2) =\pi^{-1}({\rm Aut}(\bbj^2) )
$$
of ${\rm Aut}(\bbj^2) \ss \GL(\bbj^2)$ under the surjective group homomorphism
$\pi:{\rm Aff}(\bbj^2) \to \GL(\bbj^2)$.  \medskip

Note that any affine automorphism $\wt\Phi$ can be written in the form
$$
{\wt \Phi}\ =\ \Phi   + J
\eqno{(6.8)}
$$
where $\Phi$ is a (linear) automorphism and $J$ is a section of the bundle $J^2(X)$.

Using the affine automorphism group we expand our notion of jet-equivalence in an important way.
Let $F$ and $F'$ be subequations on a manifold $X$.

\Def{6.12} Two subequations $F, F'\ss J^2(X)$ are said to be   {\bf affinely jet-equivalent} 
if there exists an affine automorphism  $\wt \Phi:J^2 (X)\to  J^2(X)$ with 
$\wt\Phi(F)=F'$.

\Def{6.13}  We say that a  subequation $F\ss J^2(X)$  {\bf is locally  affinely jet-equivalent to a euclidean 
subequation} if each point $x$ has distinguished coordinate neighborhoods $U$ such that 
$F\bigr|_U$ is affinely jet-equivalent to a euclidean (constant coefficient) subequation $U\times \bbf$ in those
distinguished coordinates.
 \medskip

\Lemma{6.14} {\sl Suppose $F$ is a subequation on a coordinate chart $U$ with
a given affine jet-equivalence  to a euclidean    subequation $U\times \bbf$.  Write the affine automorphism
$\wt\Phi$ as in (6.8) above so that
$$
J\in F_x \quad\iff\quad  \Phi_x(J) + J_x\ \in\ \bbf
$$
for $x\in U$.  Then}
$$
J\in \ft_x \quad\iff\quad  \Phi_x(J) - J_x \ \in\ \wt{\bbf}
$$
\pf 
 $
 J\in \ft_x \iff -J\notin \Int F_x \iff \Phi_x(-J) +J_x  \notin \Int\bbf
\iff -\{ \Phi_x(J) -J_x\} \notin  \Int \bbf  \iff \Phi_x(J) -J_x \in \wt\bbf
$. (Recall that  $\Int F_x$ denotes interior with respect to the fibre).
\qed

\Ex{6.15.  (The Calabi-Yau Equation)}  Let $X$ be an almost complex hermitian manifold 
(a Riemannian U$_n$-manifold), and consider the subequation $F\ss J^2(X)$ determined by 
the euclidean subequation:  
$$
A_\bbc +I \ \geq \ 0 \and \det_\bbc\{A_\bbc +I\}\ \geq\ 1.
$$
where $A_\bbc \equiv \half(A-JAJ)$ is the hermitian symmetric part of $A$.
Let $f>0$  be a smooth positive function on $X$ and write $f=h^{-2n}$.
Consider the global affine jet-equivalence of $J^2(X)$ given by
$$
\wt\Phi(r,p,A) \ =\ (r, p, (hI)A(hI)^t + (h^2-1)I)\ =\ (r,p, h^2 A+ (h^2-1)I)
\eqno{(6.9)}
$$
and set $F_f = \wt\Phi^{-1}(F)$.  Then
$$
\eqalign
{
(r,p, A) \in F_f \quad &\iff \quad   h^2(A_\bbc+I) \geq 0 
\ \ {\rm and}\ \   \det_\bbc \{h^2(A_\bbc+I)\} \geq 1  \cr
&\iff \quad   (A_\bbc+I) \geq 0 
\ \ {\rm and}\ \   \det_\bbc \{(A_\bbc+I)\} \geq f  \cr
}
$$
so the $F_f$-harmonic functions are functions $u$ with
$\Hess_\bbc u   + I \geq 0$ (quasi-plurisubharmonic)  and $\det_\bbc\{\Hess_\bbc u   + I \} = f$.
If $X$ is actually a complex manifold of dimension $n$ with K\"ahler form $\o$, this last equation
can be written in the more familiar form
$$
\left( {1\over i}\partial\overline{\partial} u  + \o \right)^n\ =\ f \o^n.
$$
\medskip
One can similarly treat the equation
$$
\left( {1\over i}\partial\overline{\partial} u  + \o \right)^n\ =\ e^u f \o^n.
$$
or the same equation with $e^u$ replaced by any non-decreasing positive  function $F(u)$.

There are many further examples illustrating the flexibility and power of 
using local affine jet-equivalence.  We present some of them in Section 18.

\vfill\eject


\centerline{\headfont \ 7.\   Strictly F-Subharmonic Functions.}
\medskip

Consider a  second order subequation $F\ss J^2(X)$ on a manifold $X$.
 For the more general case of subsets  which are just $\cp$-monotone see Remark 7.8.
Strict subharmonicity for $C^2$-functions is unambiguous.

\Def{7.1} A function $u\in C^2(X)$ is said to be {\bf strictly $F$-subharmonic on $X$} if
$$
J^2_x u  \ \in\ \Int F \fa x\in X.
\eqno{(7.1)}
$$
There is more than one way to extend this notion to functions $u\in F(X)$ which are
not $C^2$.  The   choice made in Definition 7.4  will be shown to be useful 
in the next section 8 and in discussing the  Dirichlet Problem.

Fix a metric on the vector bundle $J^2(X)$. (If $X$ is given a riemannian metric, this metric
together with the canonical splitting (4.3) determines a metric on $J^2(X)$.)

\Def{7.2}  
For each $c>0$ we define the {\bf $c$-strict}  subset  $F^c\ss F$ by
$$
 F^c_x \ =\ \{J\in F_x : \dist(J, \sim F_x) \geq c\}
\eqno{(7.2)}
$$
where $\dist$ denotes distance in the fibre $J^2_x(X)$.

\Lemma{7.3} {\sl Suppose that $F$ is a subequation.
Then the  set $F^c$ is closed and is both $\cp$- and $\cn$-monotone.
However, $F^c$ does not necessarily   satisfy (T).}

\pf Note that for $J\in J^2_x(X)$,
$$
J\in  F^c_x \qquad \iff \qquad B_x(J, c)\subseteq F_x
\eqno{(7.3)}
$$
where $B_x(J, c)\ss J^2_x(X)$ is the closed metric ball of radius $c$ about the point $J$ in the fibre.
Suppose now that $J_i \in  F^c_{x_i}$ is a sequence converging  to $J$ at $x$.
Then $B_{x_i}(J_i, c)\subseteq F_{x_i}$ for all $i$, and since $F$ is closed we
conclude that $B_x(J, c)\subseteq F_x$.  Hence $F^c$ is closed.

Suppose now that $P\in\cp_x$. Then by Condition (P) for $F$  we see  that $B_x(J, c)\subseteq F_x
\ \Rightarrow \ B_x(J+P, c) = B_x(J, c) + P     \subseteq F_x$, and so Condition (P) holds for $F^c$.
The proof that $F^c$ satisfies Condition (N) is the same.

\medskip
\noindent
{\bf Example.}  Let $K$ denote the union of the unit disk $\{|z|\leq 1\}$ with the interval $[1,2]$
on the $x$-axis in $\bbr^2$.  Let $K_c = \{ z\in \bbr^2 : \dist(x,K)\leq c\}$.  Define $F$ by requiring
 $p\in K_c$. Then $F^c$ is easily seen to be equal to $K$. Now $F$ is a subset which satisfies
(P), (N) and (T), whereas $F^c$ is a subset which satisfies (P) and (N) but not (T). \qed

\Def{7.4}  A function $u$ on $X$ is {\bf strictly $F$-subharmonic} if for each point
$x\in X$ there is a neighborhood $B$ of $x$ and $c>0$ such that $u$ is $F^c$-subharmonic
on $B$.  Let $F_{\rm strict}(X)$ denote the space of such functions.
\medskip

For a $C^2$ function it is easy to see that the two definitions of strictly $F$-subharmonicity,
given in  Definitions 7.1 and 7.4  agree   (see Remark 7.8).  Moreover, it is easy to see that Definition 7.4  is independent 
of the choice of metric on the bundle $J^2(X)$.

Strictly $F$-subharmonic functions are stable under smooth perturbations.

\Lemma{7.5. (Stability)}  {\sl Suppose $u\in F^c(X)$ and $\psi\in C^2(X)$.  For each precompact open 
subset $Y\ss\ss X$
$$
u+\d \psi \ \in \ F^{c\over 2}(Y) \qquad{\rm if\ } \d \ {\rm is\ sufficiently\ small}.
$$
}

\pf  By Definition 2.2 it will suffice to show that for all $\d>0$ sufficiently small
$$
F^c_x+\d J^2_x \psi \ \ss \ F^{c\over 2}_x \  \fa x\in Y.
\eqno{(7.4)}
$$
Choose $\d>0$ so that
$$
\d \| J^2_x \psi\|\ <\ \smfrac c 2 \  \fa x\in Y.
\eqno{(7.5)}
$$
Fix  $J\in F_x^c$ with $x\in Y$. By  (7.3) we have $B_x(J,c) \ss F_x$. Note then that 
$B_x(J+\d J^2_x \psi, {c\over 2}) = B_x(J, {c\over 2})  +\d J^2_x \psi$ is contained in 
$B_x(J,c)\ss F_x$ by (7.5).  Again using (7.3) this shows that  $J+\d J^2_x \psi \in 
F^{c\over 2}_x$.
\qed\medskip

\Cor{7.6}  {\sl Suppose $u\in F_{\rm strict}(X)$ and $\psi\in C^\infty_{\rm cpt}(X)$.  Then }
$$
u+\d \psi \ \in\ F_{\rm strict}(X)  \qquad{\sl if\ } \d \ {\sl is\ sufficiently\ small}.
$$

We shall need the following two properties.

\Lemma{7.7}  {\sl
\smallskip

(i) \ \ \ \ \ $u,v \ \in\ F_{\rm strict}(X) \quad \Rightarrow\quad \max\{u,v\}\ \in\ F_{\rm strict}(X)$
\medskip

(ii)\ \ \ \ \ If  $F$ satisfies the negativity condition (N), then
$$
u \ \in\ F_{\rm strict}(X) \ \ {\rm and}\ \ c>0 \quad \Rightarrow\quad u-c\ \in\ F_{\rm strict}(X)
$$
}
The proof is straightforward and omitted.

\Remark{7.8. ($\cp$-monotone subsets)}  The results of this section remain true for  $\cp$-monotone 
subsets and for $(\cp+\cn)$-monotone subsets
which are not  necessarily subequations (i.e., condition (T) may not be satisfied).
This is important in Section 11, on boundary convexity, where the results are applied 
to a   $\cp$-monotone subset $\Fa$ which is open.
  Everything is straightforward except the 
proof of the assertion
$$
\psi\in C^2(X)\cap F_{\rm strict}(X) \qquad\Rightarrow\qquad J^2_x\psi \in \Int F \fa x\in X.
\eqno{(7.6)}
$$
In general (7.6) cannot be proved by establishing that
$F^c \ss\Int F$.  For example, take $X=\bbr$ and define $F$ by requiring 
$p\geq 0$ if $x\leq0$ and $p\geq 1$ if $x>0$.  Then the point $x=0, p=\half$ belongs to $F^{1\over 2}$ 
but not  $\Int F$.

We prove (7.6) as follows.  Fix $x_0\in X$ and   trivialize $J^2(X) = U\times \bbj^2$ on a neighborhood
$U$ of $x_0$ via a choice of orthonormal frame field. For such a choice the 
fibre metric on $\bbj$ is  constant.   Now $\psi$ is $c$-strict for some $c>0$ in a
smaller neighborhood
 $U$ of $x_0$, i.e., 
$
B(J^2_x\psi, c) \ \ss\ F_x  
$
for $x\in U$. By continuity we have $|J^2_x\psi - J^2_{x_0}\psi| \leq {c\over 2}$ on a neighborhood $V$
of $x_0$,  and therefore
$
B(J^2_{x_0}\psi, \smfrac c 2) \ \ss\ F_x  
$
for $x\in W \equiv U\cap V$. Thus $W\times B(J^2_{x_0}\psi, \smfrac c 2) \ss F$ and in particular
$J_{x_0}^2\psi \in \Int F$.\qed
\medskip

The following elementary example shows that $J^2_x\psi \in \Int(F_x)$ for $x$ near $x_0$
(rather than $J^2_{x_0}\psi \in ( \Int F)_{x_0}$ as in (7.1)) is not sufficient to guarantee that $\psi$ is 
 $c$-strict near $x_0$.  Let $F$ be the one variable 
subequation defined by $\{|p|\leq|x|\} \cup( \{0\}\times [-1,1])$. Take $\psi\equiv0$.  
Then:

1) $\psi$ is not strictly $F$-subharmonic since $J^2_0\psi \notin \Int F$,
 
2) $\psi$ is not $c$-strict near $x=0$ since $F^c_x=\emptyset $ for $0<|x|<c$, but

3) $J^2_x\psi \in \Int(F_x)$ for all $x$.

 \vskip .5in


\centerline{\headfont \  8.\   Comparison Theory -- Local to Global.}
\medskip

In this section we begin our analysis of the uniqueness/comparison question for $F$-harmonic
functions with given boundary values on a compact subset $K$ of  $X$.
Let  us set the notation
$$
F(K) \ \equiv\ \left \{ u\in \USC(K) : u\bigr|_{\Int K} \in F(\Int K)\right\}
$$
Then the  comparison principle can be stated as follows:
$$
\eqalign
{
{\rm for \ all \ }\ u\ \in\ F(K) \quad {\rm and} \quad &-w\ \in\ \ft (K) \cr
u\ \leq\ w\ \ {\rm on\ } \ \partial K\quad\Rightarrow\quad &u \  \leq w\ \ {\rm on\ } \ K
}
\eqno{(8.1)}
$$
However, we prefer to state it in the following form which invokes duality.

\Def{8.1}  We say that {\bf  comparison   holds for $F$ on $X$} if for 
all compact sets $K\ss X$, whenever
$$
u\ \in\ F(K) \and  v\ \in\ \wt F(K),
$$
  the {\bf Zero Maximum Principle} holds for $u+v$ on $K$,
that is,
$$
u+v\ \leq \ 0\quad {\rm on\ \ } \partial K\qquad\Rightarrow \qquad
u+v\ \leq \ 0\quad {\rm on\ \ } K
\eqno{(ZMP)}
$$

\medskip

If comparison holds, it is immediate that\medskip

\noindent {\bf Uniqueness for the Dirichlet Problem holds},
that is:
\medskip
\centerline{If two $F$-harmonic functions agree on   $\partial K$, they must
agree on $K$.}
\medskip

Local comparison does not imply global comparison.  However, for a weakened form
of comparison, local does imply global.  This is the main result of this section.

\Def{8.2}  We say that   {\bf weak comparison   holds for $F$ on $X$} if for all compact subsets $K\ss X$, and functions
$$
u\in F^c(K), \quad  v\in \wt F(K), \ \ c>0
$$
   the Zero Maximum Principle (ZMP) holds for $u+v$ on $K$.
We say that {\bf local weak comparison holds for $F$ on $X$} if for all
$x\in X$, there exists a neighborhood $U$ of $x$ such that weak comparison holds
for $F$ on $U$.

\Theorem{8.3} {\sl Suppose that $F$ is a  subequation on a manifold $X$.
If local weak comparison holds for $F$ on $X$, then weak comparison
holds for $F$ on $X$.}

\pf
Suppose weak comparison fails for $F$ on $X$.
Then there exist $c>0$, $u\in F^c(X)$, $v\in \wt F(X)$ and compact $K\ss X$
with 
$$
u+v\ \leq \ 0 \quad {\rm on}\ \ \partial K\quad{\rm but}\quad \sup_K(u+v) \ >\ 0.
$$
Choose a maximum point $x_0\in \Int K$ and let $M=u(x_0)+v(x_0) >0$
denote the maximum value. Fix  local coordinates $x$ on a neighborhood
of $x_0$ containing $U=\{x: |x-x_0|<\rho\}$. By choosing  $\rho$ sufficiently small
we can assume that  weak comparison  holds for $F$ on $U$.  

By the Stability Lemma 7.5, we have $u' = u-\delta |x-x_0|^2 \in F^{c\over 2}(U)$
if $\delta$ is chosen small enough. Now $x_0$ is the unique maximum point
for $u'+v$ on $U$, with maximum value $M$.  Choose $K_0  =
\{x : |x-x_0| \leq \rho/2\} \ss U$.  Then $\sup_{\partial K_0}(u'+v) = M'<M$.
Set $u'' = u' - \max\{0, M'\}$.  Then by Condition (N) for $F^{c\over 2}$ 
(see Lemma 7.3)  we have 
$u''\in F^{c\over 2}(U)$. Moreover, 
$$
\sup_{\partial K_0} (u''+v)\ \leq \ 0\quad{\rm while}\quad u''(x_0)+v(x_0) \ =\ M-\max\{0, M'\} \ >\ 0.
$$
Thus, weak comparison fails for $F$ on $U$, contrary to assumption.\qed
\medskip
Weak comparison can be strengthened as follows.
Define
$$
F_{\rm strict}(K) \ =\ \left\{ u \in \USC(K) : u\bigr|_{\Int K} \in F_{\rm strict}(\Int K) \right\}
\eqno{(8.2)}
$$

\Lemma {8.4} {\sl Suppose weak comparison  holds for $F$ on $X$.
For all compact sets $K\ss X$, if $u\in  F_{\rm strict}(K)$ and $v\in\ft(K)$, then $u+v$ satisfies the (ZMP).}

\pf
We assume $u+v\leq 0$ on $\partial K$.  Since $u\in\USC(K)$, for each $\d>0$ the set
$$
U_\d \ =\ \{x\in K : u(x)+v(x)<\d\}  
$$
is an open neighborhood of $\partial K$ in $K$.  Exhaust $\Int K$ by compact sets $K_\e$
with $\Int K = \bigcup_\e \Int K_\e$.
Then $\partial K_\e \ss U_\d$ for $\e>0$ small.  Now $u-\d+v \leq 0$ on $\partial K_\e$.
However, $u-\d $ is $F^c$-subharmonic on $\Int K_\e$ for some $c>0$.  Hence, (WC) for $K_\e$ states that $u-\d+v \leq 0$ on $K_\e$.  Thus $u - \d+v \leq 0$ on $\Int K$. Hence $u+v\leq0$ on $K$.\qed

 \vskip .5in


\centerline{\headfont \  9.\   Strict Approximation and  Monotonicity Subequations.}
\medskip

In this section we discuss certain global approximation techniques which 
can be used to deduce comparison from weak comparison.
Consider a general second order subequation $F$ on a manifold $X$.

\Def{9.1}  We say that {\bf strict approximation } holds for $F$ on $X$ if for each compact
set $K\ss X$, each function $u\in F(X)$ can be uniformly approximated by functions in 
$F_{\rm strict}(K)$.\medskip

When strict approximation is available, it is easy to show that  weak comparison implies comparison.

\Theorem{9.2. (Global Comparison)}  {\sl  Suppose $F$ is a subequation
 on a manifold $X$.
Assume that both  local weak comparison  
and the strict approximation property hold for $F$ on $X$.
Then comparison   holds for $F$ on $X$.}

\pf 
Suppose $u\in F(K)$,  $v\in \wt F(K)$, and 
$u+v\leq 0$ on $\partial K$.
By Theorem 8.3 we can assume that weak comparison holds. 
By strict approximation, for each $\e>0$, there exists   
$\overline u \in F_{\rm strict}(K)$ with $u-\e \leq \overline u\leq u+\e$ on $K$.
Now Property (N) implies that $\overline u - \e \in F_{\rm strict}(K)$.
Note that $\overline u - \e +v \leq u+v \leq0$ on $\partial K$.
Lemma 8.4 states that   $\overline u - \e +v \leq 0$ on $K$. This
proves that $u+v$ satisfies the (ZMP).\qed

\medskip

\medskip\noindent
{\bf Example 9.3. (The Eikonal Equation)}.   This subequation $\bbf$ on $\rn$ is defined
by $|\nabla u|\leq 1$.  Given $u\in \bbf(K)$,  set $u_\e = (1-\e)u$.
Then $u_\e \in F^\e(K)$ because if $\vf $ is a test function for $(1-\e)u$ at $x_0$, then 
${1\over 1-\e} \vf$ is a test function for $u$ at $x_0$.  Thus $|\nabla \vf(x_0)|\leq 1-\e$.
\medskip

In contrast to this example, for the geometric subequations $F$ that are of primary interest in this paper,
the approximations $u_\e$ will be of the form
$$
u_\e\ =\ u+\e\psi\qquad 0< \e\leq\e_0
\eqno{(9.1)}
$$
where $\psi$ is a $C^2$-function independent of $u$.
The function $\psi$ will be referred to as an {\bf approximator} for $F$.

Suppose $M\ss J^2(X)$ is a subset such that the fibre-wise sum
$$
F+M\ \ss\ F\and \e M\ \ss\ M \qquad {\rm for}\ 0< \e\leq\e_0.
\eqno{(9.2)}
$$
If $\psi$ is a $C^2$-function which is strictly $M$-subharmonic, then 
$\psi$ is an approximator for $F$ (see the proof of Theorem 9.5).  However,
condition (9.2) implies that at each point $x$
$$
F_x +\a_1 J_1 +\a_2J_2 \ \ss  \ F_x \fa \a_1>0, \a_2>0\ \ {\rm and\ \ } J_1,J_2\in M_x.
\eqno{(9.3)}
$$
Hence we might as well assume that $M$ is a convex cone, but at this point it is convenient not 
to assume that $\overline{\bf M}$ is a subequation, i.e., that ${\bf M}$ satisfies (P) and (N).
See Remark 9.11.

 \Def{9.4} A subset $M\ss J^2(X)$ will be called {\bf a convex monotonicity cone for } $F$ if
\medskip
(1)\ \ $M$ is a  convex cone with vertex at the origin, and
\medskip
(2)\ \ \ $F+M\ \ss\ F$.
\medskip

\Theorem{9.5}  {\sl Suppose $M$ is a convex monotonicity cone for $F$ as above. 
If there exists $\psi \in C^2(X)$ which is strictly $M$-subharmonic, then strict approximation holds for 
$F$ on $X$.}
 
 \pf
 It will suffice to establish the following.
 
 \medskip
 \noindent
 {\bf Assertion 9.6.}
  {\sl For each compact subset $K\ss X$, there exists $\d>0$ such that}
 $$
 u+\e \psi\ \in \ F^{\e \d}(K) \fa u\in F(K)\quad{\rm and\ for\ all\ \ } \e>0.
 $$
To begin note   that
 $u-\vf$ has  local maximum  0 at a point $x$ if and only if $u+\e\psi - (\vf+\e\psi)$
 has   local maximum 0 at $x$. 
 Hence we must show that  under the hypothesis 
 $J^2_x\vf \in F_x$  we have that $J^2_x (\vf +\e\psi) = J^2_x \vf + \e J^2_x \psi \in F^{ \e \d}_x$,  in other words that  $F_x+\e J^2_x\psi\ss F^{\e \d}_x$.

  Since $\psi\in C^2(X)$, $J^2_x\psi$ is a continuous function of $x$. Hence $\{J^2_x\psi : x\in K\}$
 is compact.  That $\psi$ is strictly $M$-subharmonic implies that $\{J^2_x\psi : x\in K\}$ is a compact
 subset of $\Int M$.  Take $\d=$ the distance from $\{J^2_x\psi : x\in K\}$
 to $\sim \Int M$.  Then $B(J^2_x\psi, \d)\ss M_x$ for all $x\in K$ where $B$ denotes the ball in the fibre.
 Suppose now that  $J\in F_x$ and $x\in K$.  Then 
 $$
 B(J+\e J^2_x\psi, \e \d) \ =\  J +\e B(J^2_x\psi, \d)\ \ss \  F_x+M_x\ \ss\  F_x
 $$
 as desired.
 \qed
 \medskip

Combining Theorem 9.2 with Theorem 9.5 yields the version of Global Comparison
that will be used in this paper.

\Theorem{9.7}  {\sl Suppose $F$ is a subequation   on a manifold $X$.
Assume that $X$ supports a $C^2$ function which is strictly $M$-subharmonic, where $M$
is a monotonicity cone for $F$.  Then local weak comparison for $F$ implies
global comparison for $F$ on $X$.}

\Remark {9.8. (Circular Monotonicity Cones)} In a situation where the $C^2$-function $\psi$ is given,
the simplest monotonicity cone to consider is one whose fibre at each point $x$ is a circular cone
$C(J)$ about $J\equiv J^2_x\psi$.
If $\d = \dist( J^2_x\psi, \sim \Int M)$, as in the above proof, then the circular cone
can be taken to be the cone $C^\d(J)$
on the ball $B(J,\d)$.  The cross-section of this cone $C^\d(J)$
 by the hyperplane (through $J$)  perpendicular to $J$, is a ball  of radius $R$ in this hyperplane,
 where, setting $\g = 1/R$, one calculates that
 $$
 \d\ =\ {|J|\over \sqrt{1+ \g^2|J|^2}}
 $$
This same cone  will  be denoted by $C_\g(J)$ when $\g$ is to be emphasized.

\Lemma {9.9}  {\sl Suppose that $F$ is a subequation with $F_x\neq \emptyset$ and
$F_x\neq J^2_x(X)$, and fix $J\in J^2_x(X)$.  The following are equivalent.
\medskip
\item
{(1)}\ \ $F_x$ is $C_\g(J)$-monotone.\medskip

\item
{(2)}\ \  The boundary  $\partial F_x$ can be graphed over the hyperplane $J^\perp$ with graphing function
$f$ which is $\g$-Lipschitz, \ \  i.e.,  for $ J_0 \in J^\perp $
$$
J_0+tJ \in F_x \ \ \iff \ \ t\geq f(J_0),
$$ 

and  for all $\bar J_0, J_0 \in J^\perp$}
$$
-\g |J_0|\ \leq \ f(\bar J_0 + J_0)\ \leq \ \g|J_0|
$$
The elementary proof is left to the reader.

\Remark {9.10}  Suppose $\bbf$ is a universal model (i.e., a constant coefficient subequation)
which is $G$-invariant.  If $\bbf$ has a convex monotonicity cone ${\bf M}$ which is $G$-invariant,
then on a manifold $X$ with topological $G$-structure, the induced subequation
$F$ and  the induced convex monotonicity cone 
$M$ satisfy $F+M\ss F$, i.e., $M$ is a convex  monotonicity cone  for $F$ on $X$.

\Remark {9.11}   The global comparison Theorems 9.2 and 9.7 are useful even locally.
In the next Section 10 we will prove that weak comparison holds for any constant coefficient
subequation $\bbf$ on $\rn$.  Consequently, if $\bbf$ has a monotonicity cone ${\bf M}$
with non-empty interior,
then local comparison holds for $\bbf$.  To prove this, fix a point $x_0$ and  pick a point $(r,p,A) \in\Int {\bf M}$
above  $x_0$.
Let $\psi$ denote the quadratic function whose 2-jet at $x_0$ equals $(r,p,A)$.  Then $\psi$ is strictly
${\bf M}$-subharmonic in a neighborhood of $x_0$.
\medskip

Therefore, given a constant coefficient subequation $\bbf$,   the key question is:
$$
{\rm When\ does\ } \bbf\ {\rm have\  a \ monotonicity\  cone \ } {\bf M} \ {\rm with\ interior?}
\eqno{(9.4)}
$$
We might as well assume that ${\bf M} = \overline{\Int {\bf M}}$, and (as noted
prior to  Definition 9.4) that ${\bf M}$ is convex.

 Now ${\bf M}$ need not satisfy
conditions (P) or (N), i.e. ${\bf M}$ need not be a subequation.  (For a useful example let ${\bf M}= C(J)$
in Remark 9.8.)  However, such an ${\bf M}$ can always be enlarged to 
${\bf M}'= {\bf M}+(\bbr_-\times \{0\}\times \cp)$ which satisfies:
\medskip

(1) \ \ $ {\bf M}' = \overline{\Int  {\bf M}'}$ and ${\bf M}'$ is convex (as is true of $ {\bf M})$,
\medskip

(2) \ \  ${\bf M}'$ is both $\cn$-  and $\cp$-monotone, i.e., ${\bf M}'$ is a subequation,
\medskip

(3)\ \ $\bbf + {\bf M}' \ \ss\ \bbf$,   i.e., ${\bf M}'$ is a also monotonicity cone for $F$.

\medskip
\noindent
This proves that  the question (9.4) is equivalent to the question:
$$
{\rm When\ does\ } \bbf\ {\rm have\  a \  convex\ conical \ monotonicity\  subequation \ } {\bf M}?
\eqno{(9.5)}
$$
Each such ${\bf M}$ is a convex cone with interior, and $M$  contains  $\bbr_-\times \{0\}\times \cp$.  
Hence $M$
can be thought of as a fattening of $\bbr_-\times \{0\}\times \cp$ to a convex  set with interior.
The bigger   ${\bf M}$ is,  the more can be concluded about $\bbf$.  However, as illustrated by the example
${\bf M} \equiv C(J)$, (9.4) may not be easier to answer since ${\bf M}'$  cannot be described directly
without defining ${\bf M} = C(J)$ and then defining ${\bf M}'$ as the sum ${\bf M}+(\bbr_-\times \{0\}\times \cp)$.

  In the next subsecton we present a few of the basic examples of monotonicity subequations $M$
  on a manifold.  They are all based on a universal euclidean  model ${\bf M}$ satisfying Definition 9.4.
  We leave it to the reader to explicitly describe the model ${\bf M}$ in the examples.

 \bigskip
 
 \centerline{\bf Examples of  Strict Approximation Using Monotonicity.}
 \medskip

Most of the following examples are purely second order.   We use  the canonical splitting
  from Section 4.2:
$$
J^2(X)\ \cong \ \bbr\oplus T^*X\oplus \Sym(T^*X).
\eqno{(9.6)}
$$
 on the riemannian manifold $X$.
  They are the   subequations $F\ss J^2(X)$
which are  the pull-backs of subsets $F'\ss \Sym(T^*X)$ under the projection
 $J^2(X)\to \Sym(T^*X)$ induced by (9.6).  (See Subsection 4.4.)
 
 \Ex{9.12. (The Real Monge-Amp\`ere Monotonicity Subequation)} 
 For all  pure second order subsets $F$ the Positivity Condition (P) is equivalent to $F$ being  
  $P$-monotone where 
$$
P \ =\ \bbr\oplus T^*X\oplus \cp
$$
is the Monge-Amp\'ere subequation discussed in Subsection 4.6.
For a general riemannian manifold $X$, strict $P$-subharmonicity is the same as strict convexity.

\Theorem{9.13}  {\sl Suppose $X$ is a riemannian manifold which supports a 
strictly convex $C^2$  function. Then strict approximation holds for every
pure second order subequation $F$ on $X$.}

\medskip

For example in $\rn$ the function $|x|^2$ is strictly $P$-subharmonic.
More generally if $X$ has sectional curvature $\leq 0$,  the function $\delta(x) = \dist(x,x_0)^2$ is strictly 
$P$-convex up to the first cut point of $x_0$. In particular, if $X$ is complete and simply connected, then $\d$ is globally strictly $P$-convex.  Of course on any riemannian manifold $X$ the function $\d$ is strictly $P$-subharmonic
 in a neighborhood of $x_0$ since $\Hess_{x_0} \d= 2I$.
Thus strict approximation holds for all
pure second order subequations in these cases.

  \Ex{9.14.  (The Complex Monge-Amp\`ere Monotonicity Subequation)} 
 Suppose now that $(X,J)$ is a hermitian almost complex manifold and consider the projection
 $$
 J^2(X) \ \to\ \Sym_\bbc(T^*X)
 $$
 given by projecting onto $\Sym(T^*X)$ and then taking the Hermitian symmetric part (cf. 4.5). 
A subequation defined by pulling back a subset of $\Sym_\bbc(T^*X)$ will be called a {\bf complex hessian} subequation.  Each such  $F$ is
  $P^\bbc$-monotone where  $P^\bbc$ is the 
 complex Monge-Amp\`ere subequation defined in 4.6.
 
 The $C^2$ functions $\psi$ on $X$ which are strictly $P^\bbc$-subharmonic are just the classical
 $\bbc$-plurisubharmonic functions if $X$ is a complex manifold.  We will use this terminology
 even if $J$ is not integrable.
 
 \Theorem{9.15}  {\sl  Suppose $F$ is a complex hessian subequation on  a hermitian almost complex manifold $X$.  If $X$ supports a strictly $\bbc$-plurisubharmonic function of class $C^2$, then strict approximation holds for $F$ on $X$.}
 \medskip
 
  \Ex{9.16. (The Quaternionic Monge-Amp\`ere Monotonicity Subequation)} 
  We leave it to the reader to formulate the analogous result on  an almost quaternionic hermitian
 manifold  which supports a  strictly  $\bbh$-plurisubharmonic function.

 Note that $P\ss P^\bbc\ss P^\bbh$ so the corresponding subharmonic functions
 are progressively easier to find.

  \Ex{9.17. (Geometrically Defined Monotonicity Subequations)} 
  Suppose $M^{\GG}$ is geometrically defined by a subset ${\GG}$ of the Grassmann
  bundle $G(p, TX)$ as in section 4.7.  
  For example the three Monge-Amp\`ere cases over $K=\bbr, \bbc$ and $\bbh$ 
  are geometrically defined
  by taking ${\GG}$ to be  the  Grassmannian of $K$-lines $G(1,TX)$, 
  $G_\bbc(1,TX) \ss G(2,TX)$ and $G_\bbh(1,TX) \ss G(4,TX)$ respectively.
 Other important examples are given by taking ${\GG}$ to be all of $G(p,TX)$,
 $G_\bbc(p,TX)$ and $G_\bbh(p,TX)$ for general $p$.
 
 Any calibration $\phi$ of degree $p$ on a riemannian manifold $X$ determines
 the subset
 $$
 {\GG}(\phi)\ =\ \{\x\in G(p,TX) : \phi(\x)=1\}
 $$
 of calibrated $p$-planes, and hence a convex cone subequation   
 geometrically defined by ${\GG}(\phi)$.  That is, a $C^2$ function $u$ is 
 {\bf ${\GG}(\phi)$-plurisubharmonic}  if  
 $$
 \tr_\x \Hess\, u \ \geq\ 0  \fa \x\in {\GG}(\phi).
 $$
  \Theorem{9.18}  {\sl  Suppose ${\GG}$ is a  subset of the Grassmann bundle
  $G(p,TX)$ as in Section 4.7.  If $X$ supports a $C^2$ strictly ${\GG}$-plurisubharmonic function,
  then strict approximation holds for any subequation 
 on $X$ which is $M^{\GG}$ monotone.
    }
  \medskip
  Of course, one such equation is $M^{\GG}$ itself.

 \Ex{9.19. (Gradient Independent Subequations)}
These are   subsets $F\ss J^2(X)$
which are  the pull-backs of subsets $F'\ss \bbr\oplus \Sym(T^*X)$ under the projection
 $J^2(X)\to \bbr\oplus \Sym(T^*X)$ induced by (9.4). See Section 4.8.
 All such sets are $\cm^-$-monotone where
 $$
 \cm^- \ \equiv\ \bbr^-\oplus T^*X \oplus \cp.
 \eqno{(9.7)}
 $$
 if and only if they satisfy (P) and (N). Hence, gradient independent
  subequations are automatically $\cm^-$-monotone.
 
 More generally, if $M'$ is any one of the pure second
 order monotonicity subequations discussed above, then
 $$
 M\ =\ \bbr^-\oplus T^*X\oplus M'
 $$
 is a gradient independent monotonicity subequation.  Each subequation $F$ which
 is $M$ monotone must be gradient independent.
 The same $\psi$'s used in the previous examples will work for $F$.
 This is because on a compact subset $K$, $\psi-c$ is strictly 
 negative for $c>>0$.  More precisely, Theorem 9.15
 (and its complex and quaternionic versions) continues to hold for the more general gradient independent subequations.

\vfill\eject 


\centerline{\headfont \ 10.\    A   Comparison Theorem for G-Universal Subequations.}
\medskip

The main result of this section can be stated as follows.

\Theorem{10.1}  {\sl Suppose $F$ is a subequation on a manifold $X$ which is locally affinely jet-equivalent
to a constant coefficient subequation $\bbf$.  Then weak comparison holds for $F$ on $X$.}
\medskip

A case of particular interest in this paper is the following.
Suppose  $X$ is a riemannian manifold equipped  with a topological
$G$-structure, for a  closed subgroup $G\ss {\rm O}_n$.  Suppose
$
\bbf \ \ss \  \bbj^2 
$
is a euclidean  $G$-invariant  subequation on $\rn$.
Let $F\ss J^2(X)$ denote riemannian $G$-subequation with model fibre $\bbf$ on $X$
given by Lemma  5.2.

\Cor{10.2}  {\sl     
Suppose $F$ is a riemannian $G$-subequation on $X$.  Then  weak comparison holds
 for $F$-subharmonic functions on X.}
 
 \pf  On any coordinate chart $U$ with an admissible 
 local framing, $F$ is jet-equivalent to $\bbf$ (by Proposition 5.5).  \qed

\medskip\noindent
{\bf Remark.} The following  proof of local weak comparison uses the
Theorem on Sums, which is discussed in Appendix C.
The argument for equations which are honestly constant   coefficient  is particularly easy -- see Corollary C.3.

 \medskip
 \noindent
 {\bf Proof of Theorem 10.1.}
 For clarity we first present the proof in the case where $F$ is locally (linearly) jet-equivalent
 to a constant coefficient equation as in Definition 6.7.
 
 By Theorem 8.3 we need only prove weak comparison  on a chart $U$ where $F$
is jet-equivalent to $\bbf$.
Suppose    weak comparison fails on $U$.
We will derive a contradiction using the Theorem  on Sums.
Failure of comparison means that there exist
$u\in F^c(U)$ and $v\in \wt F(U)$ such that $u+v$ does not satisfy the 
Zero Maximum Principle on some compact subset $K\ss U$.
Let $h(x)$ and $L_x$  determine the field of automorphisms taking $F$ to $\bbf$.
  Theorem C.1 says that there exist a point $x_0\in \Int K$, a sequence of numbers $\e\searrow 0$
with associated points $z_\e = (x_\e, y_\e) \to (x_0, x_0)$, and
2-jets:
$$
\a_\e\ \equiv\ (x_\e, r_\e, p_\e, A_\e)\ \in \ F_{x_\e}^c
\and
\b_\e\ \equiv\  (y_\e, s_\e, q_\e, B_\e)\ \in\  \ft_{y_\e}
$$
 with the following properties.
$$
 r_\e \  =\ u(x_\e),\qquad 
 s_\e\ =\ v(y_\e),  \and r_\e+s_\e = M_\e \  \searrow \  M_0\ >\ 0
 \eqno{(10.1)}
 $$  
$$
 p_\e\ =\ {x_\e-y_\e\over \e}\ =\ -q_\e \and {|x_\e-y_\e|^2\over \e} \ \arr\ 0
 \eqno{(10.2)}
 $$  
 $$  \left( \matrix{
 A_\e & 0 \cr 0 & B_\e\cr}  \right)
 \ \leq\ {3\over\e} \left( \matrix{
 I & -I \cr -I & I\cr}  \right).
\eqno{(10.3)}
$$
The jet-equivalence of $F$ and $\bbf$ says that
$$
\a_\e'\ \equiv\ (r_\e', p_\e', A_\e')\ \in \ \bbf^c
\and
\b_\e'\ \equiv\  (s_\e', q_\e', B_\e')\ \in\  \wt \bbf
\eqno{(10.4)}
$$
where
$$
r_\e' \ =\ r_\e,\quad p_\e' \ =\ g(x_\e) p_\e, \quad A_\e' \ =\ h(x_\e)  A_\e h(x_\e)^t + L_{x_\e}(p_\e)
\eqno{(10.5)}
$$
$$
s_\e' \ =\ s_\e,\quad q_\e' \ =\ g(y_\e) q_\e, \quad B_\e'\ =\  h(y_\e) B_\e h(y_\e)^t  +  L_{y_\e}(q_\e). 
\eqno{(10.6)}
$$

Since $\a_\e' \in \bbf^c$,  we also have   that
$\a_\e'' \equiv (r_\e-M_\e, p_\e', A_\e'+P_\e)  \in   \bbf^c$ for any  $P_\e\geq 0$.

By (10.4) we have $-\b_\e' \notin \Int\bbf$.  Hence
$$
0\ <\ c\ \leq\ \dist(\a_\e'', -\b_\e')\ =\ |\a_\e''+\b_\e'|.
\eqno{(10.7)}
$$
To complete the proof we show that $\a_\e''+\b_\e'$ converges to zero.
The first component of $\a_\e''+\b_\e'$ is $r_\e-M_\e+s_\e$ which tends to zero by (10.1).

\noindent
The second component of $\a_\e''+\b_\e'$ is 
$$
p_\e'+q_\e'  \ = \ g(x_\e){(x_\e-y_\e)\over \e}
-g(y_\e){(x_\e-y_\e)\over \e} \ =\ \biggl (g(x_\e)-g(y_\e)\biggr){(x_\e-y_\e)\over \e}
$$
 which converges to zero as $\e\to 0$ by (10.2).

It remains to find $P_\e\geq0$ so that the third component of $\a_\e''+\b_\e'$,
namely $A_\e'+P_\e+B_\e'$, converges to zero.  This will contradict (10.7).

 Multiplying both sides in (10.3) by
$$
\left( \matrix { h(x_\e) & 0 \cr 0 & h(y_\e)\cr    }  \right) \ \ {\rm on\ the\ left\ and\ \ }
\left( \matrix { h(x_\e)^t & 0 \cr 0 & h(y_\e)^t\cr    }  \right) \ \ {\rm on\ the\ right}
$$
gives
$$
\left( \matrix { h(x_\e) A_\e  h(x_\e)^t & 0 \cr 0 & h(y_\e) B_\e h(y_\e)^t\cr    }  \right)   
 \ \leq\ 
 {3\over\e} \left( \matrix{  h(x_\e)h(x_\e)^t & -h(x_\e)h(y_\e)^t \cr
  -h(y_\e) h(x_\e)^t & h(y_\e) h(y_\e)^t  \cr}  \right).             
$$
Restricting these two quadratic forms to  diagonal elements $(x,x)$ then yields
$$
\eqalign
{
h(x_\e) A_\e  h(x_\e)^t + h(y_\e) B_\e h(y_\e)^t  \ &\leq\ 
 {3\over\e}\left[  h(x_\e)( h(x_\e)^t -  h(y_\e)^t)  -     h(y_\e)( h(x_\e)^t -  h(y_\e)^t)   \right]  \cr
&=\ {3\over\e}  ( h(x_\e) -  h(y_\e)) ( h(x_\e)^t -  h(y_\e)^t)   \cr
&\leq \   {\l \over\e} \left| x_\e-y_\e\right|^2\cdot I \qquad{\rm for\ some\ \ } \l>0.
}
$$
Thus there exists $P_\e\in \bbp$  so that 
$$
h(x_\e) A_\e  h(x_\e)^t + h(y_\e) B_\e h(y_\e)^t  + P_\e \ =\ 
 {\l \over\e} \left| x_\e-y_\e\right|^2\cdot I.
$$
It now follows from   the definitions in (10.5) and (10.6) that
$$\eqalign
{
 A_\e' + B_\e'   + P_\e \ &=\  {\l \over\e} \left| x_\e-y_\e\right|^2\cdot I 
+       L_{x_\e}(p_\e)   
+       L_{y_\e}(q_\e).    \cr
  }
\eqno{(10.8)}$$
However,
$$
\eqalign
{
 \left| L_{x_\e}(p_\e) + L_{y_\e}(q_\e) \right|
&=\    \left| ( L_{x_\e} -  L_{y_\e}) \left (  {x_\e-y_\e\over \e}  \right)  \right| \cr  
  & \leq \| L_{x_\e} -  L_{y_\e}\|  {|x_\e-y_\e|  \over \e}  \cr  
 =\ & \ \ O\left( {|x_\e-y_\e|^2  \over \e} \right)
}
$$
Using (10.2) this shows that 
$$
A_\e' + B_\e'   + P_\e \ \cong \ {|x_\e-y_\e|^2  \over \e}\ \to\ 0\quad{\rm as\ } \e\searrow 0.
\eqno{(10.9)}$$
This completes the proof in the case of linear jet-equivalence.

Suppose now that our local jet-equivalence is affine and can be written in the
form $\wt{\Phi}_x = \Phi_x +J_0(x)$ where $\Phi$ is a linear jet-equivalence as above (cf. (6.8)).
Then the proof above goes through essentially unchanged except that, in light of Lemma 6.9,
we must  replace (10.4) with 
$$
\eqalign
{
\a_\e' + J_0(x_\e)\ &\equiv\ (r_\e', p_\e', A_\e') + J_0(x_\e)\ \in \ \bbf^c  \cr
{\rm and}\quad   \b_\e' - J_0(y_\e)\  &\equiv\  (s_\e', q_\e', B_\e')   - J_0(y_\e) \ \in\  \wt \bbf
}
\eqno{(10.4)'}
$$
We now observe that $J_0(x_\e)-J_0(y_\e) \to 0$,  and the rest of the proof is exactly
as written above.
\qed

 \medskip

Combining  Corollary 10.2   with Theorem 9.7 yields comparison
for a wide class of riemannian $G$-subequations.

\Theorem{10.3} {\sl Suppose $F$ is a  riemannian $G$-subequation  on a manifold $X$. 
 If $X$ supports a $C^2$ strictly $M$-subharmonic
function, where $M$ is a monotonicity cone  for  $F$, then comparison holds
for $F$ on  $X$.}

 \vskip .3in


\centerline{\headfont \ 11.\   Strictly F-Convex   Boundaries and Barriers.}\medskip

In this section we introduce the notion of  a strictly $F$-convex boundary 
for a general subequation $F$. 
This notion implies  the existence of barriers,    which are crucial for our main results. 
Strictly $F$-convex boundaries have a relatively simple geometric characterization in terms of the 
second fundamental form of $\bo$ (see subsection 11.4). 
Thus even for very general equations, the geometry of the 
boundary   enters explicity into the study of the Dirichlet Problem.
 
 The idea of an $F$-convex boundary is the following.  Start with a general subequation $F$.
 For each $\l\in\bbr$, there is a reduced subequation $F_\l$ defined by setting the $r$-variable equal
 to $\l$  (e.g. $\D u \geq e^u$ becomes $\D u\geq e^\l$).  
 To each $F_\l$ we associate an {\sl asymptotic interior} $\overrightarrow{F_\l}$.
 This is an open, point-wise conical set obtained by taking rays with conical neighborhoods 
 that eventually lie in $F$.  In many examples (such as $\D u \geq e^u$) the asymptotic interiors
 $\overrightarrow{F_\l}$ all agree.  Of course this includes the case where the subequation is independent
 of $r$ to start with.  Furthermore, in this case, if $F$ is itself a cone, then $\Fa=\Int F$. In general
 these asymptotic interiors  are  simpler than $F$,  and as examples
 will show, the process of computing them is often quite straightforward.
 
 Now the boundary of a domain $\bo$  is called {\sl locally $\overrightarrow{F_\l}$-convex}   if it has 
 a  local defining function  which is  strictly $\overrightarrow{F_\l}$-subharmonic; and it  is 
said to be {\sl locally $F$-convex}  if it is locally $\overrightarrow{F_\l}$-convex for all $\l$.
 The main point of these definitions  is that $\overrightarrow{F_\l}$-convexity at a point $x$ implies the existence of a ``$\l$-barrier'' for $F$
 at $x$, and these barriers are exactly what is needed for existence of solutions to the Dirichlet
 problem.
 
 One might ask about the simpler concept where one just takes the asymptotic interior of 
 $F$ directly.  This can be done, but for many important equations involving the 
 dependent variable (e.g. the Calabi-Yau-type equations in  Example 6.15), this 
 na\"ive $F$-convexity rules out all boundaries, whereas the refined one gives just the
 right condition (see Section 19). Morever, if $F$ is independent of the $r$-variable and
 has constant coefficients, one can show that the existence of barriers is actually equivalent
 to our notion of boundary convexity.

 We note that for  constant coefficient,  pure second-order subequation on $\rn$,  
 convexity and the existence of barriers were analyzed  in [HL$_4$].
 
 We also note that under   a mild hypothesis on $F$,   the boundaries of small enough balls  in any coordinate system are strictly $F$-convex (see Proposition 11.9 below).

\vskip .3in


\noindent
{\bf 11.1  Asymptotics.}  
Consider the canonical decomposition 
 $J^2(X) =\bbr\oplus J^2_{\rm  red}(X)$ with fibre coordinates $J\equiv (r, J_0)$, and refer to a subequation
 of the form $\bbr\oplus F$ with $F\ss J^2_{\rm  red}(X)$ as a subequation {\bf of reduced type / independent 
 of the $r$ variable} (cf. Subsection 4.9).  These subequations must be handled first.

 Our first step is to replace $F$ by an (asymptotically) smaller open set $\Fa$  which is a cone
 (i.e., each fibre is a cone with vertex the origin).

 \Def{11.1}   
 Suppose that $F\ss J^2_{\rm  red}(X)$  is a subequation independent of the $r$-variable.
 The {\bf asymptotic interior} $\Fa$ of $F$ is the set of all $J\in J^2_{\rm  red}(X)$ 
 for which there exists a neighborhood $\cn(J)$ in (the total space of) $J^2_{\rm  red}(X)$
 and a number $t_0> 0$ such that 
 $$
 t\cdot \cn(J)\ \ss\ F\fa t\geq t_0
 \eqno{(11.1)}
 $$
Note that:
$$
{\rm If\ } F \ {\rm is\ a\ cone, \ then\ }  \Fa\ =\ \Int F,
 \eqno{(11.2)}
 $$
but otherwise $\Fa$ is smaller than $F$ asymptotically,  and may be empty.

\Prop{11.2}  {\sl The asymptotic interior $\Fa$ is an open cone in $J^2_{\rm red}(X)$ 
which satisfies Condition (P).}

\pf  Obviously $\Fa$ is  open and is a  cone (i.e. $t\Fa=\Fa$ for all $t>0$).
To prove (P) suppose $J\in \Fa$ and $P_{x_0}\in  \cp_{x_0}\ss J^2_{x_0}(X)$.
Extend $P_{x_0}$ to a smooth section $P$ of $\cp$ near $x_0$.
Note that the fibre-wise sum  $\cn(J)+P$ is a neighborhood $\cn(J+ P_{x_0})$ of $J+P_{x_0}$.
(The fibre at $x$ of the sum is defined to be empty if $\cn(J)_x$ is empty.)
Finally note that $t\cn(J+P_{x_0})_x = t\cn(J)_x + tP_x\ss F$ if $t\geq t_0$,
since $ t\cn(J)_x \ss F$  if $t\geq t_0$.   
\qed

\Lemma{11.3}  {\sl  Suppose $F\ss J^2_{\rm red}(X)$ is a reduced  subequation.  If
$F$ is  defined by a $G$-invariant 
universal model $\bbf$ as in Lemma 5.2,   then the asymptotic interior $\Fa$ 
of $F$ is the bundle of open cones induced by the universal model $\oa \bbf$  via (5.3), 
where $\oa \bbf$ is the  asymptotic interior  of $\bbf$. 
More generally,  jet-equivalence preserves asymptotic interiors.}

\pf Exercise.
\medskip

\Def{11.4} A $C^2$-function $u$ with $J^2_x u\in\Fa$ for all $x$ will be called {\bf strictly $\Fa$
subharmonic}.  

\medskip

This is the notion required for boundary convexity and barriers for the Dirichlet problem.
In particular, the question of when the closure of $\Fa$ is a subequation having $\Fa$ as its interior,
can be avoided.

We next observe that  for many subequations $F$ there are many  local $\Fa$-subharmonic functions.
Suppose we are in a local coordinate system for $X$ with standard fibre coordinates
$(p,A)$   for $J^2_{\rm red}(X)$.

\Prop{11.5}  {\sl 
If $\Fa_{x_0}$ is non-empty, then  for all $\e>0$, $\Fa$ contains a subset of the form
$$
B(x_0,r_0)\times B(p_0, \d_0)\times (\e I+ \Int \cp).
$$
Moreover, if $\Fa_{x_0}$ contains a point of the form $(0,A)$, then the function
$
\rho\ \equiv\ |x-x_0|^2 - R^2
$
 is strictly $\Fa$-subharmonic on $B(x_0,R)$ for $R>0$ sufficiently small.}

\pf
If $(p,A)\in \Fa_{x_0}$ then by positivity we may replace $A$ by $\l I$  (since
$\l I-A >0$ for   large $\l$). By the openness of $\Fa$ there exist  $r_0, \d>0$ such
that $B(x_0,r_0)\times B(p, \d) \times \{\l I\} \ss \Fa$,  and therefore by positivity
$B(x_0,r_0)\times B(p, \d) \times (\l I + \cp) \ss\Fa$.  Applying  the cone property proves the first
assertion.
For the second, note that if $p=0$ then
$B(x_0,r_0)\times B(0, \d_0) \times ( I + \cp) \ss\Fa$.
\qed

Without a point $(0,A_0)\in\Fa_{x_0}$ the second assertion can fail.  Consider the subequation
$A\geq p\geq 0$ on $\bbr$, or the subequation $A\geq p_0 I\geq 0$ on $\rn$.

\vskip .3in


\noindent
{\bf 11.2  Boundary Convexity.}  
Suppose  $\O$ is a domain with smooth boundary $\bo$  in $X$. (or some open subset which we also call $X$.)   \
By a {\bf  defining  function}  for $\bo$ we mean a smooth function $\rho$ defined on a neighborhood of
$\bo$ such that 
\smallskip
\centerline{$\bo  = \{x:\rho(x)\ = 0\}, \qquad d\rho \neq 0$ \ on \ $\bo$,  \and $\rho  < 0$ \ on \ $\O$.}\medskip
\smallskip
\noindent
Note that for $x\in\bo$ we have $J^2_x\rho = \{0\}\times J^2_{{\rm red},x} \rho$. 
For simplicity of notation we set $J^2_x\rho =  J^2_{{\rm red},x} \rho$.

\Prop{11.6}   Suppose  $F$ is a reduced subequation on $X$ with asymptotic interior $\Fa$.
Let  $\O\ss X$ be a domain with smooth boundary.  Then for $x\in \bo$ the following are 
equivalent.
\medskip
\item{(1)}
  There exists a local defining function $\rho$ for $\bo$ near $x$ such that
  $$
  J^2_x\rho\in \oa{F_x}
 \eqno{(11.3)}
 $$
  (so that $\rho$ is strictly $\Fa$-subharmonic near $x$).\medskip
  
  \item{(2)}
  There exists a local defining function $\rho$ for $\bo$ near $x$ and $t_0 > 0$
  such that
   $$
  J^2_x\rho  +  t(d\rho)_x \circ  (d\rho)_x       \in \oa{F_x}  \fa t\geq t_0.
 \eqno{(11.4)}
 $$

  \medskip
  \item{(3)}
Given any defining function $\rho$ for $\bo$ near $x$, there exists $t_0 > 0$ such that 
$$
J^2_x \rho + t(d\rho)_x \circ (d\rho)_x \ \in\ \Fa_x \fa t\geq t_0.
\eqno{(11.4)}
 $$
Note that if  (11.4)  holds for $t=t_0$, then it holds for all $t\geq t_0$ by positivity for $\Fa$.

\Def{11.7.(Boundary Convexity)}  Let $F$ be a reduced subequation  on $X$ and $\O\ss X$ a
smoothly bounded domain.  Then $\bo$ is {\bf strictly $F$-convex at $x\in \bo$} if the 
equivalent conditions (1), (2) and (3) hold at $x$.  The boundary $\bo$ is said to 
be {\bf strictly $F$-convex} if this holds at every point $x\in \bo$.

\medskip
\noindent
{\bf Proof of Proposition 11.6.}   {\bf (1) $\Rightarrow$ (3)}:
Assume that (1) is true for the defining function $\rho$.
Any other local defining  function for $\bo$ is of the form $\wt \rho = u \rho$ for some smooth
function $u>0$. At $x\in \bo$ we have  $d \wt \rho = u d \rho$ and 
$J^2\wt  \rho = u J^2  \rho + du\circ d\rho$.  
Hence with $\e>0$ we have
$$
\eqalign
{
J^2\wt  \rho + t d \wt  \rho \circ d \wt  \rho \ 
&= \ u J^2 \rho + du\circ d\rho +t u^2  d\rho \circ d\rho  \cr
&= \ u \left ( J^2 \rho  - \e \cdot I \right)     
+ \left( t u^2  d\rho \circ d\rho + du\circ d\rho + u\e \cdot I\right)            
}
$$
at $x$.
Now for $\e>0$ sufficiently small we have $u \left ( J^2 \rho  - \e \cdot I \right)   
\in    \Fa$ by the assumption on $\rho$, the openness of $\Fa$
 and the cone property  for $\Fa$. 
For all $t$ sufficiently large, the remaining term in the second line above lies in
$\Int \cp_x \ss \Sym(T_x^*X)$. Now apply the positivity condition (P) for $\Fa$.

 {\bf (2) $\Rightarrow$ (1)}: Assume that (2) is true for the defining function 
 $\rho$   for $\bo$ and consider $\rho_t \equiv \rho + \half t\rho^2$
for $t\geq 0$.  Then $\rho_t$ is also a defining function for $\bo$ and it has 2-jet
$$
J^2 \rho_t \ =\ J^2 \rho + t d\rho \circ d\rho  \qquad {\rm on\ \ }\bo.
$$
Hence $\rho_t$ satisfies (1) if $t\geq t_0$.

Evidently (3) $\Rightarrow$ (2), and (1) $\Rightarrow$ (2) by positivity.
\qed
\medskip

From Proposition 11.6(3) we conclude the following.

\Cor{11.8}  {\sl  Let $\O \ss \ss X$ be a pre-compact  domain with a smooth
 strictly $F$-convex boundary, and 
$\rho$ a global defining function for $\bo$.
Then for all $t>0$ sufficiently large, 
the defining function $\rho_t =\rho+\half t \rho^2$
is strictly $F$-subharmonic in a neighborhood of $\bo$.}

\pf 
Since $\rho+\half t  \rho^2$ is strictly $\Fa$-subharmonic in a neighborhood  $U_x$
of $x$ in $\bo$ for $t\geq$ some $ t_x$, we can pass  to a finite covering of  
$\bo$ by such sets $U_x$ and take $t$ to 
be the largest $t_x$ in the family.  Then $\rho_t$ is strictly  $F$-subharmonic
at all points of $\bo$ and hence in a neighborhood of $\bo$. 
\qed
\medskip

The question of a global strictly $F$-subharmonic defining function on a neighborhood 
of $\ob$ will be discussed later.

A different question, namely,  the local existence of strictly $F$-convex domains, has already
been answered in Proposition 11.5.  We restate it here in the language of $F$-convexity.

\Prop{11.9} {\sl  Let  $F$ be a reduced subequation on $X$ and fix $q\in X$.
Suppose the fibre $F_{q}$ contains a critical jet $(0,A)$.  Then in any local coordinate system,
the balls $\{x : \|x-x(q)\|< R\}$ are all strictly $F$-convex for all $R>0$ sufficiently small.
}

\medskip
\noindent
{\bf Example.}  The strictness condition: $J^2_x \rho \in \Fa$ in (1)
 is not independent of the defining function $\rho$.  Suppose
$F=\wt \bbp$, the constant coefficient \seq defined by requiring at least one eigenvalue
of $D^2u$ be $\geq 0$, and note that $F=\Fa$.  The    defining function $\rho(x) = 1-|x|^2$ for the unit sphere $S^{n-1}$
(with inwardly pointing gradient) is not   
$\wt \bbp$-subharmonic, whereas the defining function $\rho + \rho^2$
{\sl is} strictly $\wt \bbp$-subharmonic near $S^{n-1}$. 
\vskip.3in

We now complete Definition 11.7 by extending    boundary convexity from 
reduced subequations to the general case.

Associated to a general subequation $F\ss J^2(X)$ is a family of reduced subequations
$F_\l\ss J^2_{\rm red}(X)$, $\l\in \bbr$, obtained by setting the $r$-variable equal to 
the constant $\l$.  Said differently, $F_\l$ is defined by 
$$
\{\l\}\times F_\l \ = \ F\cap \biggl \{   \{ \l\}\times J^2_{\rm red}(X)\biggr\}.
$$
Each reduced subequation $F_\l$ has an asymptotic interior $\overrightarrow{F_\l}$. 
(It is not uncommon for all the asymptotic interiors $\overrightarrow{F_\l}$ to agree even though
$F$ in not a reduced subequation.)

\Def{11.10.(Boundary Convexity Continued)}  Given a general subequation $F\ss J^2(X)$ and
and a domain $\O\ss X$ with smooth boundary, we say that $\bo$ is {\bf strictly $F$-convex
at a point} $x$ if $\bo$ is strictly $F_\l$-convex at $x$ for each $\l\in \bbr$ (cf. Definition 11.7).
The boundary  $\bo$ is called (globally) {\bf $F$-convex} if it is $F$-convex at every $x\in \bo$.

\vskip .3in


\noindent
{\bf 11.3.  Barriers.} 
The existence of  barriers at a boundary point $x_0\in \bo$ of a domain $\O$
is needed for establishing boundary regularity for the Dirichlet problem.

\Def{11.11}  Given $\l\in\bbr$, $x_0\in\bo$,  and a local defining function $\rho$   for $\bo$ near $x_0$, 
we say that $\rho$ {\bf  defines a  $\l$-barrier for $F$  at $x_0 \in \bo$} if  
 there exist 
$C_0 > 0$, $\e>0$  and $r_0>0$ such that the function
$$
\b(x) \ =\ \l + C\left(  \rho(x) -   \e   {|x-x_0|^2\over 2}\right)
\eqno{(11.5)}
$$
is strictly $F$-subharmonic on $B(x_0,r_0)$ for all   
$C\geq C_0$.  If $F$ is a reduced subequation, then we say  {\bf $\rho$ defines a barrier for $F$ at $x_0$},
since the same $\rho$ works for all $\l  \in \bbr$.

\Theorem {11.12. (Existence of Barriers)} {\sl Suppose $\O\ss X$ is a domain 
with smooth boundary $\bo$
which is strictly $F$-convex at $x_0\in\bo$.
 Then for each $\l\in \bbr$ there exists a local defining function
$\rho$ for $\bo$ near $x_0$ which defines a $\l$-barrier for  $F$  at $x_0$.
}
\pf
Choose a local defining function $\rho$ for $\bo$ near $x_0$ with
$J \equiv(D_{x_0}\rho, D^2_{x_0}\rho)\in \Fa_{\l'}$ ($\l'>\l)$).
Then by the definition of $\Fa$ there exists a neighborhood $\cn(J)$
in the total reduced 2-jet space $J^2_{\rm red}(X)$, and a number  
  $C_0 > 0$ such that 
$$
C\cdot \cn(J) \ \ss \ \Int F_{\l'} \quad \ {\rm  for \ all }\ \ C\geq C_0.
$$

By shrinking we may assume that 
$\cn(J) =  B(x_0, r_0)\times B(D_{x_0}\rho, \d_0)\times B(D^2_{x_0}\rho, \e_0)$ 
is a product neighborhood in some local coordinate system.  

Set $\a(x) \equiv \rho(x) - \e {|x-x_0|^2\over 2}$ so our desired barrier will be $\b(x) = \l + C\a(x)$.  Then
$$
J_{{\rm red}, x} \a  \ =\ \left ( D_x\rho - \e(x-x_0),\, D^2_x \rho - \e I\right).
\eqno{(11.6)}
$$
By choosing $\e$ and $r_0$ sufficiently small we have
$$
J_{{\rm red}, x} \a \ \in\ \cn(J)\qquad {\rm if\ \ } x\in B(x_0, r_0).
$$
This proves that
$$
J_{{\rm red}, x} \b \ =\  CJ_{{\rm red}, x} \a\ \in\ \Int F_{\l'}
$$
if $C\geq C_0$ and $ x\in B(x_0, r_0)$.

Now by the negativity condition (N) we have that $F_{\mu}\supset F_{\l'}$ if $\mu<\l'$. 
Hence, 
$$
(-\infty, \l')\times \Int  F_{\l'}  \ \ss\ F
$$
is an open subset in the full jet space $J^2(X)$, and it contains
$J_x\b$ for all $C\geq C_0$ and all $x\in B(x_0, r_0)$. \qed

\Remark{11.13. (Existence of Barriers -- Refinements)}  Actually, less than is stated in Definition 11.10
 is required for the application
to the Dirichlet problem.  Namely, one only needs that for each $\l\in \bbr$, there is a continuous defining function 
$\rho$ such that 
there exist $C_0\geq0$, $\e>0$ and $r_0>0$ so that  
$$
\b(x) \ =\ \l + C\left(  \rho(x) -   \e   {|x-x_0|^2\over 2}\right) \ \in\ F_{\rm strict}\left(\overline{B(x_0,r_0)}\cap \ob   \right)
\eqno{(11.7)}
$$
is strictly $F$-subharmonic on $\overline{B(x_0,r_0)}\cap \ob$ for all $C\geq C_0$.
(See \S 7 for the definition of $F_{\rm strict}(K)$.)
\medskip

This applies to the intersection $\O=\O_1\cap\O_2$ where
$\O_1, \O_2 \ss X$ are two domains with smooth, strictly
 $\Fa$-convex boundaries. Fix $x_0\in \bo_1\cap\bo_2$ and $\l\in\bbr$.
 Assume that there exist local defining functions $\rho_k$ for each $\bo_k$ near $x_0$
 and constants $r_0>0$, $\e>0$ and $C_0\geq 0$ such that
 $\b_k(x) \equiv \l + C  \left(\rho_k (x) - \e|x-x_0|^2\right)
  \in\ F_{\rm strict}\left(\overline{B(x_0,r_0)}\cap \ob   \right)$  for all $C\geq C_0$.
  Then
 $$
 \max\{\b_1, \b_2\} \ \in\  \ F_{\rm strict}(\overline{B(x_0,r_0)}\cap \ob_1\cap\ob_2)
 $$
 for all $C\geq C_0$  by Lemma 7.7(i).\medskip

 This remark is particularly useful in establishing existence for parabolic equations.

  Our notion of strict boundary convexity for a subequation $F$ is sufficient to insure the existence
of a family of $F$-barriers (Theorem 11.12).  Might there be other different conditions which also 
insure the existence of barriers?  For reduced euclidean subequations   the answer is essentially no.  
 
\Prop{11.14}
  {\sl   Let  $F$ be a  reduced constant coefficient subequation on  $\rn$  
  and consider a domain $\O\ss\ss \rn$.  Suppose 
$\rho$ is a defining function for $\bo$ at $x_0$ which  
 defines a barrier for $F$   at $x_0$.  Then $J^2_{x_0}\rho \in \overrightarrow
{F_{x_0}}$, i.e., $\bo$ is strictly $F$-convex at $x_0$.
}
\medskip
\noindent
{\bf Proof (Outline).}  By hypothesis we have $C(D_x\rho-\e(x-x_0), \, D^2_x\rho-\e I) \in \Int F$
if $|x-x_0|$and $\e$ are small and $C$ is large.  
One must show that $(D_x\rho-\e(x-x_0), \, D^2_x\rho-\e I)$
fills out a neighborhood $\cn(J)$ of $J=(D_{x_0}\rho,\, D^2_{x_0}\rho)$ in the total space
$J^2_{\rm red}(X)$.  The only difficulty is in the $p$-variable.  
The map $x\mapsto D_x\rho-\e(x-x_0)$ has derivative at $x_0$ equal to $D^2_{x_0} -\e I$
which is invertible for almost all  $\e$.  Hence this map fills out a ball about $p_0=D_{x_0}\rho$.
\qed

 \medskip

Finally we point out that the asymptotic interior $\Fa$ of a reduced subequation
$F\ss J^2_{\rm red}(X)$  could have been defined using another section
$J_1$ of $J^2_{\rm red}(X)$ as the vertex rather than the zero section.
Define the 
 {\bf asymptotic interior of $F$ based at $J_1$}, denoted $\overrightarrow{F_{J_1}}$,
 to be the set of  $J$ for which there exists a neighborhood $\cn(J)$ of $J$ in $J^2_{\rm red}(X)$
and $t_0   >  0$ 
such that  $J_1 + t\cn(J) \ \ss\ F$ for all  $t\geq t_0$.  Then
$$
 \overrightarrow{F_{J_1}} \ =\ \overrightarrow{F_{J_2}} \quad {\rm  for\ 
  any\  two\  sections \ \ } J_1\ \ {\rm and\ \ }   J_2.
\eqno{(11.8)}
$$

The proof  is left to the reader.

\vskip .3in


\noindent
{\bf 11.4.  A Geometric Characterization of  Boundary Convexity.} 
In this subsection we see that on a riemannian manifold $X$, Proposition  11.6 enables 
us to characterize the $\Fa$-convexity of a boundary $\bo$  in terms of 
 its  second fundamental form $II_{\bo}$ with respect to  the outward-pointing normal $n$. 
  Recall the canonical decomposition  of the 2-jet bundle given in Section 4.2:
  $$
  J^2(X) = \bbr\oplus T^*X\oplus
 \Sym(T^*X).
  \eqno{(11.9)}
 $$
 
 \Prop {11.15}  {\sl  The boundary $\bo$ is strictly $F$-convex at a point $x\in \bo$
 if and only if   
 $$  
 \left( 0, n, \,  tP_n \oplus II_{\bo}\right) \ \in\ \oa{{F_x}} \ \ \ 
 \fa t  \geq {\rm some\ } t_0.
 \eqno{(11.10)}
 $$
 where $P_n$ denotes orthogonal projection onto the normal line $\bbr\cdot n$ at $x$.}
\medskip\noindent
{\bf Note.} 
Blocking   with respect to the decomposition
 $T_xX = \bbr\cdot n \oplus T_x(\bo)$, (11.10) can be rewritten
  $$  
 \left( 0, (1,0), \left(\matrix {t & 0 \cr 0 & II_{\bo}\cr}\right)\right) \ \in\      \oa{{F_x}} \ \ \ \ 
 \fa t \geq {\rm some\ } t_0.
 \eqno{(11.10)'}
 $$
\pf 
Choose $\rho$ to be the signed distance function  on a neighborhood $U$ of $\bo$. That is,
$$
\rho(y)  \   =\ \cases{-\dist(y,\bo)\quad {\rm if\ \ } y\in U^-  \cr
+\dist(y,\bo)\quad {\rm if\ \ } y\in U^+  \cr}.
$$
where $U-\bo= U^+\cup U^-$ and signs are chosen so that $\nabla \rho = n$ on the boundary
$\bo$. Then it is a standard calculation (cf. [HL$_2$, (5.7)]) that 
 $$
 \Hess_x \rho \ =\ \left(  \matrix{
 0&0\cr
 0& II_{\bo}
 }
 \right).
 $$
 with respect to the splitting $T_x X =(\bbr \cdot \n \rho) \oplus T_x \bo$.
 Since $\nabla \rho = n$ at $x$ the assertion follows directly from 
   Proposition 11.6 and Definition 11.7. \qed

\vskip .3in


\noindent
{\bf 11.5. Example: Equations Involving Curvatures of the Graph.} The condition of $\Fa$-convexity is nicely illustrated by the following 
class of equations.  Fix a closed subset $\bbs\ss \rn$, invariant under permutation of coordinates
and satisfying $S+(\bbr^+)^n \ss S$. Consider a function $u:\O\to \bbr$ on a smoothly
bounded domain $\O\ss\rn$.  At each point $x\in \O$ let $\kappa(x) 
= (\kappa_1(x),...,\kappa_n(x))$ be the principal curvatures of the graph of $u$
in $\bbr^{n+1}$.  We define our subequation $F$ by the condition that $\kappa(x)\in \bbs$
for all $x$.  The simplest case, $\bbs = \{\sum_j \kappa_j\geq 0\}$ corresponds to the graph
having non-negative mean curvature, and the resulting equation on $u$ is the classical 
minimal surface equation.

This subequation determined by $\bbs$ 
involves only first and second derivatives $(Du, D^2 u)=(p,A)$ and can be defined
as follows. Consider $E\in\Symn$ given by
$$
E\ \equiv \ I - \smfrac 1 {\nu (1+\nu)} p\circ p \qquad {\rm where\ \ } \nu = \sqrt{1+|p|^2}
$$
Then the principal curvatures $\kappa_1,...,\kappa_n$ of the graph of $u$
are the eigenvalues of the linear transformation
$\smfrac 1 \nu EAE$.
Thus we have that
$$
F\ =\ \left \{  (p,A) : {\rm the\ eigenvalues\ of \ } \smfrac 1 \nu EAE \ \ {\rm lie\ in\ } \bbs\right\}.
$$
Straightforward computation shows that 
$$
\Fa \ \supseteq \ \left \{ (p,A) :  {\rm the\ eigenvalues\ of \ } P^\perp AP^\perp  \ \ {\rm lie\ in\ } \Int\, \bbs\right\}.
$$
where $P^\perp$ is orthogonal projection onto the hyperplane $p^\perp$. Thus we find the following.

\Prop{11.17} {\sl The boundary $\bo$ is strictly $\Fa$-convex if its own principal curvatures 
$\kappa_1^{\bo},...,   \kappa_{n-1}^{\bo}$,    satisfy}
$$
(0, \kappa_1^{\bo},...,   \kappa_{n-1}^{\bo}) \in \Int\, \bbs\qquad {\rm at\ each\ } x\in \bo.
$$

Thus for the minimal surface equation we obtain the classical condition that the
 boundary have positive mean curvature at each point. 
 More generally, let $\bbs_k$ be the closure of the component  of
 $ \{\s_1(\kappa)>0,...,\s_k(\kappa)>0\}$ containing $(1,...,1)$, and let $F_k$ be 
  the corresponding subequation.  Then we get the condition 
  that $\s_1(\kappa)>0,...,\s_k(\kappa)>0$ for the principal curvatures of $\bo$.
  Similarly one could restrict further to the set $\s_k/\s_\ell >1$ (or $\s_\ell/\s_k >1$) as in [LE],
  and the corresponding condition comes out for the boundary. Of course Proposition 11.17  applies to 
  very general sets $\bbs$.
 
 Existence for the Dirichlet problem for these examples is the topic of section 17.

\vfill\eject


\centerline{\headfont \   12.\   The Dirichlet Problem -- Existence.}
\medskip

Throughout this section we assume that $F$ is a  subequation on a manifold $X$ 
and that  $\O\ss\ss  X$ is a domain  with smooth boundary $\bo$. 
Furthermore, we assume that both $F_{\rm strict}(\ob)$ and $\ft_{\rm strict}(\ob)$ contain at least one function bounded below.  (See  (8.2) for this notation.)  This assumption is ``minor'', for example it is obvious locally
for all subequations with $F_x \neq\emptyset $ and $\ft_x\neq\emptyset$ for all $x$.
Our key assumption in the existence theorems is that $\bo$ is both $F$ and $\ft$ strictly convex.

\Def{12.1} Given a boundary function $\vf\in C(\bo)$, consider the {\sl Perron family}
$$
\cf(\vf) \ \equiv\   \left  \{ u\in \USC(\ob) : u\bigr|_\O \in F(\O) \ \ 
{\rm and \  \      u\bigr|_{\bo} \leq \vf  \ \ }      \right\}
$$
and define the {\sl Perron function} 
$$
U(x) \ \equiv\ \sup \{ u(x) : u\in \cf(\vf)\}.
$$
to be the upper envelope of the Perron family.

\medskip

We shall begin by isolating all the conclusions that hold only under the assumption of weak 
comparison.  This is done in the next theorem
and its corollary.  In the two subsequent  theorems the remaining gap in existence
is filled in two different ways. In the first we see that, assuming comparison, the gap is easy to fill.  
In the second we assume constant coefficients and apply  an argument of Walsh [W]  to fill the gap.

The method of proof for the next theorem is the classical barrier argument.  In the case
where $F=\cp_\bbc$ on $\bbc^n$ these arguments can be found as far back as 
Bremermann [B] except for the ``bump argument''  for part (3), given in the proof of 
Lemma  $\wt {\bf F}$ below, which is  due to Bedford and Taylor
[BT]. This argument was rediscovered by Ishii [I].

\Theorem{12.2}  {\sl  Suppose that $\bo$ is both $F$ and $\ft$ strictly convex,
and that weak comparison holds for $F$ on $X$.
Given $\vf\in C(\bo)$, the Perron function $U$ satisfies:
\medskip

\hskip .2in (1)\ \ $U_*\ =\ U\ =\ U^*\ =\ \vf$ on $\bo$,\smallskip

\hskip .2in (2)\ \ $U\ =\ U^*$ is $F$-subharmonic on $\O$, \smallskip

\hskip .2in (3)\ \ $-U_*$ is $\wt F$-subharmonic on $\O$.  \smallskip

}

\Cor{12.3}  {\sl If $U$ is lower semicontinuous on $\O$, i.e., if $U_*=U$, then: \medskip

\hskip .2in  (a)\ \ $U\ \in\ C(\ob)$,\medskip

\hskip .2in  (b)\ \ $U$  is $F$ harmonic on $\O$,\medskip

\hskip .2in  (c)\ \ $U\ = \ \vf$ on $\bo$.\medskip
\noindent
i.e., $U$ solves the Dirichlet problem on $\ob$ for boundary values $\vf$.
}

\Theorem{12.4}  {\sl Assume that comparison holds for $F$, and suppose that
$\bo$ is both $F$ and $\ft$ strictly convex. Then for each 
$\vf\in C(\bo)$ the Perron function $U$ solves the Dirichlet problem on $\O$ 
for boundary values $\vf$.
}

\Theorem{12.5}  {\sl Suppose that  $F$ is a constant coefficient subequation on $X=\rn$,
or more generally suppose that $X= K/G$ is a riemannian homogeneous space and
that $F$ is a subequation  which is invariant under the natural action of the Lie group
$K$ on $J^2(X)$.

If $\bo$ is both $F$ and $\ft$ strictly convex, then existence holds 
for the  Dirichlet problem on $\O$ as above.
}

\Remark{12.6} There is a hidden hypothesis in these existence results: the sets
 $\overrightarrow{F_\l}$ and $\overrightarrow{\ft}_\l$ for $\l\in\bbr$ 
 (see \S 11) must be non-empty in order for there to exist
any  domain $\O$ whose boundary is both $F$ and $\ft$ strictly convex.
For example, for the Eikonal subequation $|p|\leq 1$  this fails  and,  as is well known, existence
fails for general $\vf\in C(\partial B)$ where $B$ is a ball in $\rn$.
However, Theorem 12.5 does apply to the infinite Laplacian, defined by taking
$F$ to be the closure of the set $\{(r,p,A)\in J^2(\rn): \langle Ap,p\rangle >0\}$.
This example is self-dual, i.e., $\ft=F$, and it is also a cone with $\Fa=\Int F$.
It is easy to check that $(r,p,A)\in \Int F$ if and only if either $\langle Ap,p\rangle >0$   or  $A>0$ and $p=0$. 
Hence, strictly  convex functions are strictly $F$-subharmonic and define
domains for which existence holds.

\medskip

Now for the proofs.

\Lemma {F}
$$
U^*\bigr|_\O \ \in\ F(\O)
$$

 \Lemma{$\wt {\bf F}$}
 $$
- U_* \bigr|_\O \ \in\ \ft(\O)
$$

\Prop{ F}  {\sl Suppose $\bo$ is strictly $F$-convex at $x_0\in\bo$.  For each $\d>0$ small, there exists
$\uu\in \cf(\vf)$ with the additional properties: }
\smallskip 

$\bullet$ \ \  $\uu$ is continuous at $x_0$,

 \smallskip 

$\bullet$ \ \ $\uu(x_0) = \vf(x_0)-\d$,
\smallskip 

$\bullet$ \ \ $\uu \in F_{\rm strict}(\ob)$.

\Prop{$\wt {\bf F}$}  {\sl Suppose $\bo$ is strictly ${\ft}$-convex at $x_0\in\bo$.  For each $\d>0$ small, there exists
$\ou\in \wt \cf(-\vf)$ with the additional properties: }
\smallskip 

$\bullet$ \ \  $\ou$ is continuous at $x_0$,

 \smallskip 

$\bullet$ \ \ $\ou(x_0) = - \vf(x_0)-\d$,
\smallskip 

$\bullet$ \ \ $\ou \in \ft_{\rm strict}(\ob)$.
 
 \Cor{F} 
 $$
 \vf(x_0) \ \leq \ U_*(x_0)
 $$
 
 \Cor{$\wt {\bf F}$} 
 $$
  \ U^*(x_0)\ \leq \  \vf(x_0)
 $$
 
 \bigskip
 \noindent
 {\bf Conclusion 1. (Boundary Continuity)} {\sl We have 
 $
 U_*= U = U^* = \vf  \ \  {\sl on\ } \bo
 $.
 In particular, $U$ is continuous at each point of $\bo$.}
 \pf  By Corollaries F and $\wt {\rm F}$ we have
 $\vf (x_0) \leq U_*(x_0) \leq U(x_0) \leq U^*(x_0) \leq \vf(x_0) \ \forall\, x_0$.\qed

 \medskip
 \noindent
 {\bf Conclusion 2.} \ \ $U^*\ \in\ \cf(\vf)$.
 \pf
Corollary $\wt {\rm F}$ together with Lemma F.\qed

 \medskip
 \noindent
 {\bf Conclusion 3.} \ \ $U=U^*$ on $\ob$.
 \pf
 We have $U^*\leq U$ on $\ob$ since $U^*\in \cf(\vf)$.\qed
 \vskip.3in
 
 \noindent
{ \bf Proof of Lemma F.}
Because of the families locally bounded above property, it suffices to show that 
$\cf(\vf)$ is uniformly bounded above.  Pick $\psi\in \ft_{\rm strict}(\ob)$ bounded below.
Pick $c>>0$ so that $\psi-c \leq -\vf$ on $\bo$.
By Lemma 7.7(ii) we have $\psi-c \in \ft_{\rm strict}(\ob)$.  Given $u\in \cf(\vf)$ we have 
$u+(\psi-c)\leq0$ on $\bo$.
Weak comparison for $\ft$ implies that 
$u+\psi-c\leq0$ on $\ob$.  Hence, $u\leq -\inf_{\ob}\psi +c$ for all $u\in \cf(\vf)$.\qed

\medskip
 \noindent
{ \bf Proof of Proposition F.}
Assume $\psi\in F_{\rm strict}(\ob)$ is bounded below on $\ob$.
Since $\bo$ is strictly $F$-convex at $x_0$, Theorem 11.12 states that there exist a local 
defining function $\rho$ for $\bo$ near $x_0$, and $r>0, \e_0>0$ and $C_0>0$ such that
in some local coordinates
$$\eqalign
{
\b(x) \ \equiv\ \vf(x_0) -\d +&C(\rho(x) - \e|x-x_0|^2) \ \in\ F_{\rm strict}(\overline{B(x_0, r)}\cap \ob)\cr
 &\forall\, C\geq C_0\ \ {\rm and\ \ } \forall \e\leq \e_0.
 }
\eqno{(12.1)}
$$
Shrink $r>0$ so that
$$
\vf(x_0)-\d \ <\ \vf(x)\quad {\rm on}\ \bo\cap B(x_0,r).
\eqno{(12.2)}
$$
 Pick  $N > \sup_{\bo} |\vf| + \sup_{\ob}\psi$ so that
 $$
 \psi - N\ <\ \vf - \d \quad{\rm on\ } \bo.
\eqno{(12.3)}
$$
Choose $C$ so large that on $A\equiv (B(x_0,r) \sim B(x_0,r/2))\cap\ob$ we have
$$
 \b \ <\ \psi - N  
\eqno{(12.4)}
$$
This is possible since $\rho(x) -\e|x-x_0|^2$ is strictly negative on $A$ and $\psi-N$ is bounded below on $A$.

By (12.4) we have that
$$
\uu(x) \ =\ \max\{ \b, \psi-N\}
$$
is a well defined   function on $\ob$ which is equal to $\psi-N$ outside $B(x_0,r/2)$.

Now we have  $\psi -N\in \ F_{\rm strict}(\ob)$, because Condition (N) is satisfied. 
Thus Lemma 7.7(i) and condition (12.1) imply that $\uu\in\ F_{\rm strict}(\ob)$.

Outside the set  $B(x_0,r/2)\cap \ob$ we have $\uu = \psi -N$ 
which is $\leq \vf$ on $\bo$ by (12.3), while
on  $B(x_0,r/2)\cap \bo$ we still have $\psi-N \leq \vf$, but also
$$
\b\ =\ \vf(x_0)-\d+C(\rho(x)-\e|x-x_0|^2) \ =\ 
         \vf(x_0)-\d-C \e|x-x_0|^2 \ \leq \vf(x_0)-\d \ \leq \ \vf(x)
$$
by (12.1).  Thus $\uu\bigr|_{\bo}\leq \vf$.  This proves that $\uu\in\cf(\vf)$. Finally note that
$\b(x_0) =  \vf(x_0)-\d$ which is $>\psi(x_0)-N$ by (12.3).
 Continuity of $\b$ and the upper semicontinuity 
of  $\psi$  implies that $\b> \psi-N$ in a neighborhood of $x_0$.  That is $\uu =\b$ in a neighborhood of $x_0$.  Thus $\uu$ is continuous at $x_0$, and $\uu(x_0)=\vf(x_0)-\d$.\qed

\medskip
 \noindent
{ \bf Proof of Proposition $\wt {\bf F}$.} This is merely Proposition F with an exchange of roles.\qed

\medskip
 \noindent
{ \bf Proof of Corollary F.}
 Since  $\uu\in\cf(\vf)$, we have $\uu\leq U$.  Hence, $\uu_*\leq  U_*$.  
 By the continuity of $\uu$ at $x_0$ and the fact that 
 $\uu(x_0) = \vf(x_0)-\d$, we have $\vf(x_0) - \d\leq U_*(x_0)$ for all $\d>0$ small.\qed

\medskip
 \noindent
{ \bf Proof of Corollary $\wt {\bf F}$.}
 Here we use weak comparison again.  Choose $u\in \cf(\vf)$.  Since $\ou\leq -\vf$
 on $\bo$, we have $u+\ou\leq 0$ on $\bo$.  Since $\ou\in \ft_{\rm strict}(\ob)$, 
 weak comparison implies that  $u+\ou\leq 0$ on $\ob$.  Therefore, 
 $U+\ou\leq 0$ on $\ob$, i.e. $U\leq -\ou$ on $\ob$.  Since $\ou$ is continuous at $x_0$
 and $\ou(x_0) = -\vf(x_0) -\d$, this implies that $U^*(x_0) \leq \vf(x_0)+\d$ for all $\d>0$ small.\qed

\medskip
 \noindent
{ \bf Proof of Lemma $\wt {\bf F}$.}
 Note that the Conclusions 1, 2 and 3 are now established.
 Suppose $-U_*\bigr|_\O \notin \ft(\O)$. Then by Lemma 2.4 there exist $x_0\in \O$, $\e>0$,
 and $\psi \in C^2$ near $x_0$ so that in local coordinates
 \smallskip
$$
\eqalign
{
&(1)\ \ -U_* -\psi \ \leq \ -\e|x-x_0|^2 \quad {\rm near\ } x_0  \cr
&(2)\hskip .75in = \ \ 0 \quad\quad\quad\quad\quad {\rm at\ } x_0  \cr
}
$$
but
$$
J^2_{x_0}\psi \ \notin \ \wt F_{x_0},
\qquad {\rm i.e., }\qquad
-J^2_{x_0}\psi \ \in \ \Int F_{x_0}.
$$
Then there exist $r>0, \d>0$ so small that
$$
u \ \equiv \ -\psi +\d \ \ {\rm is \ } F\!-\!{\rm subharmonic\ on\ } B(x_0,r)
$$
Moreover, (1) implies that for $\d>0$ sufficiently small, 
\medskip
\centerline{    (1)$'$ \ \ $u\ <\ U_*$ on a neighborhood of $\partial B(x_0,r)$}
\medskip
\noindent
Since $U_* \leq U$, statement (1)$'$ implies that the function
$$
u'\ \equiv\  \cases{
\quad U \qquad\qquad \  {\rm on \ } \ \ \ob-B(x_0,r)  \cr
\max\{U, u\}\ \ \  \ \  {\rm on} \ \ \overline {B(x_0,r)}  \cr
}
$$
is $F$-subharmonic on $\O$.
Conclusion 1 says that   $U=\vf$ on $\bo$, which implies $u'\in \cf(\vf)$.
Therefore, $u'\leq U$ on $\ob$, which implies in turn that $u\leq U$ on 
$B(x_0,r)$.

Now statement (2) above says $U_*(x_0) = u(x_0)-\d$.
Pick a sequence $x_k\to x_0$ with $\lim_{k\to \infty} U(x_k) = U_*(x_0)$.
Then
$$
\eqalign
{
&(i) \ \ \ \lim_{k\to\infty} U(x_k)\ =\ u(x_0)-\d    \cr
&(ii) \ \ \ \lim_{k\to\infty} u(x_k)\ =\ u(x_0).    \cr
}
$$
This implies that $u(x_k)>U(x_k)$ for all $k$ large,   contradicting  the
fact that $u\leq U$ on 
$B(x_0,r)$.
\qed

 \bigskip
  This completes the proof of Theorem 12.2  and its Corollary 12.3.
 Interior continuity is all that remains in showing that $U\bigr|_\O$ is $F$-harmonic.
 \medskip
 
 \noindent
 {\bf Proof of Theorem 12.4. (Assuming Comparison):}  By Corollary F and Lemma $\wt{\rm F}$
 $$
 -U_*\ \in\ \wt{\cf}(-\vf).
 $$
 In particular
 $$
 U-U_* \ \leq\ 0\qquad{\rm on\ }\bo.
 $$
Since $U\bigr|_\O \in F(\O)$ and   $-U_*\bigr|_\O \in \ft(\O)$, comparison implies
$$
U-U_* \ \leq\ 0\qquad{\rm on\ }\ob,
$$
 that is, $U\leq U_*\leq U$.  Since $U=U^*$ we are done. \qed

  \bigskip

 \noindent
 {\bf  Proof of Theorem 12.5.  (Assuming Constant Coefficients):} 
 We suppose that $X=\rn$ and that $F$ has constant coefficients.
 Let $\O_\d \equiv \{x\in \O : \dist(x, \bo) >\d\}$ and 
 $C_\d  \equiv \{x\in \ob : \dist(x, \bo <\d\}$.
 Suppose $\e>0$ is given.  By the continuity of $U$ at points of $\bo$ and the compactness of 
 $\bo$ it follows easily that there exists a $\d>0$ such that 
 $$
 {\rm if }\ |y|\leq \d, \ \ {\rm then\ } U_y\leq U+\e\ \ {\rm on\ }   C_{2\d}
 \eqno{(12.5)}
 $$
 where $U_y(x) \equiv U(x+y)$ is the $y$-translate of $U$ and where we define $U$ to 
 be $-\infty$ on $\rn-\ob$.  We claim that 
 $$
 {\rm if }\ |y|\leq \d, \ \ {\rm then\ } U_y\leq U+\e\ \ {\rm on\ } \ob
 \eqno{(12.6)}
 $$ 
 Setting $z=x+y$, this implies that 
 $$
 {\rm if }\ z\in \ob, \ x\in \ob, \  {\rm and}\ |z-x|\leq \d,   \ \ {\rm then\ } U(z)\leq U(x) +\e
 $$
 and therefore by symmetry that $|U(z)-U(x)|\leq \e$. Thus the proof is complete once
 (12.6) is established.
 
 To prove (12.6), note first that $U_y-\e  \in  F(\O_\d)$ for each $|y|<\d$ by the translation
invariance of $F$ and property (N).  Since $U_y \leq U+\e$ on the collar $C_{2\d}$, one has
 $$
 g_y \ \equiv\ \max\{ U_y -\e, U\} \ \in\ F(\O).
 $$
 Now (12.5) implies that $g_y = U$ on $C_{2\d}$. Therefore, 
 $$
 g_y\ \in\ \cf(\vf),
 $$
 and hence $g_y\leq U$ on $\ob$.  This proves that
 $$
 U_y-\e\ \leq g_y\ \leq\ U\qquad{\rm on\ \ } \O_\d.
 $$
 Combined with (12.5) this proves (12.6). \qed

\medskip\noindent
This Walsh argument  evidently applies to domains
$\O$ in any riemannian homogeneous space $X=G/H$ provided that the 
equation $F$ is invariant under the action of $G$ on $J^2(X)$. \qed

\medskip

Our existence result for general constant coefficient subequations can be restated as follows.

\Theorem{12.7} {\sl  Suppose that $F\ss J^2(\rn)$ is a constant coefficient subequation,
and that $\O$ is a domain with smooth boundary $\bo$ which is both $F$ and $\ft$
strictly convex.  Given $\vf\in C(\bo)$, let $U$ denote the $F$-Perron function on $\ob$ with boundary
values $\vf$, and let $V$ denote the $\ft$-Perron function on $\ob$ with boundary
values $-\vf$.  Then both $-V$ and $U$ are solutions to the Dirichlet problem for $F$ with boundary values
$\vf$.  Furthermore, any other such solution $u$  satisfies
$$
-V\ \leq\ u\ \leq\ U\qquad {\sl on}\ \ \ob.
$$
 }

 \pf  This follows easily from Theorem 12.5.
 \qed
 
 \Ex{12.8. (Existence without Uniqueness)}  Consider the dual subequations
 $$
 \eqalign
 {
 F \ &\equiv\  \left \{(r,p,A) : A-\half|p|^{1\over 2}\left ( I+P_{[p]}\right) \geq 0\right\}
 \and  \cr
 \ft \ &\equiv\  \left \{(r,p,A) : A+\half|p|^{1\over 2}\left ( I+P_{[p]}\right) \geq 0 \right\}
 \cr
 }
\eqno{(12.7)}
$$
 (where $P_{[p]}$ us orthogonal projection onto the $p$-line). First note that the 
 asymptotic interior of both $F$, and of $\ft$, is $\Int \cp$.  Consequently, the boundaries
 $\bo$ which are strictly $F$- and ${\ft}$-convex are precisely the boundaries which
 are classically strictly convex.  Consider now the Dirichlet problem for the $R$-ball
 $\O = \{x:|x|<R\}$ with boundary-value  function $\vf = 0$ (hence $\wt \vf = -\vf = 0$ also).
 Then Theorem 12.7 applies.  Moreover,
 \smallskip
 
 $\bullet$ \ \ The $F$-Perron function is $U\equiv 0$.

 \smallskip
  
 $\bullet$ \ \ The $\ft$-Perron function is $V(x) \equiv {1\over 12}(R^3 -|x|^3)$.

\smallskip
 \noindent
 Thus, in particular, while existence and weak comparison hold (cf. Theorem 10.1 or 
 Corollary C.3), local comparison fails here.

\pf 
Obviously $U\equiv 0$ is $F$-harmonic.  Moreover, 
$$
-U(x) + \smfrac \e 2 |x|^2 \ =\ \smfrac \e 2 |x|^2 \ \ {\rm is \ strictly\ }  \ft \ {\rm subharmonic}.
\eqno{(12.8)}
$$
Assuming (12.8), weak comparison  implies that on $\ob$:
$$
u(x) + ( - U(x) + \smfrac \e 2 |x|^2) \leq \smfrac \e 2 R^2
\fa u\in\cf(\vf).
$$
 Thus, $u\leq U$ for all  $u\in\cf(\vf)$, i.e., $U$ is the $F$-Perron function for $\vf=0$.

Similarly, calculation shows that $V(x) = {1\over 12}(R^3 -|x|^3)$ is $\ft$-harmonic.
Moreover, 
 $$
-V(x) + \smfrac \e 2 |x|^2  \ \ {\rm is \ strictly\ }  F \ {\rm subharmonic}.
\eqno{(12.9)}
$$
Assuming (12.9), weak comparison  implies that on $\ob$:
$$
(-V(x)  + \smfrac \e 2 |x|^2) + v(x) \ \leq \ \smfrac \e 2 R^2
\fa v\in  \wt\cf(-\vf).
$$
 Thus, $v\leq V$ for all  $v\in \wt\cf(\vf)$, i.e., $V$ is the $\ft$-Perron function for $-\vf=0$.

To prove (12.8) simply note that the Hessian of $-U(x)+ \smfrac \e 2 |x|^2$ is $A=\e \cdot I$
and apply the definition. To prove (12.9) note that for the function $-V(x) + \smfrac \e 2 |x|^2$
$$
p\ =\ \left(\smfrac {|x|} 4 +\e \right )x\and A\ =\ \left(\smfrac {|x|} 4 +\e \right )\cdot I + \smfrac {|x|} 4 P_{[x]}.
$$
The eigenspaces for $A-\half |p|^{\half} (I+P_{[p]})$ are:

\centerline{
$[x]$ with eigenvalue 
${|x|\over 2} +\e- \left ({|x|^2\over 4} +\e|x| \right )^{1\over2}$
  and\   $[x]^\perp$ with eigenvalues 
${|x|\over 4} +\e-  \half\left ({|x|^2\over 4} +\e|x| \right )^{1\over2}$.
}

\noindent
It is easy to see that these eigenvalues are $>0$ on $\ob$. \qed
\medskip

Despite the above examples of strict approximation,  it must fail in general
for this equation since otherwise uniqueness would hold.  More precisely, 
$U$ has no $F$-strict approximation and $V$ has no $\ft$-strict approximation.

 \medskip\noindent
 {\bf Remark 12.9. (Refinements).}
 
 1) \   Note that in the proof of Proposition F we only needed  
 the barrier $\b$ defined in (12.1)  to be strict on the $\O$-portion of the ball,
  that is we only 
really  needed $$\beta \in F_{\rm strict}(\O\cap B(x_0, r)).$$
This can sometimes be established even though $\bo$ is not strictly $F$-convex.
 For example, suppose $F$ is defined by
 $$
 A\ \geq\  0  \and \det A\  \geq\  e^r
\eqno{(12.10)}
 $$
 and note that $\Fa$ is defined by
 $$
A\ >\ 0\and r\ <\ 0.
 $$
Consider $\O = \{x : |x| < R\}$ and set $\rho(x) = \half(|x|^2-R^2)$. Note that with 
$\l = \vf(x_0)-\d$    the barrier satisfies $\b(x)\leq \l$ on the $\O$-portion of $B(x_0,r)$.
Then it is easy to see that for $C >0$ sufficiently large,
we have   $\beta \in F_{\rm strict}(\O\cap B(x_0, r))$.  This is because $J^2_x \beta = (\b(x), C(x-2\e(x-x_0)), 
 C(1-2\e)I)$, so that $\det A - e^r \geq C^n(1-2\e)^n - e^\l > 0$ for all 
 $x\in \O\cap B(x_0, r)$.  Thus Proposition F and Corollary F are valid for this $F$ in spite of the 
 fact that $\bo$ is not strictly $\Fa$-convex.  In fact, Proposition F and Corollary F are false for the subequation
 $H\equiv {\rm closure}\{\Fa\}$ unless $\vf(x_0)\leq 0$.  It is now easy to show that existence and uniqueness 
 hold  for the subequation $F$ defined by (12.10), since $\psi(x) = \half(|x|^2-R^2)$ is a good approximator for $F$.
 
 \medskip

 2)\ Notice that the weakened form of Proposition F with $\underline u \in F(\ob)$ but not necessarily strict,
 is all that is needed to prove Corollary F.   However,  to prove Corollary $\wt{\rm F}$ from Proposition $\wt{\rm F}$ 
 weak comparison was used so that the strictness of $\overline u\in {\wt F}_{\rm strict}(\ob)$ cannot be dropped.

 \vskip .3in


\centerline{\headfont \ 13.\   The Dirichlet Problem -- Summary Results.}
\medskip

We present here several summary results which follow from the work above.
We make the following standing hypotheses in this section.

\medskip

(i)  \ \ \ $F$ is a subequation on a riemannian manifold $X$.

\medskip

(ii)  \ \ $\O \ss\ss X$ is a domain with smooth boundary $\bo$.

\medskip

(iii)  \ Both $F_{\rm strict}(\ob)$ and  ${\ft}_{\rm strict}(\ob)$ have at least one function
         bounded below.

 \medskip
 Consider the following\medskip
 \medskip
 \noindent
  {\bf DIRICHLET PROBLEM for $F$:}    
 \medskip
 \noindent
 Given $\vf\in C(\bo)$, consider the  Perron function 
$$
U\ \equiv \sup_{u\in \cf(\vf)} u \qquad {\rm where\ \ }  
\cf(\vf) = \left \{ u \in F(\ob) :  u\bigr|_{\bo}\leq \vf\right\}.
$$

\noindent
{\bf EXISTENCE.} For each $\vf \in C(\bo)$ the Perron function satisfies  
\medskip
 
$\bullet$\ \  $U$ is $F$-harmonic 
 
 \medskip
 
$\bullet$\ \   $U = \vf$\ \  on $\bo$
 
 \medskip
 
$\bullet$\ \   $U\in C(\ob)$.
  \bigskip
 
\noindent
{\bf UNIQUENESS.} $U$ is the only function with these three properties.

\bigskip

In what follows, $M_F$ denotes a monotonicity cone for $F$.

\Theorem{13.1} {\sl Suppose $F$ is a riemannian $G$-subequation 
where  $X$  is provided with a topological $G$-structure.
Suppose there exists a 
$C^2$ strictly $M_F$-subharmonic function on $X$.

Then for every domain $\O\ss\ss X$ whose boundary is strictly
$F$- and ${\ft}$-convex, both existence and uniqueness
hold for the Dirichlet problem.}

\Theorem{13.1$'$} {\sl Theorem 13.1 also holds for any subequation $F$ which is 
locally affinely jet-equivalent to a riemannian $G$-subequation on $X$.}

\Theorem{13.2} {\sl Suppose $F$ is a   subequation for which 
weak comparison holds.
Suppose there exists a 
$C^2$ strictly $M_F$-subharmonic function on $X$.

Then for every domain $\O\ss\ss X$ whose boundary is strictly
$F$- and ${\ft}$-convex, both existence and uniqueness
hold for the Dirichlet problem.}

\Theorem{13.3} {\sl Suppose comparison holds for a   subequation   $F$ on $X$.

Then for every domain $\O\ss\ss X$ whose boundary is strictly
$F$- and ${\ft}$-convex, both existence and uniqueness
hold for the Dirichlet problem.}

\Theorem{13.4} {\sl Let $\bbf$  be a constant coefficient subequation on $\rn$.

{\bf (Existence.)}  Then for every domain $\O\ss\ss X$ whose boundary is strictly
$ \bbf$- and ${\wt \bbf}$-convex,   existence  
holds for the Dirichlet problem.
 \smallskip
 
{\bf (Comparison.)}  If $\bbf$ is pure second-order, or more generally, independent of the gradient,
then comparison holds for $\bbf$ on $\rn$.}
 
 \medskip
 
 This provides a new proof of existence and uniqueness for pure second-order subequations, established
 in [HL$_4$] using a result of Slodkowski [S].
 
 \pf The existence is just a restatment of Theorem 12.7.
 Since weak comparison holds for any constant coefficient subequation,
 comparison will follow from strict approximation.  A closed subset $\bbf\ss\bbj^2$ is a pure 
 second-order subequation if and only if $M=\bbr\times\rn\times\cp$ is a monotonicity set for $\bbf$.
 Since $|x|^2$ is strictly $M$-subharmonic, strict approximation holds for $\bbf$.  More generally, recall
 that a closed subset $\bbf\ss\bbj^2$ is independent of the gradient if and only if 
 $M\equiv \bbr_-\times\rn\times\cp$ is a monotonicity set for $\bbf$. In this case the function
 $|x|^2- R^2$ is strictly $M$-subharmonic on the ball of radius $R$ about the origin.  Since each
 compact subset $K\ss \rn$ is contained in such a ball, strict approximation holds for $\bbf$.\qed

 \medskip
More generally we have the following (see Theorem 9.13).

\Theorem{13.5} {\sl Let $X = K/G$  be a riemannian homogeneous space and
suppose $F$ is a subequation which is invariant under the natural action of 
the Lie group $K$ on $J^2(X)$.

{\bf (Existence.)}  Then for every domain $\O\ss\ss X$ whose boundary is strictly
$F$- and ${\ft}$-convex,   existence  
holds for the Dirichlet problem.
 \medskip
 {\bf (Comparison.)}  If $X$ supports a strictly convex $C^2$-function,
then comparison holds for  $F$ on $X$.}

\vfill\eject


\centerline{\headfont \   14.\   Universal Riemannian Subequations.}
\medskip

In this section we consider  subequations defined on any riemannian manifold  by the requirement that
 $\Hess_x u \in  \bbf$   (for   $u\in C^2$), where $\bbf$ is a closed subset of $\Symn$ which is 
 O$_n$-invariant.  Recall from Section 4.4, that for such  purely second-order closed subsets of $J^2(X)$, 
 condition (N) is automatic and  
condition (P) implies condition  (T).
 
 Each O$_n$-invariant closed subset $\bbf \ss \Symn$  determines a closed subset $\L$,
invariant under the permutation group $\pi_n$ and consists of all 
$n$-tuples of eigenvalues of $A$ where $A\in \bbf$.  Conversely,  each closed subset $\L\ss \rn$ invariant
under $\pi_n$ determines a closed O$_n$-invariant subset $\bbf \ss \Symn$ namely
$\bbf \equiv \{A : \l(A) \in\L\}$.   Moreover,
 $$
 \bbf \ \ {\rm satisfies\ (P)} \quad\iff\quad \L+\rn_+ \ \ss\ \L
 \eqno{(14.1)}
 $$
 where $\bbr_+ = [0, \infty)$.
 The implication from left to right in (14.1) is obvious since $\L$ can be taken to be the  subset 
 of diagonal elements in $\bbf$ and $\rn_+$ the set of diagonal elements in $\cp$.
 To prove the reverse implication consider the {\sl ordered eigenvalues}
 $$
\l_1(A) \ \leq \ \cdots \ \leq \ \l_n(A)
 \eqno{(14.2)}
 $$
 of $A\in \Symn$.  Here $\l_k(A)$, the $k$th smallest eigenvalue, is a well defined continuous function 
 on $\Symn$.  The standard fact needed is the {\sl monotonicity of the ordered eigenvalues}:
  $$
A\ \leq \ B \quad \Rightarrow\  \l_k(A) \ \leq \   \l_k(B) \fa k
 \eqno{(14.3)}
 $$
  which follows from the minimax definition of $\l_k(A)$.
 
  Set 
 $$
 \l(A)\ =\ (\l_1(A),...,\l_n(A)).
 $$
Now assume $A\in \bbf$ and $P\geq 0$ so that $B=A+P\geq A$.  Monotonicity (14.3) says that $\l(B)$ equals
$\l(A)$ plus a vector in $\rn_+$.  Since $\l(A)\in \L$, the assumption $\L+\rn_+\ss \L$ implies 
that $\l(B)\in \L$.  Hence $B\in \bbf$, that is, $\bbf$ satisfies (P), and  (14.1) is proved.

 \Def{14.1}  Suppose $\L$ is a closed subset of $\rn$ invariant under the permutation group $\pi_n$.
 If $\L$ satisfies 
 $$
 \L+\rn_+\ \ss\ \L,
  \eqno{(14.4)}
 $$
then $\L$ will be referred to as a {\bf positive} (or {\bf  $\rn_+$-monotone}) set in $\rn$.
 \medskip
 
In the final Remark 14.11 we give a canonical description of all possible 
$\rn_+$-monotone sets.  Otherwise,  the remainder of this section is devoted to examples. 
Our discussion is confined to
 subsets $\bbf\ss \Symn$, but each $\bbf$ corresponds to a universal subequation
 in riemannian geometry, and our language   reflects that fact.
 
 \Ex{14.2. (Monge-Amp\`ere -- The Principal Branch)}   The Monge-Amp\`ere equation det$(A)=0$
 gives rise to several subequations or branches.  The principal branch is just the subequation $\cp$.
 In terms of the ordered eigenvalues it is given by:
 $$
 \l_1(A)\ \geq \ 0.
 $$
 We call $\cp$ a {\sl branch} because $\partial \cp\ss\{\det = 0\}$. However, there are further natural
 subequations with this property as we shall see below.
 
 One can show that  on a connected open set $X\ss \rn$ a function $u$ is $\cp$-suharmonic if and only if 
 $u$ is convex (or $\equiv -\infty$).  
 
 The dual subequation $\cpt$ is defined by 
 $$
 \l_n(A)\ \geq \ 0.
 $$
 In [HL$_4$] we proved that on open sets $X\ss\rn$, a function $u$ is $\cpt$-subharmonic
  if and only if $u$ is {\sl subaffine}, that is,
 $$
 u\leq a \ \ {\rm on\ } \partial K \ \ \ \Rightarrow\ \ \  u\leq a \ \ {\rm on\ } K
 $$
 for each affine function a and each compact subset $K\ss X$.
 Thus $u$ is $\cp$-harmonic if and only if $u$ is convex and $-u$ is subaffine.
 Since $\partial \cpt \ss\{\det=0\}$, the subequation $\cpt$ is another branch of the Monge-Amp\`ere equation.

 \Ex{14.3. (The Other Branches of the Monge-Amp\`ere Equation)}   
 Define the $q$th branch $\cp_q$ by the condition that
 $$
 \l_{q}(A)\ \geq\ 0 \qquad  ({\rm at\ least\  } n-q+1\ {\rm eigenvalues\ } \geq 0).
 $$
 This gives  $n$ branches or subequations
 $$
 \cp \ =\ \cp_{1} \ \ss\ \cdots\ \ss\ \cp_q\ \ss\ \cdots\ \ss\ \cp_n\ =\ \cpt
 $$
 for the Monge-Amp\`ere equation since, for each $q$,
 $$
 \partial \cp_q \ \ss\ \{\det = 0\}
 $$
 and the positivity condition
 $$
 \cp_q + \cp \ \ss\ \cp_q
 $$
  is satisfied by (14.3). Duality becomes
  $$
  \cpt_q\ =\ \cp_{n-q+1}
  $$
 since $\l_q(-A) \ =\ -\l_{n-q+1}(A)$.  The principle branch $\cp$ is a  monotonicity subequation
for each branch $\cp_q$.  Since each $\cp_q$ is a cone, we have 
$$
\oa{\cp_q} \ =\ \cp_q.
$$
 Thus the boundary of a domain $\O\ss X$ is strictly $\cp_q$-convex at a point
 $x$  iff its second fundamental form $II_{\bo}$  (with respect to the interior-pointing normal)
 has at least $n-q$ principle curvatures $> 0$. Thus, $\bo$ is strictly
 $\cpt_q=\cp_{n-q+1}$-convex if $II_{\bo}$ has at least $q-1$ principle curvatures $> 0$.
 
 Theorem 13.1 gives us the following result.
 
 \Theorem {14.4} {\sl Let $\O\ss\ss X$ be a domain with smooth boundary in a riemannian manifold $X$.
 Suppose $\O$ admits a smooth strictly convex global defining function.  Then the Dirichlet problem
 for   every branch of the real Monge-Amp\`ere equation is uniquely solvable for all continuous
 boundary functions.}\medskip

 \Ex{14.5. (Geometric $p$-Plurisubharmonicity)}   
 Let $G(p,\rn)$ denote the grassmannian of all $p$-dimensional subspaces of $\rn$.
 Define $\bbf(G(p,\rn))$ by requiring that 
 $$
 \tr_W  A   \ \geq \ 0  \fa W\in G(p,\rn)
 \eqno{(14.5)}
 $$
  As with all geometrically defined subequations, this is a convex cone
 subequation.  Functions which are $F(G(p,\rn))$-subharmonic
 are more appropriately called {\sl geometrically $p$-plurisubharmonic} since
 one can prove that a function $u$ is $\bbf(G(p,\rn))$-subharmonic if and only if 
 its restriction to every minimal $p$-dimensional submanifold of $X$ is subharmonic
 with respect to the induced riemannian metric (or $\equiv -\infty$).

 The subequation $\bbf(G(p,\rn))$ can be written in terms of the ordered eigenvalues as
 $$
 \l_{1}(A)+\cdots+\l_p(A)\  \geq\  0,
 $$
 i.e.,  all $p$-fold sums of the eigenvalues are $\geq0$.
 
 The dual subequation  $\wt \bbf(G(p,\rn))$ can be described by either of the two equivalent conditions:
 $$
\tr_W A\ \geq \ 0 \quad {\rm for\ some\ \ } W\in G(p,\rn)
$$
 $$
 \l_{n-p+1}(A) + \cdots + \l_n(A) \ \geq\ 0.
 $$

 Both the subequation $\bbf(G(p,\rn))$  and its dual  $\wt \bbf(G(p,\rn))$ are branches of a 
 polynomial equation.  Define a polynomial $M_p$ on $\Symn$ by
 $$
 M_p(A)\ =\ {\prod_{|I|=p}}' \left( \l_{i_1}(A) + \cdots + \l_{i_p}(A)  \right)
  \eqno{(14.6)}
 $$
 where the prime indicates that the product is over all multi-indices $I=(i_1,...,i_p)$ with
 $i_1< i_2 <  \cdots < i_p$.  Then 
 $$
 A\in \partial \bbf(G(p,\rn)) \ \ {\rm or}\ \  A\in \partial \wt \bbf(G(p,\rn)) \qquad\Rightarrow\qquad 
 M_p(A) \ =\ 0.
 $$

 The equation $ M_p(A) \ =\ 0$ has many branches, or subequations.
 Let 
 $$
 \l_I(A) \ \equiv\ \l_{i_1} + \cdots +\l_{i_p}
 $$
 denote the $I$th $p$-fold sum of ordered eigenvalues.  These $N= {n\choose p}$
 real numbers can be ordered as
 $$
 \l_1(p,A) \ \leq \ \cdots \ \leq \ \l_N(p,A).
 $$
 Define $\bbf_k, k=1,..., N$ by the condition that
 $$
 \l_k(p,A)\ \geq\ 0.
 $$
 That is, $A\in \bbf_k$ has at least $N-k+1$ $p$-fold sums of its eigenvalues $\geq 0$.
 Note that if $B\in \bbf_1 = \bbf(G(p,\rn))$ (the principal branch) and $A\in \bbf_k$, then
 $A+B\in \bbf_k$.  That is
 \bigskip
 \centerline{ $\bbf(G(p,\rn))$ is a monotonicity subequation for each branch $\bbf_k$.}
 \bigskip
 Since $\bbf_k$ is already a cone, we have that $\oa \bbf_k = \bbf_k$. The dual is given by
 $$
 \wt \bbf_k \ =\ \bbf_{N-k+1}.
 $$
 Note  that  by (11.6)$'$ strict boundary convexity of a domain $\O$ means that 
 $$
\left( \matrix {t & 0 \cr 
0 &  II_{\bo}
}\right) \  \in \ \bbf_k \fa t>0\ \ {\rm sufficiently\ large.}
 $$
  Note that the top ${n-1\choose p-1}$
 $p$-fold sums of eigenvalues of this matrix are automatically positive for large $t$.
 Hence, if ${n\choose p}-k+1 \leq {n-1\choose p-1}$ every boundary is automatically
 strictly $\bbf_k$-convex.  Otherwise, $\bbf_k$-boundary convexity means that 
 $II_{\bo}$ has at least ${n\choose p}-k+1 - {n-1\choose p-1}$ $p$-fold sums of
 its principal curvatures $ \geq 0$.

 We leave it to the reader to formulate corollaries of Theorem 13.1 for these equations
 (as we did for the Monge-Amp\`ere equation above).
 
  \Remark{14.6}  Each $A\in\Symn$ acts as a derivation $D_A$ on $\L^p\rn$
 and $D_A\in \Sym(\L^p\rn)$.  The polynomial $M_p$ is the restriction of the determinant 
  on  $\Sym(\L^p\rn)$ to Image$(D_A)$
 $$
 M_p(A)\ =\ \det\left(D_A\right).
 $$
 
%
 %

  \Ex{14.7. (Elementary Symmetric Functions)}
Recall that the cone $\cp$ can be defined by requiring that
$$
\s_1(A)\ \geq\ 0, \ \s_2(A)\ \geq\ 0,\    ... \ , \ \s_n(A)\ \geq\ 0 
$$
 where $\s_k(A)$ is the $k$th elementary symmetric function of the eigenvalues of $A$.
 For each $k$, $1\leq k\leq n$ the condition
$$
\s_1(A)\ \geq\ 0, \ \s_2(A)\ \geq\ 0,\    ... \ , \ \s_k(A)\ \geq\ 0 
$$
 defines a convex cone subequation $\bbf(\s_k)$. One can show that this set $\bbf(\s_k)$
 is exactly the closure of the connected component of the complement
 of $\{\s_k(A)=0\}$ which contains the identity $I$. 
 
 The equation $\{\s_k(A)=0\}$ has $k-1$ other branches 
 $$
 \bbf(\s_k)\ =\  \bbf_1(\s_k) \  \ss\  \bbf_2(\s_k) \  \ss\   \cdots \ss\  \bbf_k(\s_k) 
,
 \eqno{(14.7)}
 $$
 each of which is $\bbf(\s_k)$-monotone, and for which $\wt \bbf_j(\s_k) = \bbf_{k-j+1}(\s_k)$. 
 (See the general discussion in [HL$_7$].)  
 The last branch $\bbf_k(\s_k)$, which is the dual to $\bbf_1(\s_k)$,  is  given by the condition
 $$
 \s_1(A) \ \geq\ 0 \ \ \  {\rm or}\ \ \  - \s_2(A) \ \geq\ 0 \ \ \  {\rm or}\ \ \  \s_3(A) \ \geq\ 0 \ \ \  {\rm or}\ \ \ ... \ \ \ 
 {\rm or} \ \ \  (-1)^{k-1}\s_k(A)\ \geq\ 0.
 $$
 
 \Theorem{14.8}  {\sl  Let  $\O\ss\ss X$ be a  domain with smooth boundary in
 a riemannian manifold $X$.   Suppose $\O$ is globally strictly $\bbf(\s_k)$-convex,
 that is, suppose there is a strictly $\bbf(\s_k)$-subharmonic defining function for 
 $\O$.  Then for every branch of the equation $\s_k(\Hess \, u)=0$,
  the Dirichlet problem is uniquely solvable for all continuous boundary data.}
  
  \pf
  The existence of a strictly $\bbf(\s_k)$-subharmonic defining function $\rho$  implies that
  the boundary is strictly $\bbf(\s_k)$-convex,  and from the inclusions (14.7) it 
  is also $\bbf_j(\s_k)$-convex for each $j$.  Since $\wt \bbf_j(\s_k) = \bbf_{k-j+1}(\s_k)$,
  the boundary convexity hypothesis of Theorem 13.1 is satisfied. Since $\bbf(\s_k)$
  is a monotonicity subequation of $\bbf_j(\s_k)$ and $\rho$ is strictly  $\bbf(\s_k)$-subharmonic,
  the other hypothesis is satisfied and Theorem 13.1 applies.\qed
 
 \medskip

  \Ex{14.9. (The Special Lagrangian Potential Equation)}
This is the subequation $\bbf_c$ defined by the condition
$$
\tr \arctan A \ =\ \arctan(\l_1(A)) + \cdots + \arctan(\l_n(A))\ \geq\  {c \pi \over 2}.
$$
 for $c\in (-n, n)$. 
  Since $\arctan$ is an odd function, the dual equation is again of this form.  
 That is,
 $$
 \wt \bbf_c \ =\ \bbf_{-c}
 $$
 For values of $c$ where $\bbf_c$  is convex, the Dirichlet problem for this equation
 was studied in detail by Caffarelli, Nirenberg and Spruck [CNS] who established
 existence, uniqueness and regularity.  Existence, uniqueness and continuity
 for all other branches $\bbf_c$ were established  in [HL$_4$].
  
One computes  that the asymptotic interior of $\bbf_c$ is
 $$
\oa{{\bbf_c}} \ =\ \Int \cp_q \ = \ \{A : \l_q(A) >0\}
 \eqno{(14.8)}
$$
($\l_q(A)$ is the $q$th ordered eigenvalue) where $q$  is the unique integer such that
$$
{n-c\over 2} \ \leq \  q\ < \ {n-c\over 2}+ 1
\eqno{(14.9)}
$$

 \Theorem{14.10}  {\sl Let  $X$ be
 a riemannian $n$-manifold $X$ on which there exists some global strictly convex function.
Let  $\O\ss\ss X$ be a  domain with smooth boundary and suppose $\bo$ is strictly 
$\cp_q$-convex for an integer $q$ satisfying 
 $$
 1\ \leq \ q\ <\ {n\over 2}+1
 \eqno{(14.10)}
$$
Then the Dirichlet problem for $\bbf_c$-harmonic functions is uniquely solvable on
$\O$ for all continuous $\vf\in\bo$ and for all $c$ with}
$$
|c|  \   < \ n-2q+2
\eqno{(14.11)}
$$
 
 \pf
 Recall that strict $\cp_q$-convexity implies strict  $\cp_{q'}$-convexity for $q\leq q' \leq n$.
  Suppose $c\geq 0$ and (14.9) holds.  Let $\wt q$ be the integer satisfying (14.9) with $c$ replaced by $-c \leq 0$.  
 Then $\wt q \geq q$. Hence, by the first remark and (14.8) the hypothesis on the boundary
 is satisfied and Theorem 13.1 applies to $\bbf_c$. As $c$ descends to zero, the $q$ in (14.9) increases,
 so by the initial remark, the boundary hypothesis continues to be satisfied.\qed
 \medskip
 
 The subequations $\bbf_c$  are related to the equation
 $$
 {\rm Im}\left\{ e^{-i\theta} \det_\bbc (I+ i\Hess\, u)    \right\}\ =\ 0.
 \eqno{(14.12)}
$$
  (for fixed $\theta$) which arises in Special Lagrangian geometry [HL$_1$]. If $u$ satisfies
  (14.12), then the graph $\{y=\nabla u\}$ in $\rn\times \rn$ is absolutely volume-minimizing.
  Note that 
  $
   {\rm Im}\left\{ e^{-i\theta} \det_\bbc (I+ iA)    \right\}
   =
    {\rm Im}\left\{ e^{-i\theta} \prod_k (I+ i\l_k(A))    \right\} =0
  $ yields the congruence
  $$
  \sum_k \arctan(\l_k(A))\ \equiv \ \theta \ ({\rm mod} \ \pi)
  $$
  Thus the equation (14.12) has many disjoint connected sheets corresponding
  to the subequations $\bbf_{2({\theta\over \pi} +k)}$
  for either $n$ or $n-1$ integer values of $k$.  The maximal (and minimal) values were treated
  in [CNS] and the other values in the euclidean case in [HL$_4$].

  \medskip
  An interesting case where all of the above applies is the cotangent bundle $T^*K$ of a Lie group $K$. Fixing an orthonormal
  framing $e_1,...,e_n$ with respect to a left-invariant metric gives a splitting
  $$
  T^*K\ =\ K\times \rn
  $$
  which determines in an obvious way  a hermitian  almost complex structure and an
  invariant $(n,0)$-form, i.e., an almost Calabi-Yau structure. As above we can solve
  the Dirichlet problem for the special Lagrangian potential equation 
  on $\cp_q$-convex  domains  $\O\ss K$ for $q$ as above.
  For each solution $u$  the graph
  of $\nabla u$ in $T^*K$ will be a special Lagrangian submanifold.

  \Remark{14.11. (A canonical form for general   $\rn_+$-monotone sets $\L$)}  
  Consider the hyperplane $H$ in $\rn$ perpendicular to $e=(1,...,1)$. 
  The boundary of the positive orphant   $ \rn_+$ is a graph over $H$.  More precisely for $\mu = (\mu_1,..., \mu_n) \in H$, 
 set $\| \mu \|^+ \equiv - \mu_{\rm min}$  where 
 $\mu_{\rm min} = \min\{\mu_1,..., \mu_n\}$.
 Then
 $$
 \partial \rn_+\ =\ \{\mu+\|\mu\|^+ e : \mu \in H\}.
 $$
 Similarly
  $$
 \partial \rn_-\ =\ \{\mu-\|\mu\|^- e : \mu \in H\}.
 $$
 where $\|\mu\|^- = \mu_{\rm max} = \max\{\mu_1,...,\mu_n\}$.
 
 One can characterize the positive sets $\L$ as follows.  
 There exists a function $f:H\to \bbr$, invariant under the action of $\pi_n$ on $H$, which is 
 $\|\cdot\|^{\pm}$-Lipschitz, i.e., 
  $$
-\|\mu\|^-  \ \leq f(\l + \mu) - f(\l)\ \leq\ \|\mu\|^+\fa \l, \mu \in H
  \eqno{(14.13)}
 $$
such that
 $$
 \L\ =\ \{\mu + te : \mu\in H \ {\rm and\ } t\geq f(\mu)\}.
  \eqno{(14.14)}
 $$
Finally note that $\| \l+\mu\|^{\pm} \leq \| \l\|^{\pm} + \| \mu\|^{\pm}$ and $\| \mu \|^-= \| -\mu \|^+$.
 \medskip

 \vskip .3in
 

\centerline{\headfont \   15.\   The Complex and Quaternionic Hessians.}
\medskip

Virtually all  the  results of the previous section carry over directly to the complex 
and quaternionic hessians, that is, to the case of U$_n$ and ${\rm Sp}_1\cdot {\rm Sp}_n$ invariant 
equations.
\medskip
\noindent
{\bf The Complex Case.} Consider $\bbc^n= (\bbr^{2n}, J)$ where $J:\bbr^{2n}\to \bbr^{2n}$ denotes the
standard complex structure.  Then we have the set of  {\sl hermitian symmetric matrices}
$$
\Sym_\bbc(\bbc^n) \ =\ \{A\in \Sym(\bbr^{2n}) : AJ = JA\}
$$
 and the natural projection
 $$
 \pi:   \Sym(\bbr^{2n}) \arr \Sym_\bbc(\bbc^n) \qquad{\rm given \ by\ } \pi(A) = \half(A-JAJ).
 $$
 If  $A\in \Sym_\bbc(\bbc^n)$ and $A(e)= \l e$, then $A(Je) = \l Je$, and so 
 $\bbc^n$ decomposes into a direct sum of $n$ complex eigenlines  with
 eigenvalues $\l_1(A) \leq \cdots\leq  \l_n(A)$.  Let $\cp^\bbc = \pi(\cp) = \{\l_1(A) \geq 0\}$
 be the cone of positive hermitian symmetric matrices.
 
 Note that closed subsets $\bbf\ss  \Sym_\bbc(\bbc^n)$ invariant under U$_n$ are uniquely
  determined by their eigenvalue
 set $\L\ss\rn$ (invariant under $\pi_n$).  Moreover, one again has monotonicity
 of the ordered eigenvalues, so that positivity for invariant sets can be expressed as
 $$
\pi^{-1}(\bbf) +\cp \ \ss\  \pi^{-1}(\bbf) \quad \iff 
\quad \bbf +\cp^\bbc \ \ss\ \bbf \quad \iff \quad \L+\rn_+ \ \ss\ \L.
 $$ 
 Thus each positive $\pi_n$ invariant set $\L\ss \rn$ determines a  pure second order
 U$_n$-invariant subequation $\pi^{-1}(\bbf)$. All these equations carry over to any complex 
 (or almost complex) manifold $X$ with a hermitian metric.
 
 The complex Monge-Amp\`ere equation
 $$
 \det_\bbc A \ =\ \l_1(A) \cdots \l_n(A)
 $$
 behaves exactly as in the real case with principle branch $\cp^\bbc$ defined by
 $\{\l_1(A)\geq 0\}$ and the remaining branches defined by $\l_k(A)\geq 0\}$.
 The cone $\cp^\bbc$ is a monotonicity subequation for each of the branches.
  The   $\cp^\bbc$-subharmonic functions are exactly the classical plurisubharmonic functions
 on $X$.    The boundary
 of a domain is strictly $\cp^\bbc$-convex if and only if it is classically strictly pseudoconvex.
 Theorem 13.1 now gives the following.
 
 \Theorem {15.4} {\sl Let $\O\ss\ss X$ be a domain with smooth boundary
  in a almost complex hermitian manifold $X$.
 Suppose $\O$ admits a smooth strictly plurisubharmonic global defining function.  Then the Dirichlet problem
 for   every branch of the  complex  Monge-Amp\`ere equation is uniquely solvable for all continuous
 boundary functions.
 
 Furthermore, if $X$ carries some strictly plurisubharmonic function on a neighborhood of
 $\ob$, then one can uniquely solve the Dirichlet problem for the branch $\cp_q$
 if and only if the second fundamental form of $\bo$ at least $\max\{q, n-q+1\}$
 principal curvatures $>0$ at each point.}\medskip
 
 This generalizes a result of Hunt and Murray [HM] and Slodkowski [S]
 to the case of almost complex manifolds.

 The discussion of geometric $p$-plurisubharmonicity, elementary symmetric functions 
 and the special Lagrangian potential equation all carry over, and analogues of 
 Theorems 14.8 and 14.10 hold.  
 \bigskip
\noindent
{\bf The Quaternionic Case.} Consider $\bbh^n= (\bbr^{4n}, I,  J, K)$ where $I,J,K :\bbr^{4n}\to \bbr^{4n}$ denote 
right scalar multiplication by the unit quaternions $i,j,k$.
 Then we have the set of  {\sl  quaternionic hermitian symmetric matrices}
$$
\Sym_\bbh(\bbh^n) \ =\ \{A\in \Sym(\bbr^{4n}) : AI=IA, AJ = JA \ {\rm and\ } AK=KA\}
$$
 and the natural projection
 $$
 \pi:   \Sym(\bbr^{2n}) \arr \Sym_\bbc(\bbc^n) \qquad{\rm given \ by\ } \pi(A) = \smfrac14(A-IAI-JAJ-IKI).
 $$
 If  $A\in \Sym_\bbh(\bbh^n)$ and $A(e)= \l e$, then $A(Ie) = \l Ie$, $A(Je) = \l Je$ and $A(Ke) = \l Ke$, and so 
 $\bbh^n$ decomposes into a direct sum of $n$ quaternionic eigenlines  with
 eigenvalues $\l_1(A) \leq \cdots\leq  \l_n(A)$.  Let $\cp^\bbh = \pi(\cp) = \{\l_1(A) \geq 0\}$
 be the cone of positive quaternionic hermitian symmetric matrices.
The discussion now completely parallels the one given for the complex case above.
The discussion  carries  over with no change.

\vfill\eject


\centerline{\headfont \   16.\   Geometrically Defined Subequations -- Examples.}
\medskip

The subequations $\bbf(\GG)$ which are geometrically defined by a subset $\GG$ of the Grassmann
 bundle $G(p, TX)$ (as in Section 4.7) are always convex cone subequations.  As we have seen the
 universal subequations $\bbf(\GG)$ where $\GG$ is one of the grassmannians
 $
 G(p,\rn), G_\bbc(p, \bbc^n)\ \ {\rm or}\ \ G_\bbh(p, \bbh^n) \quad{\rm for\ } p=1,...,n
 $
 are particularly interesting.  These are principal branches of a polynomial equation $M_p=0$.
 However, there are many additional interesting examples. For some of them there is no known
 polynomial operator.
 
 We begin with the general result.  
 Fix a closed subset $\GG \ss G(p, \rn)$.
 Set $G=\{g\in {\rm O}_n : g(\GG)=\GG\}$
 and let $X$ be a riemannian manifold with topological
 $G$-structure so that the riemannian $G$-subequation $F(\GG)$ is defined on $X$.
 A domain $\O\ss\ss X$ is {\sl strictly $\GG$-convex} if it admits a strictly
 $\GG$-plurisubhamonic defining function.  
 The existence and topological structure of such domains has been studied in [HL$_{2,5}$]. 
 It is shown there that if $\bo$ is strictly $\GG$-convex 
 and if $\ob$ supports a strictly $\GG$-plurisubharmonic function,
 then $\O$ is itself strictly $\GG$-convex.
 From  Theorem 13.1 we have the following.
 
 \Theorem{16.1}  {\sl Let $X$ be a riemannian manifold with topological
 $G$-structure.
 Then for any strictly $\GG$-convex domain $\O\ss\ss X$  
 the Dirichlet problem for $F(\GG)$-harmonic functions 
 is uniquely solvable for all continuous boundary data.}
 \medskip

\Ex{16.2. (Calibrations)}  Fix a form $\phi\in \L^p \rn$ with comass 1, i.e., 
$$
\|\phi\|_{\rm comass} \ \equiv\ \sup  \left \{ \left|  {\phi|_W \over {\rm vol}_W}  \right|  :   W\in G(p, \rn) \right\}
\ \leq\ 1.
$$
Given  such a form , we define 
$$
\GG(\phi) \ \equiv\  \left \{ W\in G(p, \rn) : \left| {\phi|_W \over {\rm vol}_W}  \right| =1\right\}
$$
Let $G_\phi = \{g\in\ {\bf SO}_n : g^*\phi=\phi\}$ and note that $G_\phi$ preserves $\GG(\phi)$.
Therefore,  when $X$ has a topological $G_\phi$-structure,  
Theorem 16.1 applies. 
In this case the 
  $\GG(\phi)$-plurisubharmonic and $F)\GG(\phi))$-harmonic functions   are simply
  called $\phi$-{\sl plurisubharmonic} and $\phi$-{\sl harmonic} functions.
Specific examples are given at the end of the introduction.

  Note that if $X$ has a topological $G_\phi$-structure, the form  $\phi$ determines  a 
global smooth $p$-form $\phi$ on all of $X$ with  $\|\phi\|_{\rm comass}   \equiv 1$.
If $d\phi = 0$,  
then $\phi$ is a standard calibration [HL$_1$], and all $\phi$-{submanifolds} are automatically minimal.
For this case an analysis of $\phi$-plurisubharmonic functions, $\phi$-convexity, $\phi$-positive currents,
etc. is carried out in detail in [HL$_{2,3}$].

\Ex{16.3. (Lagrangian Plurisubharmonicity)}  Take $\GG\equiv \LAG \ss G_\bbr(n\, \bbc^n)$
to be the set of Lagrangian $n$-planes in $\bbc^n$. These planes have many equivalent descriptions.
They occur, for example,  as tangent planes to graphs of gradients over $\rn\ss\bbc^n$
and also as rotations of  $\rn\ss\bbc^n$ by an element of the unitary group U$_n$.
In particular,  the invariance group which fixes $\LAG$ is  ${\rm U}_n$.
There exists a polynomial $M$ on $\Sym(\bbr^{2n})$ which is $I$-hyperbolic and of degree
$2^n$ with the top branch of $\{M=0\}$ equal to $\bbf(\LAG)$.  This is discussed in a separate
paper [HL$_6$].

\Ex{16.4.  $\bbf(\GG(\f))$ as a Monotonicity Subequation}  
To be specific consider the associative calibration $\phi$ on $\bbr^7= {\rm Im} \bbo$
and the set $\GG(\phi)$ of associative 3-planes.
The convex cone subequation $\bbf(\GG(\f))$ is not the principal branch of a polynomial equation
$\{M=0\}$.  Nevertheless, one can construct subequations $\bbf$ which are  $\GG(\phi)$-monotone.
For example, define $\bbf$ by the condition that $A\in \bbf$ if 
$$
\exists \ \x,\eta\ \in\ \GG(\phi)\ \  {\rm with}\ \x\perp\eta\ \ {\rm and\ }  \ \tr_\x A \ \geq\ 0, \ \ \tr_\eta A\ \geq\ 0.
$$

 
 \vfill\eject
 

\centerline{\headfont 17. \ Equations Involving The Principal Curvatures of the Graph.}\medskip

Let $F_\bbs$ be one of the constant coefficient  subequations in $\rn$ discussed in Section 11.5.
This is defined by the requirement that for $C^2$-functions $u$, the principal curvatures
$(\kappa_1,...,\kappa_n)$ of the graph of $u$ lie in a given subset $\bbs\ss\rn$
which is symmetric and satisfies positivity. Now one can check that 
$$
\wt{F_\bbs}  \  =\ F_{\wt \bbs}
$$
 and in many interesting cases (such as those arising from G\aa rding hyperbolic 
 polynomials) we have $\bbs\ss\wt\bbs$ (or the reverse), and so $F_\bbs$-convexity implies $F_{\wt \bbs}$-convexity.
 From Theorem 12.7 and Proposition 11.17 we have the following.
 
 \Theorem{17.1} {\sl Let $F_\bbs$ be as above, and suppose $\O\ss\ss\rn$ is a smoothly bounded
 domain which satisfies the geometric convexity condition:
 $$
  (0, \l_1(x),...,\l_{n-1}(x)) \ \in\ \bbs \cap \wt \bbs \quad{\rm at\ each\ \ } x\in\bo,
 $$
 $\l_1(x),...,\l_{n-1}(x)$ are the  principal curvatures of the boundary at  $x$. Then existence holds
 for the Dirichlet problem for $F_\bbs$ for all continuous boundary data.  In fact the further statements
 in Theorem 12.7 hold.}

 \Remark{17.1}  For many interesting   equations involving the elementary
 symmetric functions $\s_k(\kappa) = \s_k(\kappa_1,.., \kappa_n)$ one also
 has uniqueness.  See the very nice paper [LE] for example.

 \vskip .3in
 
 

\centerline{\headfont 18. \ Applications of  Jet-Equivalence -- Inhomogeneous Equations.}\medskip

The notion affine jet-equivalence is flexible and powerful. Given Theorem  10.1
it is possible to treat a vast array of equations on manifolds.  One important case is the Calabi-Yau
equation examined in  6.15.  That example    can be greatly generalized as follows.

\Ex{18.1}  Let $F$ be the riemannian 
$G$-subequation for $G=$ O$_n$ (or U$_n$ or Sp$_n$) defined by
$$
\s_1(A+I)\ \geq \ 0, \  \s_2(A+I)\ \geq \ 0,\ ...\ , \s_{k-1}(A+I)\ \geq\ 0 \quad {\rm and} \quad \s_k(A+I)\geq 1
$$
where $\s_\ell(A)$ is the $\ell^{\rm th}$ elementary symmetric function in the eigenvalues of $A$
(or the complex or quaternionic hermitian symmetric part of $A$).
Let $\wt\Phi$ be the global affine jet-equivalence 
$$
\wt\Phi(r,p,A) \ =\ (r, p, (hI)A(hI)^t + (h^2-1)I)\ =\ (r,p, h^2 A+ (h^2-1)I)
$$
defined in (6.9) and set $f=h^{-2k}$. Let 
$F_f = \wt\Phi^{-1}(F)$.  Then a smooth $F_f$-harmonic function is one  which satisfies:
$$
\s_1(\Hess\, u+I)\ \geq \ 0,\ ...\ , \s_{k}(\Hess\, u+I)\ \geq\ 0  \qquad{\rm and}
$$
$$
\s_k(\Hess\, u+I)\ =\ f.
$$

This example immediately generalizes by replacing $\s_k(A)$ with any homogeneous
O$_n$-invariant polynomial on $\Symn$ which is hyperbolic with respect to the identity $I$
(see  18.4 below).

\Ex{18.2}  There is another (simpler) version of the examples above.
Let $F$ be the riemannian 
$G$-subequation for $G=$ O$_n$ (or U$_n$ or Sp$_n$) defined by
$$
\s_1(A)\ \geq \ 0, \ ...\ , \s_{k-1}(A)\ \geq\ 0 \quad {\rm and} \quad \s_k(A+I)\geq 1
$$
with $\s_\ell(A)$ as above. Let $\wt\Phi$ be the global affine jet-equivalence 
$
\wt\Phi(r,p,A) \ =\ (r, p, (hI)A(hI)^t)\ =\ (r,p, h^2 A)
$
and set $f=h^{-2k}$. Then the smooth $F_f = \wt\Phi^{-1}(F)$ harmonic functions satisfy
$$
\s_1(\Hess\, u)\ \geq \ 0,\ ...\ , \s_{k}(\Hess\, u)\ \geq\ 0  \qquad{\rm and}
$$
$$
\s_k(\Hess\, u)\ =\ f.
$$

\Ex{18.3.  (Inhomogeneous Equations for $\l_q(\Hess\,u)$)}  Let $F$ be a subequation on a manifold
$X$ and   $J$  any section of   $J^2(X)$.  Then $F_J \equiv F+J$ (fibre-wise translation)
is  a subequation affinely equivalent to $F$.  These two subequations have the same asymptotic 
interior (cf. \S 11).  Furthermore, any   monotonicity cone for one is a monotonicity cone for the other
(cf. \S 9).

A simple but  illustrative example,  for $X$ riemannian, is given by using the canonical splitting (4.3)
and taking $J_x = f(x)$Id where $f$ is a smooth
function on $X$ and Id is the identity section of $\Sym(T^*X)$.
Let $F=\cp_q$ the $q^{\rm th}$ branch of the Monge-Amp\`ere subequation (see Example 14.3).
Then the $F_J$-harmonic functions are solutions of the equation
$$
\l_q(\Hess\,u)\ =\ f(x)
$$
where $\l_q(\Hess\,u)$ is the $q^{\rm th}$ ordered eigenvalue of $\Hess\,u$.

Another  example comes from $F=F(\GG)$ as in section 4.7. Translating $F(\GG)$ by 
$J= {1\over p}f\cdot {\rm Id}$ gives the subequation
$$
\tr_\x( \Hess_x u )\ \geq\ f(x)    \fa \x\in \GG  \ {\rm and \ all\ } x\in X.
$$
The $F_J$-subharmonic functions are {\sl quasi-$\GG$-plurisubharmonic} with variable right hand side.

Further interesting examples are given by translating the subequations arising as branches of 
any G\aa rding $I$-hyperbolic polynomial on $\Symn$.

\Ex{18.4.  (G\aa rding Hyperbolic Polynomials)} The examples  above can be vastly generalized
by using G\aa rding's theory of hyperbolic polynomials [G].
(For details of what follows the reader is referred to [HL$_7$].)  
Let $M:\Symn \to \bbr$ be a homogeneous
polynomial of degree $\ell$  which is   hyperbolic  with respect to the identity $I$. 
Then the G\aa rding cone $\bbf^M$  is defined to be the connected component of $\{M>0\}$ containing $I$.
 It is a convex cone in $\Symn$ and is a monotonicity cone for the sets
 $\bbf^M(c) \equiv \{ A\in \bbf^M : M(A)\geq c\}$ and also for the branches $\bbf^M = \bbf^M_1\ss
 \bbf^M_2 \ss \bbf^M_3\ss \cdots$ described in the introduction.

Suppose now that this G\aa rding polynomial $M$
is invariant under a subgroup $G\subseteq  {\rm O}_n$ and that $X$ 
is a riemannian manifold with topological $G$-structure.  Then each of the sets above
determines a pure second-order subequation on $X$. 
For each of these,    the basic one $F^M$,  associated to $\bbf^M$,  is a 
monotonicity cone. Suppose now that $B$ is a smooth section of $F^M \ss \Sym(T^*X)$
and $h$ is a positive smooth function.  Consider the affine jet-equivalence $\Phi$ defined by
$\Phi(r,p,A) = (r,p, hI(A+B(x))hI)$ using the canonical splitting of  $J^2(X)$. Set 
$F_\Phi \equiv \Phi^{-1}(F^M(1))$ and let $f= h^{-2\ell}$. Then  one
finds as above that at any point $x\in X$
$$
A \in F_\Phi \qquad \iff\qquad A+B(x) \in F^M \ \ {\rm and}\ \ M(A+B(x))\geq f(x).
$$
 Then a smooth $F_\Phi$-harmonic function $u$ is one that satisfies
 $$
\Hess_x u +B(x) \in F^M \qquad {\rm and}\qquad  M(\Hess_x u +B(x)) = f(x)
 $$
 
 One can also translate the various branches $F^M_k$ of the main subequation
 by the section $B$.

\vfill\eject


\centerline{\headfont  19.  \ Equations of Calabi-Yau type in the Almost Complex Case.}\medskip

Let $X$ be an almost complex hermitian manifold and consider the subequation 
$$
\Hess_\bbc u + I   \ \geq\ 0 \and \det_\bbc  \left(\Hess_\bbc u + I\right) \ \geq \  F(u)f(x)
\eqno{(19.1)}
$$
with $F, f > 0$ and $F$ non-decreasing (e.g. $F(u)=e^u$) (cf. Example 6.15 and \S 15).
Let $\O\ss\ss X$ be a domain with smooth boundary $\bo$.  To address the Dirichlet problem for 
(19.1) on $\O$ we need to analyze the condition of $F$-convexity for $\bo$. 
For this we fix $\l\in\bbr$ and consider the subequation $F_\l$ (independent of the $r$-variable)
given in $(p,A)$-coordinates by
$$
A_\bbc   + I   \ \geq\ 0 \and \det_\bbc  \left(A_\bbc  + I\right) \ \geq \  F(\l)f(x).
\eqno{(19.2)}
$$
Recall the complex positivity cone $\cp_\bbc = \{A: A_\bbc >0\}$.

\Def{19.1}  The domain $\O\ss X$ is  {\bf strictly pseudo-convex} if there exists a defining function
$\rho$ for $\O$ which is strictly $\cp_\bbc$-subharmonic on $\ob$
 (i.e., $\Hess_\bbc \rho>0$ on  $\ob$).

\Lemma {19.2} {\sl
The boundary of a   strictly pseudo-convex domain is strictly $F$  and $\ft$ convex.
}
\pf
One computes directly that the open cone $\Int \cp_\bbc$ lies is the asymptotic
interior $\overrightarrow{F_\l}$ for any value of $\l$. One also sees directly that $\Int \cp_\bbc$ 
lies is the asymptotic interior $\oa{\ft}_\l$ for any value of $\l$. (Details are left as
an exercise.) The assertion now follows from Definition 11.10.\qed

\Theorem{19.3}  {\sl
Let $\O$ be a strictly pseudo-convex domain in an almost complex hermitian manifold $X$.
Then the Dirichlet problem for the Calabi-Yau equation (19.1) is uniquely solvable for 
any continuous function on the boundary $\bo$.
}
\pf
First recall from Example 6.15 that $F$ is locally affinely jet-equivalent to a constant coefficient 
subequation on $X$.  It is straightforward to see that $\cp_\bbc$ is a monotonicity cone
for $F$, and by hypothesis there exists a global defining  function for $\O$ which is  strictly 
$\cp_\bbc$-subharmonic.  We have already seen that $\bo$ is strictly $F$ and $\ft$ convex.
Hence Theorem 13.1$'$ applies. \qed

\medskip

Admittedly the Calabi-Yau equation holds more interest in the (integrable) complex
manifold case. However, it is surprising that the Dirichlet problem is uniquely solvable 
in this very general setting.

\Ex{19.4}  Many examples can be obtained by starting with a strongly pseudo-convex domain 
$\O_0$ in a complex manifold $X_0$ and then deforming the complex structure $J_0$ to a nearby
almost complex structure which is not integrable. Any given hermitian metric on $X_0$ 
can be continuously deformed by averaging over the nearby  $J$'s. For small enough 
deformations any given strongly plurisubharmonic defining function for $\O_0$ will
remain $\cp_\bbc$-subharmonic. Thus, the Dirichlet problem for (19.1) 
can be solved for all sufficiently nearby structures.

\vfill\eject


\centerline{\headfont Appendix A. Equivalent Definitions  F-Subharmonic.}\medskip

In this appendix we assume the following  for each fibre
$F_x\ss \bbJt = \bbr\times\rn\times \Symn$. 
The Positivity Condition is not required, nor is $F$ required to be closed.
Our assumption  is that for some $P>0$:
$$
{\rm  If} \ (r,p,A+\a  P)\in F_x \fa \a>0  \ \ {\rm   small,\  then\ \ }
(r,p,A)\in F_x.  
\eqno{(A.1)}
$$

 \medskip\noindent
 {\bf   Remark.}  Assuming the Positivity Condition
(P), Condition (A.1) is equivalent to requiring that  {\sl  the $J^1$-fibres of $F$,
that is, the fibres of $F$ under the natural map $J^2(X)\to J^1(X)$,  are closed. }
In terms of the standard coordinates this means that
 $$
 {\rm Each\ \ } F_{x,r,p}\ =\ \{ A\in\Symn : (x,r,p,A) \in F\}\ \ {\rm is\ closed.}
 \eqno{(A.1)'}
 $$
 For the proof note that if $A\in {\overline F}_{x,r,p}$,
 then for each $\e>0$ there exists $A_\e\in F_{x,r,p}$ with $A_\e-A\leq \e P$.
 By (P) this implies $A+ \e P \in F_{x,r,p}$.  Assuming (A.1) this proves that 
 $A\in F_{x,r,p}$, i.e., this proves (A.1)$'$. The converse obviously holds for all $P>0$.
 Thus, assuming positivity, (A.1) holds for one $P>0$ if and only if it holds for all $P>0$.

 \Prop{A.1}  {\sl  Suppose   that $u\in\USC(X)$ and $x_0\in X$.   Let  $x=(x_1,...,x_n)$ be local
 coordinates on a neighborhood of $x_0$.
 Then the following Conditions I, II, III, and IV are equivalent.}
 
 \medskip
 
 \item {I.}  For all $\vf \in C^2$ near $x_0$, 
 $$
(1) \ \  \left\{ 
 \eqalign
{  
u-\vf \  &\leq \ 0\quad{\rm near}\ x_0  \cr
  \  &= \ 0\quad \ \ {\rm at}\ x_0  
}
\right\} \quad\Rightarrow \quad \jt_{x_0}\vf\, \in \, F_{x_0}
$$

 \medskip
 
 \item {II.}  For all  $(r,p,A)\in \bbJt$, 
 $$
(2) \ \  \left\{ 
\eqalign
{
u(x)-\bigl[r+\langle p, x-x_0\rangle +\half \langle A(x-x_0), x-x_0 \rangle\bigr] \ &\leq \ 0  \ \ 
{\rm near\ } x_0   \cr
&=\ 0 \ \  {\rm at}\ \ x_0
}
\right\} 
 $$
$$
\Rightarrow \quad(r,p,A) \in F_{x_0}
$$

 \medskip
 
 \item {III.}  For all  $(r,p,A)\in \bbJt$, 
 $$
(3)\ \  \left\{ 
\eqalign
{
u(x)-\bigl[r+\langle p, x-x_0\rangle +\half \langle A(x-x_0), x-x_0 \rangle
 + o(|x-x_0|^2)\bigr] \ &\leq \ 0 \ \ 
{\rm near\ } x_0   \cr
&=\ 0 \ \  {\rm at}\ \ x_0
}
\right\} 
 $$
$$
\Rightarrow \quad(r,p,A) \in F_{x_0}
$$

 \medskip
 
 \item {IV.}  For all  $(r,p,A)\in \bbJt$ and $\a>0$, 
 $$
(4)\ \  \left\{ 
\eqalign
{
u(x)-\bigl[r+\langle p, x-x_0\rangle +\half \langle A(x-x_0), x-x_0 \rangle\bigr] \ &\leq \ -\a|x-x_0|^2 \ \ 
{\rm near\ } x_0   \cr
&=\ 0 \ \ \ \  \qquad\qquad  {\rm at}\ \ x_0
}
\right\} 
 $$
$$
\Rightarrow \quad(r,p,A) \in F_{x_0}
$$
\pf
\ (I $\Rightarrow$ II):  Given $(r,p,A)\in \bbJt$ satisfying (2), set 
\smallskip
\centerline{$\vf =  r+\langle p, x-x_0\rangle +\half \langle A(x-x_0), x-x_0 \rangle$.}
\smallskip
\noindent
Since the quadratic function $\vf$ satisfies (1), the condition I implies 
$(r,p,A)= \jt_{x_0}\vf \in F_{x_0}$.\medskip

\noindent (II $\Rightarrow$ III):  Given $(r,p,A)\in \bbJt$ satisfying (3),  it follows that $\forall\, \a>0$
$$
\eqalign
{
u(x)-\bigl[r+\langle p, x-x_0\rangle +\half \langle A(x-x_0), x-x_0 \rangle\bigr] \ &\leq \ \a|x-x_0|^2 \ \ 
\  {\rm near\ } x_0   \cr
&=\ 0 \ \ \  \qquad\qquad  {\rm at}\ \ x_0
}
$$
or equivalently
$$
\eqalign
{
u(x)-\bigl[r+\langle p, x-x_0\rangle +\half \langle (A+2\a I)(x-x_0), x-x_0 \rangle\bigr] \ &\leq \ 0
\qquad\qquad \ \  {\rm near\ } x_0   \cr
&=\ 0 \ \   \qquad\qquad  {\rm at}\ \ x_0.
}
$$
However, this is just (2) for $(r,p, A+2\a I)$.  Hence by II 
we have $(r,p, A+2\a I)  \in F_{x_0}$ for all $\a>0$, proving
that $(r,p, A)  \in F_{x_0}$ because of the assumption (A.1).
 \medskip

\noindent 
 (III $\Rightarrow$ I):  Given $\vf$ satisfying (1), the Taylor series for $\vf$ satisfies (3).
 \medskip

\noindent 
 (II $\Rightarrow$ IV):  This holds since (4) $\Rightarrow$ (2).
 \medskip

\noindent 
 (IV $\Rightarrow$ II):  Suppose that $(r,p,A)\in \bbJt$ satisfies (2).  Equivalently,
 $$
 \eqalign
{
u(x)-\bigl[r+\langle p, x-x_0\rangle +\half \langle (A+2\a I)(x-x_0), x-x_0 \rangle\bigr] \ &\leq \ -\a|x-x_0|^2 \ \ 
{\rm near\ } x_0   \cr
&=\ 0 \ \ \ \  \qquad\qquad  {\rm at}\ \ x_0
}
$$
That is, $(r,p,A +2\a I)$ satisfies (4).
Therefore, by IV, $(r,p,A + 2\a I)\in F_{x_0}$ for all $\a>0$.
Finally (A.1) implies that 
$(r,p,A)\in F_{x_0}$. \qed

\vfill\eject


\centerline{\headfont Appendix B. Elementary Properties of  F-Subharmonic Functions.}\medskip

The proof of Theorem 2.6, which lists  the elementary properties of $F$-subharmonic 
functions, is given in this appendix.  As explained in Remark 2.6, the positivity condition
is not needed in these proofs.
For Property (A) only condition A.1 is  required.
For Property (B) we only require that the $J^1$-fibres of $F$ be closed (cf. (A.1)$'$).
Not surprisingly, for (C), (D) and (E) the full hypothesis that $F$ be closed is used.

\medskip
\noindent
{\bf Proofs.}   

\noindent
(A)\ The condition that $\max\{u,v\} -\vf\leq 0$ near $x_0$ with equality at $x_0$ 
implies that for one of the 
functions $u,v$, say $u$, we have $u(x_0)=\vf(x_0)$.  In this case, $u-\vf\leq 0$ near $x_0$ 
 with equality at $x_0$.  Hence, $\jt_{x_0}\vf \in F_{x_0}$.

\noindent
(B)\  This follows from Condition III in Proposition A.1.

The remaining properties are proved by a common method which uses Lemma 2.4.

\noindent
(C)  Recall the basic fact that if 
$\{v_j\}\ss\USC(X)$ is a decreasing sequence with limit $v$, 
then
$$
\lim_{j\to\infty}  \biggl\{ \sup_K v_j\biggr\} \ =\ \sup_K v
\eqno(B.1)
$$
This is proven as follows.  Given $\d>0$, the upper semi-continuity of $v_j$ implies that
the set $K_j = \{x\in K : v_j(x) \geq \sup_K v +\d\}$ is compact.
The sets $K_j$ are decreasing since $\{v_j\}$ is decreasing.  The pointwise convergence 
of $\{v_j\}$ to $v$ implies that $\bigcap_j K_j =\emptyset$.  Hence, there exists $j_0$ with
 $K_j =\emptyset$ for all $j\geq j_0$.  That is, $v_j(x) < \sup_K v+ \d$ for all $j\geq j_0$ and $x\in K$.

Suppose now that $u\notin F(X)$.  Then by Lemma  2.4 there exists $x_0\in X$,
 local coordinates $x$  about $x_0$, $\a>0$, and a quadratic function 
 $\vf(x) = r+\langle p, x-x_0\rangle +\half \langle A(x-x_0), x-x_0 \rangle$ such that 
$$
\eqalign
{
u(x)- \vf(x) \ &\leq \ -\a|x-x_0|^2 \qquad
{\rm near\ } x_0 \qquad{\rm and}  \cr
&=\ 0\qquad\qquad\qquad\ \ \ {\rm at}\ \ x_0
}
\eqno{(B.2)}
$$
but $J^2_{x_0} \vf = (r,p, A) \notin F_{x_0}$.

Pick a small closed ball $B$ about $x_0$ so that $u-\vf$ has a strict maximum over $B$
at $x_0$.   Choose a maximum point $x_j \in B$ for the function $u_j-\vf$ over $B$.
Fix a  smaller open ball $B'$ about $x_0$ and let $K = B-B'$ denote the corresponding compact 
annulus.  Apply (B.1) to the decreasing sequence $v_j \equiv u_j-\vf$.
Since $\sup_K(u-\vf) < 0$, this proves that  $\sup_K (u_j-\vf) <0$ for all $j$ sufficiently large.

Since the maximum value of $u_j-\vf  \geq u-\vf$ over $B$ is $\geq0$, this proves that $x_j \notin K$ i.e., $x_j \in B'$ for $j$ large.  Since $B'$ was arbitrary, we have
$$
\lim_{j\to\infty} x_j \ =\ x_0.
\eqno{(B.3)}
$$
In particular, $x_j\in\Int B$ is an interior maximum point  for $u_j-\vf$. Set 
$$
r_j = u_j(x_j),\ \ \ p_j = D_{x_j}\vf = p+A(x_j-x_0), \  \ \ A_j = D^2_{x_j}\vf = A.
$$
This proves that 
$$
(r_j, p_j, A_j) \in\ F_{x_j}
$$
since $u_j\in F(X)$.

Applying (B.1) to $K=B$ yields that 
$r_j = \sup_B  u_j \searrow \sup_B u = r$.
Hence, $\lim_{j\to\infty} r_j = r$. This proves that
$$
\lim_{j\to\infty} (x_j,  r_j, p_j, A_j) \ =\ (x_0, r, p, A).
\eqno{(B.4)}
$$
Since $F$ is closed, this proves that $(x_0, r, p, A) \in F_{x_0}$, contrary to hypothesis.
\medskip

\noindent
(D)\ Since (B.1) is valid if $\{v_j\}$ converges uniformly to $v$, the proof of (D) is essentially 
the same as the proof of (C), except easier.
\medskip

\noindent
(E) \ Suppose $u\notin F(X)$. Then exactly as in the previous proofs we have (B.2) for $v^*$.
Since 
$$
r  \ =\ v^*(x_0)  \ = \ \lim_{k\to\infty} \ \sup_{|y-x_0|\leq {1\over k}}\left\{ \sup_{f \in \cf} f(y)
\right\},
$$
it follows easily that there exists a sequence $y_k\to x_0$ in $\rn$ and a sequence
$f_k\in \cf$ such that
$$
\lim_{k\to\infty} f_k(y_k)\ =\ r.
$$
Choose a maximum point $x_k$ for $f_k-\vf$ over $B$, a small closed ball about $x_0$
on which the condition (B.2) holds.
By taking a subsequence we may assume that $x_k\to \bar x\in B$.  Now
$$
f_k(y_k) -\vf (y_k)\ \leq \ f_k(x_k) -\vf (x_k).
\eqno{(B.5)}
$$
The left hand side has limit zero.  Hence,
$$
\eqalign
{
0 \  &\leq \liminf_{k\to \infty}f_k(x_k) -\vf (\bar x)  \cr
&\leq \limsup_{k\to \infty}v(x_k) -\vf (\bar x) \ \leq\ v^*(\bar x) -\vf(\bar x).\cr
}
\eqno{(B.6)}
$$
since by definition $f_k\leq v$.  This proves that $\bar x=x_0$ by (B.2). In particular, each $x_k$ is interior to $B$
for $k$ large.  Therefore,  since each $f_k$ is $F$-subharmonic, we find that
$$
\left( f_k(x_k) , \ D_{x_k}\vf,\ D_{x_k}^2 \vf\right)\ \in\ F_{x_k}.
$$
Note that (B.6) implies $\lim_{k\to\infty}f_k(x_k) = v^*(x_0) = \vf(x_0)$.
Since $F$ is closed we conclude that
$$
\jt_{x_0}\vf = \lim_{k\to\infty}\left( f_k(x_k) , \ D_{x_k}\vf,\ D_{x_k}^2 \vf \right) \  \in \ F_{x_0},
$$
which is a contradiction.\qed

 \vfill\eject


\centerline{\headfont \ Appendix C.\   The Theorem on Sums.}
\medskip

 In this appendix we recall the fundamental Theorem on Sums, which
 plays a key role in viscosity theory  ([CIL], [C]).  We restate the result in a form 
 which is particularly suited to our use.
 
 Fix an open subset $X\ss\rn$, and let $F, G\ss J^2(X)$
 be two second order partial differential subequations.
   Recall that a function $w\in\USC(X)$ satisfies the {\sl Zero Maximum Principle}
on a compact subset $K\ss X$  if
$$
w\ \leq \ 0 \ \ {\rm on\ } \partial K \qquad\Rightarrow\qquad w\ \leq \ 0 \ \ {\rm on\ } K.
\eqno{(ZMP)}
$$

 \Theorem {C.1}  {\sl   Suppose $u\in F(X)$ and $v\in G(X)$, but $u+v$ does not 
 satisfy the Zero Maximum Principle (ZMP) on a compact set $K\ss X$.
 Then there exist a point $x_0\in \Int K$ and a sequence of numbers $\e\searrow 0$
 with associated points $z_\e = (x_\e, y_\e) \to (x_0, x_0)$ in $X\times X$, and there exist
 $$
 J_{x_\e} \ \equiv \ (r_\e, p_\e, A_\e)\ \in\  F_{x_\e}
 \and 
 J_{y_\e} \ \equiv \ (s_\e, q_\e, B_\e)\ \in \ G_{y_\e} 
 $$
 such that
   $$r_\e \  =\ u(x_\e),\ \ \  s_\e\ =\ v(y_\e), \and
  r_\e + s_\e \ =\ M_\e \ \searrow \ M_0 >0,
  \leqno{(1)}$$     
 $$ p_\e\ =\ {x_\e-y_\e\over \e},\ \ \   q_\e\ =\ {y_\e-x_\e\over \e} \ =\ -p_\e, 
 \and {|x_\e-y_\e|^2\over \e    } \ \to\ 0
 \leqno{(2)}$$  
  $$-{3\over\e} I\ \leq\ \left( \matrix{
 A_\e & 0 \cr 0 & B_\e\cr}  \right)
 \ \leq\ {3\over\e} \left( \matrix{
 I & -I \cr -I & I\cr}  \right)
\leqno{(3)}$$
 }

\Remark{C.2}   In fact, we have
$$
 J_{x_\e} \ \in\ \overline{J_{x_\e}^{2,+} }u
 \and
  J_{y_\e} \ \in\ \overline{J_{y_\e}^{2,+} }v,
 $$
where $J_{x}^{2,+} u$ denotes the upper 2-jet of $u$ at $x$.
\medskip  

Restricting the right hand inequality in (3) to the diagonal yields
$$
A_\e+B_\e \ =\ -P_\e\qquad{\rm where\ \ } P_\e\ \geq\ 0.
\leqno{(3)'}
$$
This is enough to prove a result which lies between weak comparison and comparison 
for  a constant coefficient subequation  $\bbf$ on $\rn$.

\Cor{C.3} {\sl Suppose $\bbh$ and $\bbf$ are constant coefficient subequations with
$\bbh\ss \Int \bbf$.
Suppose $K$ is a compact subset of $\rn$.  If $u\in \bbh(K)$ and $v\in \wt\bbf(K)$, then}
$$
u+v\ \leq \ 0 \quad {\rm on}\ \partial K\qquad\Rightarrow\qquad u+v\ \leq \ 0  \quad {\sl on}\  K
$$

\pf By (1), (2), and (3)$'$ we have
$$
(-s_\e, -q_\e, -B_\e) \ =\ (r_\e-M_\e, p_\e, A_\e+P_\e).
$$

Now $(r_\e, p_\e, A_\e)\in \bbh$ implies $(r_\e-M_\e, p_\e, A_\e+P_\e)\in\Int \bbf$.

However, $(s_\e, q_\e, B_\e)\in \wt{\bbf}$ says that $(-s_\e, -q_\e, -B_\e)\notin \Int \bbf$.\qed

 \vfill\eject


\centerline{\headfont \ Appendix D.\   Some Important Counterexamples.}
\medskip

One might  wonder whether comparison, or at least uniqueness for the 
Dirichlet problem, can be established if one weakens the assumption that there exist a 
strictly $M$-subharmonic function where $M+F\ss F$.
Consider for example the case of domains which are strictly $F$-convex
( i.e., for which there exists a  globally strictly $F$-subharmonic defining function)  but  the condition
 $F+F\ss F$ is not satisfied.  
 We give here an   example of a strongly $P_3$-convex
 domain in a non-negatively curved space  where comparison  and, in fact,
 uniqueness for the Dirichlet problem fail.
 
 Consider the standard riemannian 3-sphere $S^3 = \{x\in \bbr^4 : \|x\|=1\}$and the great 
 circle $\g = \{(x_1, x_2, 0,0) : x_1^2+x_2^2=1\} \ss S^3$.  Define
 $$
 \d : S^3 \ \to \bbr \quad{\rm by}\quad  \d(x) \ \equiv\ \dist(x, \g)
 $$
where distance is taken in the 3-sphere metric.  The level sets $T_s \equiv \d^{-1}(s)$
for $0<s< \pi/2$ are flat tori, which are orbits of the obvious $T^2$ torus action
on $\bbr^4=\bbr^2\times\bbr^2$. The eigenvalues $\l_1(s), \l_2(s)$ of the second fundamental
form $II_s$ of $T_s$ are constant on $T_s$ and,
by  the Gauss curvature equation, satisfy 
$$
\l_1(s) \ =\ -{1\over \l_2(s)}.
$$
Straightforward calculation (cf. [HL$_2$, (5.7)]) shows that the riemannian Hessian
$$
\Hess\, \d \ =\ \left ( \matrix{0&0\cr 0& II_s\cr}\right)  \qquad {\rm where \ \ }  \d = s  \in (0,\pi/2).
$$
from which it follows easily that
$$
\Hess\, \left ( \smfrac 1 2 \d^2\right) \ =\ \left ( \matrix{1&0\cr 0& \d \cdot II_s\cr}\right) 
\ =\  \left ( \matrix{1&0 &0\cr 0& s \l(s) & 0\cr 0&0& -{ s\over \l_(s)}}\right) 
$$
where $\d = s  \in (0,\pi/2)$. As $s\searrow 0$, $\l(s) = {1\over s} + O(1)$, and so
$$
\Hess \,  \left ( \smfrac 1 2 \d^2\right)\bigr|_{s=0}  \ =\ 
 \left ( \matrix{1&0 &0\cr 0& 1 & 0\cr 0&0&0}\right) 
$$
In particular we see the following. 
 Let 
 $$
 U\equiv S^3 - \wt \g 
 $$
  where $\wt\g = \{(0,0,x_3, x_4) : x_3^2+x_4^2=1\} $ is the ``opposite'' or ``focal'' geodesic  to $\g$.
  Then
$$
\Hess\, \left ( \smfrac 1 2 \d^2\right) \ \ 
{\rm has\  2\  strictly\  positive\  eigenvalues\  everywhere\  on \ } U, \ {\rm and}
\eqno{(D.1)}
$$
$$
\Hess\, \left ( - \smfrac 1 2 \d^2\right) \ \ 
{\rm has\  1\     eigenvalue\ } \geq 0 {\ \rm everywhere\  on \ } U.
\eqno{(D.2)}
$$

\medskip\noindent
{\bf Conclusion D.1. (Co-convex $\not\Rightarrow$ Maximum Principle)}.   Example (D.2) shows that on 
spherical domains the Maximum
Principle fails for co-convex, i.e. $\cpt$-subharmonic functions 
(where $\cpt$ means at least one eigenvalue $\geq 0$).   
 In euclidean space co-convex functions do 
  satisfy the maximum principle. They are called
{\sl subaffine functions} and play an important role in [HL$_4$]. \medskip

Consider now the product
$$
U\times U \ \ss \  S^3\times S^3
$$
and  define $\d_k = \d\circ \pi_k$ where $\pi_k : U\times U \to U$ denotes projection onto 
the $k$th factor.  Set $$\rho \equiv \smfrac 1 2  \d_1^2 +  \smfrac 1 2  \d_2^2$$
and note that 
$$
\Hess\, \rho \ =\ \Hess\, \left ( \smfrac 1 2 \d_1^2\right) \oplus \Hess\, \left ( \smfrac 1 2 \d_2^2\right).
$$
In particular, $\rho$ has four strictly positive eigenvalues everywhere on $U\times U$.
In the terminology of \S 7 this means that $\rho$ {\bf  is strictly $P_2$-subharmonic on $U\times U$.}
(Recall that $P_q$-subharmonic means that  $\Hess\,u$ 
has at least $n-q=6-q$ eigenvalues $\geq0$.)
Therefore the domain
$$
\O_c \ \equiv\ \{(x,y)\in U\times U: \rho(x,y)\leq c  \}\ \ {\rm is\ strictly\ } P_2 \ {\rm convex}
\eqno{(D.3)}
$$
for $0<c <{\pi^2\over 8}$.

Consider now the functions $$u_1 = -\smfrac 1 2 \d_1^2 \and u_2 = -\smfrac 1 2 \d_2^2.$$
By (D.2) have the following. Recall that $u$ is $P_q$-harmonic if $\l_{q+1}\equiv0$ where
$\l_1\leq \cdots \leq \l_6$ are the eigenvalues of $\Hess\,u$.

\Lemma{D.2}  {\sl Each $\Hess\, u_k$, $k=1,2$  has  two negative eigenvalues,  three zero eigenvalues
and one non-negative eigenvalue at every point of $U\times U$. In particular, 
$$
u_1, u_2 \ \in\ P_2(U\times U) \and -u_1, -u_2 \ \in\  P_1(U\times U).
$$
Furthermore, each $u_k$ is $P_2$, $P_3$ and $P_4$-harmonic, whereas
each $-u_k$ is $P_1$, $P_2$ and $P_3$-harmonic on $U\times U$.
}

\medskip\noindent
{\bf Conclusion D.3. (Comparison Fails)}.  Each domain $\O_c$  is 
strictly $P_2$-convex. Furthermore there are  smooth functions
$$
u_1\ \in\ P_2(\O_c)\and u_2\ \in\ \wt P_2(\O_c)
$$
(since $P_2(\O_c) \ss   P_3(\O_c) = \wt P_2(\O_c)$) such that
$$
u_1 + u_2 \ \equiv\  -c \ <\ 0 \ \ {\rm on\ }\partial \O_c \and 
\sup_{\O_c} (u_1 + u_2)\ =\ 0.
$$
Hence,  the Comparison Principle fails for $P_2$ on each $\O_c$

\medskip\noindent
{\bf Conclusion D.4. (Uniqueness fails for the Dirichlet problem)}.   Each domain $\O_c$
is  strictly  $P_3$-convex (since $P_2\ss P_3$, and so $P_2$-subharmonic
$\Rightarrow$ $ P_3$-subharmonic).  By the lemma
we have that 
$$
{\rm both\ functions\ \ } u_1 \ \ {\rm and}\ \  (c-u_2) \  \ {\rm are \ } P_3\!-\!\!{\ \rm harmonic\ on }\ \O_c,
$$
and 
$$
u_1 \ \equiv \  c-u_2 \ \  {\rm on\ \ } \partial \O_c.
$$

\medskip\noindent
{\bf Remark  D.5. (Existence without uniqueness (again))}.   Note that the domains $\O_c$
are  strictly  $\oa {P_2}$ and $\overrightarrow { \wt {P_2}}$-convex (since $P_k = \oa {P_k} $
and $P_2\ss  P_3$).   Since the isometry group of $S^3\times S^3$ is transitive and preserves
these subequations, Theorem 12.5 guarantees 
the  existence of solutions to the Dirichlet Problem for all continuous boundary data on each
$\O_c$.  However, as  seen above, uniqueness fails in general.

\vfill\eject

\def\b{{\beta}}



\centerline{\bf References}

\vskip .2in

\noindent
\item{[A$_1$]}   S. Alesker,  {\sl  Non-commutative linear algebra and  plurisubharmonic functions  of quaternionic variables}, Bull.  Sci.  Math., {\bf 127} (2003), 1-35. also ArXiv:math.CV/0104209.  

\smallskip

\noindent
\item{[A$_2$]}   S. Alesker,  {\sl  Quaternionic Monge-Amp\`ere equations}, 
J. Geom. Anal., {\bf 13} (2003),  205-238.
 ArXiv:math.CV/0208805.  

\smallskip

\noindent
\item{[AV]}    S. Alesker and M. Verbitsky,  {\sl  Plurisubharmonic functions  on hypercomplex manifolds and HKT-geometry}, arXiv: math.CV/0510140  Oct.2005

\smallskip

\noindent
\item{[AFS]}  D. Azagra, J. Ferrera and B. Sanz, {\sl Viscosity solutions to second order partial differential
equations on riemannian manifolds}, ArXiv:math.AP/0612742v2,  Feb. 2007.
 \smallskip

\noindent
\item{[BT]}   E. Bedford and B. A. Taylor,  {The Dirichlet problem for a complex Monge-Amp\`ere equation}, 
Inventiones Math.{\bf 37} (1976), no.1, 1-44.

\smallskip

 \item{[B]}  H. J. Bremermann,
    {\sl  On a generalized Dirichlet problem for plurisubharmonic functions and pseudo-convex domains},
          Trans. A. M. S.  {\bf 91}  (1959), 246-276.
\medskip

 \item{[BH]}  R. Bryant and F. R. Harvey,
    {\sl  Submanifolds in hyper-K\"ahler geometry},
          J. Amer. Math. Soc. {\bf 1}  (1989),  1-31.
\medskip

\noindent
 \item{[CNS]}   L. Caffarelli, L. Nirenberg and J. Spruck,  {\sl
The Dirichlet problem for nonlinear second order elliptic equations, III: 
Functions of the eigenvalues of the Hessian},  Acta Math.
  {\bf 155} (1985),   261-301.

 \smallskip

\noindent
\item{[C]}   M. G. Crandall,  {\sl  Viscosity solutions: a primer},  
pp. 1-43 in ``Viscosity Solutions and Applications''  Ed.'s Dolcetta and Lions, 
SLNM {\bf 1660}, Springer Press, New York, 1997.

 \smallskip

\noindent
\item{[CCH]}  A. Chau, J. Chen and W. He,  {\sl  Lagrangian mean curvature flow for entire Lipschitz graphs},  
ArXiv:0902.3300 Feb, 2009.

 \smallskip

\noindent
\item{[CGG$_1$]}  Y.-G. Chen, Y. Giga and S. Goto  {\sl  Uniqueness and existence of viscosity solutions 
of generalized mean curvature flow equations},  Proc. Japan Acad. Ser. A. Math. Sci  {\bf 65} (1989), 207-210.

 \smallskip

\noindent
\item{[CGG$_2$]}  Y.-G. Chen, Y. Giga and S. Goto  {\sl  Uniqueness and existence of viscosity solutions 
of generalized mean curvature flow equations},  J. Diff. Geom. {\bf 33} (1991), 749-789..

 \smallskip

\noindent
\item{[CIL]}   M. G. Crandall, H. Ishii and P. L. Lions {\sl
User's guide to viscosity solutions of second order partial differential equations},  
Bull. Amer. Math. Soc. (N. S.) {\bf 27} (1992), 1-67.

 \smallskip

\noindent
\item{[DK]}   J. Dadok and V. Katz,   {\sl Polar representations}, 
J. Algebra {\bf 92} (1985) no. 2,
504-524.

\smallskip

\noindent
\item{[E]}   L. C. Evans,  {\sl   Regularity for fully nonlinear elliptic equations
and motion by mean curvature},  pp. 98-133 
in ``Viscosity Solutions and Applications''  Ed.'s Dolcetta and Lions, 
SLNM {\bf 1660}, Springer Press, New York, 1997.

 \smallskip

\noindent
\item{[ES$_1$]}   L. C. Evans and J. Spruck,  {\sl   Motion of level sets by mean curvature, I},  
J. Diff. Geom. {\bf 33}  (1991), 635-681.

 \smallskip

\noindent
\item{[ES$_2$]}   L. C. Evans and J. Spruck,  {\sl   Motion of level sets by mean curvature, II},  
Trans. A. M. S.  {\bf 330}  (1992),  321-332.

 \smallskip

\noindent
\item{[ES$_3$]}   L. C. Evans and J. Spruck,  {\sl   Motion of level sets by mean curvature, III},  
J. Geom. Anal.   {\bf 2}  (1992),  121-150.

 \smallskip

\noindent
\item{[ES$_4$]}   L. C. Evans and J. Spruck,  {\sl   Motion of level sets by mean curvature, IV},  
J. Geom. Anal.   {\bf 5}  (1995),   77-114.

 \smallskip

\noindent
\item{[G]}   L. G\aa rding, {\sl  An inequality for hyperbolic polynomials},
 J.  Math.  Mech. {\bf 8}   no. 2 (1959),   957-965.

 \smallskip

\noindent
\item{[Gi]}   Y. Giga, {\sl  Surface Evolution Equations -- A level set approach},
Birkh\"auser,  2006.

 \smallskip

 \noindent 
\item {[H]}   F. R. Harvey,  { Spinors and Calibrations},  Perspectives in Math. vol.9, Academic Press,
Boston, 1990.

 \smallskip

 \noindent 
\item {[HL$_1$]}   F. R. Harvey and H. B. Lawson, Jr,  {\sl Calibrated geometries}, Acta Mathematica 
{\bf 148} (1982), 47-157.

 \smallskip

\item {[HL$_{2}$]} F. R. Harvey and H. B. Lawson, Jr., 
 {\sl  An introduction to potential theory in calibrated geometry}, Amer. J. Math.  {\bf 131} no. 4 (2009), 893-944.  ArXiv:math.0710.3920.

\smallskip

\item {[HL$_{3}$]} F. R. Harvey and H. B. Lawson, Jr., {\sl  Duality of positive currents and plurisubharmonic functions in calibrated geometry},  Amer. J. Math.    {\bf 131} no. 5 (2009), 1211-1240. ArXiv:math.0710.3921.

\smallskip

\item {[HL$_{4}$]}  F. R. Harvey and H. B. Lawson, Jr., {\sl  Dirichlet duality and the non-linear Dirichlet problem},    Comm. on Pure and Applied Math. {\bf 62} (2009), 396-443.

\smallskip

\item {[HL$_{5}$]} F. R. Harvey and H. B. Lawson, Jr.,  {\sl  Plurisubharmonicity in a general geometric context},  Geometry and Analysis {\bf 1} (2010), 363-401. ArXiv:0804.1316.

\smallskip

\item {[HL$_{6}$]} F. R. Harvey and H. B. Lawson, Jr., {\sl  Lagrangian plurisubharmonicity and convexity},  Stony Brook Preprint (2009).

\smallskip

\item {[HL$_{7}$]} F. R. Harvey and H. B. Lawson, Jr.,  {\sl  Hyperbolic polynomials and the Dirichlet problem},   ArXiv:0912.5220.
\smallskip

\item {[HL$_{8}$]} F. R. Harvey and H. B. Lawson, Jr.,  {\sl  Potential Theory on almost complex manifolds}, (in preparation).
\smallskip

%
   %
%
%

   \noindent
\item{[HM]}    L. R. Hunt and J. J. Murray,    {\sl  $q$-plurisubharmonic functions 
and a generalized Dirichlet problem},    Michigan Math. J.,
 {\bf  25}  (1978),  299-316. 

\smallskip

   \noindent
\item{[I]}    H. Ishii,    {\sl  Perron's method for Hamilton-Jacobi equations},    Duke  Math. J.{\bf 55} (1987),  369-384.

\smallskip

   \noindent
\item{[K]}    N. V. Krylov,    {\sl  On the general notion of fully nonlinear second-order elliptic equations},    Trans. Amer. Math. Soc. (3)
 {\bf  347}  (1979), 30-34.

\smallskip

   \noindent
\item{[LM]}    H. B. Lawson, Jr and M.-L. Michelsohn,    {Spin Geometry},    Princeton University  Press, 
Princeton, NJ, 1989.

\smallskip

   \noindent
\item{[LE]}    Y. Luo and A. Eberhard,    {An application of $C^{1,1}$ approximation to
comparison principles for viscosity solutions of curvatures euqations},   Nonlinear Analysis
 {\bf 64} (2006), 1236-1254.

\smallskip

\noindent
\item{[PZ]}   S. Peng and D. Zhou, 
{\sl Maximum principle for viscosity solutions on riemannian manifolds},    
ArXiv:0806.4768, June 2008.

\smallskip

\item {[RT]} J. B. Rauch and B. A. Taylor, {\sl  The Dirichlet problem for the 
multidimensional Monge-Amp\`ere equation},
Rocky Mountain J. Math {\bf 7}    (1977), 345-364.

\smallskip

\item {[S]}  Z. Slodkowski, {\sl  The Bremermann-Dirichlet problem for $q$-plurisubharmonic functions},
Ann. Scuola Norm. Sup. Pisa Cl. Sci. (4)  {\bf 11}    (1984),  303-326.

\smallskip

\item {[W]}   J.  B. Walsh,  {\sl Continuity of envelopes of plurisubharmonic functions},
 J. Math. Mech. 
{\bf 18}  (1968-69),   143-148.

\smallskip

\end

Suppose in particular that we choose the framing 
${\partial\over \partial x_1},..., {\partial\over \partial x_n}$ for local coordinates defined on $U$.
Then from (5.2)$'$ we have that the local trivialization is 
$$
(\Phi\circ\Psi)(u)\ =\ (u, Du, D^2 u - \Gamma_x  (Du)).
\eqno{(5.5)}
$$

\medskip
\noindent
{\bf 5.4.  Transformations associated to a change of framing.}
Suppose we are given a change of framing $e = h \overline e$ over the neighborhood $U$, 
that is,
$$
e_i(x) \ =\ \sum_{j=1}^n h_{ij}(x) \overline e_j(x)
\eqno{(5.6)}
$$
for a smooth function $h:U\to {\rm GL}_n(\bbr)$. Then the 
associated trivializations are related by
$$
(r, \,p,\,H)\ =\ (\overline r\ , h(x) \cdot \overline p, \ h(x) \cdot \overline H  \cdot h^t(x))
\eqno{(5.7)}
$$
as one can see by noting for example that 
$$
(\Hess\, u)(e_i, e_j) 
 = (\Hess\, u)\left (\sum_k h_{ik}\overline e_k, \sum_{\ell} h_{j\ell}\overline e_\ell\right) =
 \sum_{k,\ell} h_{ik} (\Hess\, u)(\overline e_k, \overline e_\ell) h_{j\ell}.
 $$

\pf
For clarity we  give the proof when $F$  is a constant coefficient subequation on an open set
$X\ss\rn$.  The hypothesis is that : there exist $C_0\geq0$ $\e_0>0$ and $r_0>0$ such that
for $C\geq C_0$, $0\leq \e\leq \e_0$ and $|x-x_0|\leq r_0$, the reduced 2-jet of the function
$\b(x) = C(\rho(x) -{\e\over 2} |x-x_0|^2)$ belongs to $\Int F$.  Set 
$$
J^2_x\rho = (p(x), A(x))\and J_0 = J_{x_0}^2\rho = (p_0, A_0).
$$
Then by hypothesis
$$
C(p(x)-\e(x-x_0), A(x) - \e{\rm Id}) \ \in\ \Int F \fa C\geq C_0, 0<\e\leq \e_0, \ {\rm and\ } |x-x_0|\leq r_0.
\eqno{(11.11)}
$$

Define a neighborhood $\cn(J_0)$ by requiring that 
$$
|p-p_0|\leq \bar r \and A-A_0\geq  - \bar \e {\rm Id}
$$
where $J=(p,A)$.  We show that there are sufficient 2-jets given by (11.11) to fill out
$C\cn(J_0)$, thus proving that $C\cn(J_0)\ss \Int F$ if $C\geq C_0$.

Consider the map $\Phi(x) = p(x) -\e(x-x_0)$, and note that $\Phi(x_0)= p_0$.
The derivative $\Phi'(x_0) = A_0-\e{\rm Id}$ is non-singular for a dense set  of $\e$,
so we may choose $0<\e\leq \e_0$ so that $\Phi$ is a diffeomorphism in a neighborhood
of $x_0$.  Pick $0<\bar \e<\e$ and then pick $0<r<r_0$ so that  
\medskip

(1)
$\Phi$ restricted to $B(x_0,r_0)$ is a diffeomorphism, and\medskip

(2)
$A_0-A(x) \geq -(\e-\bar \e){\rm Id}$\ \  if $|x-x_0|   <  r$. \medskip

Now pick $\bar r>0$ so that $B(p_0,\bar r)\ss \Phi(B(x_0,r_0))$.
Then given $p\in B(p_0, \bar r)$, there exists a unique $x\in B(x_0, r_0)$ satisfying $p=p(x)-\e(x-x_0)$.
By (11.11)  this proves that 
$$
C(p, A(x) - \e{\rm Id}) \ \in\ \Int F \fa C\geq C_0.
$$
Since $(p,A) \in \cn(J_0)$ defined above, we have $A\geq A_0-\bar \e {\rm Id}$. Also by (2) we have
$A_0-\bar \e {\rm Id} \geq A(x) - \e {\rm Id}$.  The positivity condition for $\Int F$ now implies that 
$C(p,A) \in \Int F$ for all $C\geq C_0$. \qed